\documentclass[leqno,11pt]{article}  
\usepackage{amsmath}
\usepackage{amssymb}

\usepackage[german,english]{babel}
\usepackage{amsthm}
\usepackage{xypic}  

\author{\textsc{Elmar Grosse-Kl\"onne}}
\date{}

\theoremstyle{plain} 
\newtheorem{satz}{Theorem}[section]  
\newtheorem{kor}[satz]{Corollary}  
\newtheorem{lem}[satz]{Lemma}  
\newtheorem{pro}[satz]{Proposition}  
 
\newcommand{\id}{\mbox{\rm id}}

\theoremstyle{remark}

\theoremstyle{definition}

\newcommand{\w}{\omega}

\newcommand{\0}{\ensuremath{\overrightarrow{0}}}

\begin{document}

\begin{center}{\bf From pro-$p$ Iwahori-Hecke modules to
    $(\varphi,\Gamma)$-modules I}\\by Elmar Grosse-Kl\"onne
\end{center}

\begin{abstract} Let ${\mathfrak o}$ be the ring of integers in a finite
  extension $K$ of
${\mathbb Q}_p$, let $k$ be its residue field. Let $G$ be a split reductive group over ${\mathbb
  Q}_p$, let $T$ be a maximal split torus in $G$. Let ${\mathcal
  H}(G,I_0)$ be the pro-$p$-Iwahori Hecke ${\mathfrak o}$-algebra. Given a
semiinfinite reduced chamber gallery (alcove walk) $C^{({\bullet})}$ in the $T$-stable apartment, a period  $\phi\in N(T)$ of
$C^{({\bullet})}$ of length $r$ and a homomorphism $\tau:{\mathbb
  Z}_p^{\times}\to T$ compatible with $\phi$, we construct a functor from the category ${\rm
  Mod}^{\rm fin}({\mathcal
  H}(G,I_0))$ of finite length ${\mathcal
  H}(G,I_0)$-modules to \'{e}tale $(\varphi^r,\Gamma)$-modules over Fontaine's
ring ${\mathcal
  O}_{\mathcal E}$. If $G={\rm GL}_{d+1}({\mathbb Q}_p)$ there are essentially
two choices of ($C^{({\bullet})}$, $\phi$, $\tau$) with $r=1$, both leading to a functor from ${\rm
  Mod}^{\rm fin}({\mathcal
  H}(G,I_0))$ to \'{e}tale $(\varphi,\Gamma)$-modules and hence to ${\rm
  Gal}_{{\mathbb Q}_p}$-representations. Both induce a bijection between the
set of absolutely simple supersingular
${\mathcal H}(G,I_0)\otimes_{\mathfrak o} k$-modules of dimension $d+1$ and
the set of irreducible representations of ${\rm Gal}_{{\mathbb Q}_p}$
over $k$ of dimension $d+1$. We also compute these functors on modular reductions of
tamely ramified locally unitary principal series representations of $G$ over
$K$. For $d=1$ we recover Colmez' functor (when restricted to ${\mathfrak
  o}$-torsion ${\rm
  GL}_{2}({\mathbb Q}_p)$-representations generated by their pro-$p$-Iwahori invariants).

\end{abstract}

\tableofcontents

\section{Introduction}

In his remarkable opus \cite{colhaupt} on the $p$-adic local Langlands correspondence for ${\rm GL}_2({\mathbb
  Q}_p)$, Colmez established a bijection between certain representations of
${\rm GL}_2({\mathbb
  Q}_p)$ and certain two-dimensional representations of the absolute Galois
group ${\rm Gal}_{{\mathbb Q}_p}$ of the field ${\mathbb Q}_p$ of $p$-adic numbers. These representations have coefficients either in a finite
extension of ${\mathbb F}_p$, or in a finite extension of ${\mathbb Q}_p$. In
either case, the theory of $(\varphi,\Gamma)$-modules as developed by Fontaine
\cite{fon} provides the required intermediate objects in order to pass from
one side to the other. Prior to Colmez' work the characteristic $p$ correspondence
had been suggested by Breuil as an explicit "by hand" matching between the
objects on either side; it was then astonishing to see this
correspondence being realized even by a {\it functorial} relationship between ${\rm GL}_2({\mathbb
  Q}_p)$-representations and $(\varphi,\Gamma)$-modules. A certain functor ${\bf
  D}$ from ${\mathfrak o}$-torsion representations of ${\rm GL}_2({\mathbb
  Q}_p)$ to
$(\varphi,\Gamma)$-modules over ${\mathfrak o}$ constitutes one half of this
relationship. Here ${\mathfrak o}$ is the ring of integers in a finite
extension $K$ of ${\mathbb Q}_p$. Although Colmez does not phrase it in these
terms, his functor ${\bf
  D}$ may be viewed as factoring through a functor from certain
coefficient systems on the Bruhat Tits tree ${\mathfrak X}$ of ${\rm PGL}_2({\mathbb
  Q}_p)$ to $(\varphi,\Gamma)$-modules. The purpose of the present paper is to
suggest an extension of this latter functor to {\it certain} coefficient systems on
the Bruhat Tits building $X$ of a general split reductive group $G$ over ${\mathbb
  Q}_p$. Such coefficient systems can in particular be attached to
(${\mathfrak o}$-finite-length) modules
over the pro-$p$-Iwahori Hecke ${\mathfrak o}$-algebra ${\mathcal
  H}(G,I_0)$, formed with respect to a pro-$p$-Iwahori subgroup $I_0$ in $G$. The entire
construction depends on a certain choice, and for each such choice we end up (Theorem \ref{exakt}) with an exact functor
from such ${\mathcal
  H}(G,I_0)$-modules to $(\varphi^r,\Gamma)$-modules, with $r\in{\mathbb
  N}$ depending on that choice.

Let ${\mathfrak
  v}_0$ denote the vertex of ${\mathfrak X}$ fixed by ${\rm GL}_2({\mathbb
  Z}_p)$. In ${\rm GL}_2({\mathbb
  Q}_p)$ consider the element $\varphi=\left( \begin{array}{cc} p & 0 \\ 0 &
    1\end{array}\right)$ and the subgroups ${\mathfrak N}_0=\left( \begin{array}{cc} 1 & {\mathbb Z}_p \\
    0 & 1\end{array}\right)$ and $\Gamma=\left( \begin{array}{cc}
    {\mathbb Z}_p^{\times} & 0 \\
    0 & 1\end{array}\right)$. The orbit of ${\mathfrak
  v}_0$ under the submonoid $\lfloor {\mathfrak N}_0,\varphi,\Gamma\rfloor$ of $G$
generated by ${\mathfrak N}_0$, $\varphi$ and $\Gamma$ defines a halftree $\overline{\mathfrak
X}_+$ inside ${\mathfrak X}$: its edges are those whose both vertices belong to
that orbit. Adding the unique edge with only one vertex (namely ${\mathfrak
  v}_0$) in that orbit we
obtain the half tree ${\mathfrak
X}_+$. Let ${\mathcal V}$ be a $\lfloor {\mathfrak N}_0,\varphi,\Gamma\rfloor$-equivariant ${\mathfrak o}$-torsion coefficient system on ${\mathfrak
X}_+$. Let $D({\mathcal V})=H_0(\overline{\mathfrak
X}_+,{\mathcal V})^*$ and $D'({\mathcal V})=H_0({\mathfrak
X}_+,{\mathcal V})^*$ (Pontryagin duals). Under a suitable finiteness conditions these are
compact ${\mathcal O}_{\mathcal E}^+={\mathfrak o}[[{\mathfrak N}_0]]$-modules
and the natural map $D'({\mathcal V})\otimes_{{\mathcal O}_{\mathcal
    E}^+}{\mathcal O}_{\mathcal E}\to D({\mathcal V})\otimes_{{\mathcal
    O}_{\mathcal E}^+}{\mathcal O}_{\mathcal E}=:{\bf D}({\mathcal V})$ is bijective, where
${\mathcal O}_{\mathcal E}$ denotes the $p$-adic completion of ${\mathcal
    O}_{\mathcal E}^+$ with respect to the complement of $\pi_K{\mathcal
    O}_{\mathcal E}^+$, where $\pi_K\in {\mathfrak o}$ is a uniformizer. The
  actions of $\varphi$ and $\Gamma$ then provide ${\bf D}({\mathcal V})$ with the
  structure of an \'{e}tale $(\varphi,\Gamma)$-module. (The ${\mathcal
    O}_{\mathcal E}^+$-lattice $D'({\mathcal V})$ carries the $\varphi$-operator, the ${\mathcal
    O}_{\mathcal E}^+$-lattice $D({\mathcal V})$ carries the $\psi$-operator.)

Now suppose we are given a
$G$-equivariant coefficient system ${\mathcal V}$ on $X$. By what we said, in order to pass from ${\mathcal V}$ to an
\'{e}tale $(\varphi,\Gamma)$-module we might try to assign to ${\mathcal V}$ a
$\lfloor {\mathfrak N}_0,\varphi,\Gamma\rfloor$-equivariant
coefficient system on ${\mathfrak
X}_+$. Thinking sheaf theoretically, a very naive pattern would simply be:
choose an embedding (a notion to be clarified) $\iota:{\mathfrak
X}_+\to X$ and take the pull back $\iota^{-1}{\mathcal V}$ of ${\mathcal V}$ to ${\mathfrak
X}_+$. However, in order that $\iota^{-1}{\mathcal V}$ obtains a
$\lfloor{\mathfrak N}_0,\varphi,\Gamma\rfloor$-action from the $G$-action on
${\mathcal V}$ one must ask an equivariance property of $\iota$. At first one
would think of such an equivariance property with respect to a chosen
embedding of $\lfloor{\mathfrak N}_0,\varphi,\Gamma\rfloor$ into $G$, but if $G\ne {\rm GL}_2({\mathbb
  Q}_p)$ such a construction is apparently not available. Now the main point of the
present paper is that this embedding idea can nevertheless be implemented if
the equivariance property is negociated down to the bare minimum which just
suffices for such an equivariant pull back of coefficient
systems --- but only for equivariant coefficient
systems of a very specific type: coefficient
systems "of level $1$" in our terminology. 

Let $T\subset G$ be a split maximal torus, $A\subset X$ the corresponding
apartment, $C\subset A$ a chamber, $I\subset G$ the corresponding Iwahori
subgroup, $I_0\subset I$ its pro-$p$-Iwahori subgroup, $NT$ a Borel subgroup with unipotent radical $N$ and $N_0=N\cap
I$. We fix an infinite reduced chamber gallery $C^{({\bullet})}=(C=C^{(0)}, C^{(1)},
C^{(2)},\ldots)$ in $A$ with $|N_0\cdot C^{(i)}|=p^i$ for $i\ge0$. Given $C^{({\bullet})}$, we
choose an auxiliary datum, which we call $\Theta$: for all $i\ge0$ an
identification of the orbit $N_0\cdot C^{(i)}$ with
the set of edges of ${\mathfrak
X}_+$ "at distance $i$", and of the orbit $N_0\cdot ({C}^{(i)}\cap {C}^{(i+1)})$ with
the set of vertices of ${\mathfrak
X}_+$ "at distance $i$". This should respect face inclusions.\footnote{By convention, our chambers
are the {\it closures} of those open subsets usually referred to as the chambers
of $X$.} Carefully chosen, such a $\Theta$ identifies the
action of $N_0$ on the neighbourhood of any $n\cdot({C}^{(i)}\cap
{C}^{(i+1)})$ (for $n\in N_0$) --- i.e. on the chambers containing $n\cdot({C}^{(i)}\cap
{C}^{(i+1)})$ --- with the action of ${\mathfrak N}_0$ on the neighbourhood of
the corresponding vertex in ${\mathfrak
X}_+$. See Theorem \ref{baumeinbettung}. Such a $\Theta$ given, 
$N_0$-equivariant coefficient systems ${\mathcal V}$ on $X$ of level $1$ give
rise to ${\mathfrak N}_0$-equivariant coefficient systems $\Theta_*{\mathcal V}$ on ${\mathfrak
X}_+$: this is essentially built into the level $1$ property which exactly says that the
$N_0$-action on ${\mathcal V}(n\cdot({C}^{(i)}\cap
{C}^{(i+1)}))$ is insensitive to the subgroup fixing the above neighbourhood of $n\cdot({C}^{(i)}\cap
{C}^{(i+1)})$
pointwise.\footnote{Notice that we only need values of ${\mathcal V}$ on all the $n\cdot{C}^{(i)}$ and
$n\cdot({C}^{(i)}\cap{C}^{(i+1)})$, not on smaller facets. Therefore, deviating from
usual terminology, a 'coefficient system' in this paper is ony defined on
facets of codimension $0$ or $1$.} See Theorem \ref{pufo}. The ambiguity in the particular choice of $\Theta$
is immaterial for our further purposes, as long as $C^{({\bullet})}$ is fixed. 

If $\phi\in N(T)$ is a period of length $r\in{\mathbb N}$ for
$C^{({\bullet})}$, i.e. if $\phi(C^{(i)})=C^{(i+r)}$ for all $i\ge0$, then $\Theta$ can
be chosen in such a way that $\Theta_*{\mathcal V}$ carries an action of
$\lfloor{\mathfrak N}_0,\varphi^r\rfloor$. Similarly, if for an embedding $\tau:{\mathbb
Z}_p^{\times}\to T$ the image $\tau({\mathbb
Z}_p^{\times})$ commutes with $\phi$ and acts semilinearly on ${\mathcal V}$, then $\Theta_*{\mathcal V}$
carries an action of $\lfloor{\mathfrak N}_0,\varphi^r,\Gamma\rfloor$.

For equivariant coefficient systems ${\mathcal V}$ of level $1$ on $X$ satisfying a certain finiteness condition we thus obtain
an \'{e}tale $(\varphi^r,\Gamma)$-module ${\bf D}(\Theta_*{\mathcal V})$ by the
functor analogous to the one discussed at the beginning (now with a
$\varphi^r$-action instead of a $\varphi$-action).

To any admissible ${\mathfrak o}$-torsion representation $V$ of $G$ one can
associate a $G$-equivariant coefficient system ${\mathcal V}$ on
$X$ which indeed satisfies the said finiteness condition and which to a chamber $C'$ of $X$ with corresponding
pro-$p$-Iwahori subgroup $I'_0$ assigns ${\mathcal
  V}(C')=V^{I'_0}$. But there is even a functorial construction of such
coefficient systems depending only on the ${\mathcal
  H}(G,I_0)$-module $V^{I_0}$ and not on the $G$-representation $V$. In other words, there is a functor $M\mapsto{\mathcal V}_M$ from the category ${\rm
  Mod}^{\rm fin}({\mathcal
  H}(G,I_0))$ of ${\mathcal
  H}(G,I_0)$-modules which (as ${\mathfrak o}$-modules) are of finite length, to
coefficient systems on $X$ which satisfy the said finiteness condition, cf. Proposition \ref{jalgheck}. One gets a
functor $$M\mapsto{\bf D}(\Theta_*{\mathcal V}_M)$$ from ${\rm
  Mod}^{\rm fin}({\mathcal
  H}(G,I_0))$ to the category of \'{e}tale
$(\varphi^r,\Gamma)$-modules; it is {\it exact}. See Theorem \ref{exakt}.  If the ${\mathfrak o}$-action
on $M$ factors through $k$ then the action of ${\mathcal O}_{{\mathcal E}}$ on
${\bf D}(\Theta_*{\mathcal V}_M)$ factors through the residue field
$k_{\mathcal E}$($\cong k((t))$) of ${\mathcal O}_{{\mathcal E}}$ and we have$${\rm dim}_{k_{\mathcal E}}{\bf D}(\Theta_*{\mathcal V}_M)\le{\rm dim}_kM.$$The functor depends on the choice of $C^{({\bullet})}$, $\phi$ and $\tau$. It is of course of interest to find $C^{({\bullet})}$, $\phi$, $\tau$ such that the length $r$ of the period $\phi$ is small.

 Let $\Gamma_0$ be the maximal pro-$p$-subgroup of
  $\Gamma$. A slight simplification of the above construction produces an \'{e}tale $(\varphi^r,\Gamma_0)$-module ${\bf
    D}(\Theta_*{\mathcal V}_M)$ for {\it any} choice of $(C^{({\bullet})},\phi)$,
  not requiring a cocharacter $\tau$ as above. (The existence of such $\tau$ apparently is a
  fairly restrictive assumption on $(C^{({\bullet})},\phi)$.) Thus, any pair
  $(C^{({\bullet})},\phi)$ gives rise to a functor $M\mapsto{\bf
    D}(\Theta_*{\mathcal V}_M)$ to \'{e}tale $(\varphi^r,\Gamma_0)$-modules, cf. Proposition \ref{natgam}, Corollary \ref{torsch} (a), Theorem \ref{exakt} (a).

Let $G={\rm GL}_{d+1}({\mathbb Q}_p)$ for some $d\in{\mathbb N}$. Then, asking for $C^{({\bullet})}$, $\phi$, $\tau$ with $r=1$ one has essentially just two choices, namely where $C^{({\bullet})}$ is interpolated by the translates of $C$ by one of the two (if $d>1$) extreme (in the Dynkin diagram) simple coroots. As $r=1$ we obtain usual \'{e}tale $(\varphi,\Gamma)$-modules ${\bf
    D}(\Theta_*{\mathcal V}_M)$, hence we can pass to their corresponding ${\rm Gal}_{{\mathbb Q}_p}$-representations $W({\bf
    D}(\Theta_*{\mathcal V}_M))$. We compute the resulting functor \begin{gather}M\mapsto W({\bf D}(\Theta_*{\mathcal V}_M))\label{semendfr}\end{gather}from  ${\rm
  Mod}^{\rm fin}({\mathcal
  H}(G,I_0))$ to ${\rm Gal}_{{\mathbb Q}_p}$-representations over ${\mathfrak o}$ in important cases. Let ${\mathcal H}(G,I_0)_k={\mathcal H}(G,I_0)\otimes_{{\mathfrak o}}k$.\\

{\bf Theorem \ref{supersinghaupt}:} {\it If $G={\rm GL}_{d+1}({\mathbb Q}_p)$ the functor (\ref{semendfr}) induces a bijection between 

(a) the set of isomorphism classes of absolutely simple supersingular
${\mathcal H}(G,I_0)_k$-modules of dimension $d+1$ and 

(b) the set of isomorphism classes of smooth irreducible representations of ${\rm Gal}_{{\mathbb Q}_p}$ over $k$ of dimension $d+1$.}\\

The required input about ${\mathcal H}(G,I_0)_k$-modules is provided by
work of Ollivier \cite{oll} and Vign\'{e}ras \cite{vigneras} while the
facts needed about ${\rm Gal}_{{\mathbb Q}_p}$-representations we found in work of
Berger \cite{berger}. These ingredients allow us to make the correspondence
in Theorem \ref{supersinghaupt} completely explicit, similarly to Breuil's matching correspondence in the case of ${\rm GL}_{2}({\mathbb Q}_p)$. 

Still concentrating on $G={\rm GL}_{d+1}({\mathbb Q}_p)$, we next evaluate the functor on {\it reduced standard ${\mathcal H}(G,I_0)_k$-modules} (or {\it ${\mathcal H}(G,I_0)_k$-modules of $W$-type}). Such modules admit a $k$-basis indexed by the finite Weyl group $W$ in terms of which the ${\mathcal H}(G,I_0)_k$-action can be given very neatly. One may expect reduced standard modules to play a similar role for the representation theory of ${\mathcal H}(G,I_0)_k$ as standard modules do in similar and more classical contexts. If $Y$ is a locally
unitary tamely ramified principal series representation of $G$ over
$K$ then the ${\mathcal H}(G,I_0)_k$-module arising by
modular reduction from the ${\mathcal H}(G,I_0)\otimes_{{\mathfrak o}}K$-module $Y^{I_0}$ of $I_0$-invariants is a reduced standard ${\mathcal H}(G,I_0)_k$-module, see \cite{wty}. (Here, by local unitarity of $Y$ we precisely mean that $Y^{I_0}$ admits an ${\mathcal
  H}(G,I_0)$-stable ${\mathfrak o}$-lattice.) For any reduced standard ${\mathcal
  H}(G,I_0)_k$-module $M$ we explain how the ${\rm Gal}_{{\mathbb Q}_p}$-representation $W({\bf
    D}(\Theta_*{\mathcal V}_M))$ can be filtered in such a way that the
  subquotients take the form ${\rm
    ind}(\omega_{m+1}^h)\otimes\omega^s\mu_{\beta}$ with varying $m\ge0$,
$h, s\in{\mathbb Z}$, $\beta\in k^{\rm alg}$ (in usual notations,
  e.g. as in \cite{berger}). Such a filtration is induced from a suitable $k$-vector space filtration on $M$ which itself is induced from a suitable (set theoretical) filtration on $W$ associated with $M$. In the special case where
  $M=V^{I_0}$ for a tamely ramified principal series representation $V$ of $G$
  {\it over $k$} one can be even more precise. In that case we describe a filtration
  of $W({\bf D}(\Theta_*{\mathcal V}_M))$ with $d$-dimensional subquotients of the form ${\rm
    ind}(\omega_{d}^h)\otimes\omega^s\mu_{\beta}$ with varying $h, s, \beta$,
  and we find ${\rm dim}_{k}W({\bf D}(\Theta_*{\mathcal V}_M))=d!d$ (whereas ${\rm dim}_k(M)=(d+1)!$, as for any reduced standard ${\mathcal H}(G,I_0)_k$-module). See
  Theorem \ref{descmodprincser} for the precise statement. It generalizes to $G={\rm GL}_{d+1}({\mathbb Q}_p)$ Colmez' computation for $G={\rm GL}_{2}({\mathbb Q}_p)$.

The key step in computing the ${\rm Gal}_{{\mathbb Q}_p}$-representation $W({\bf
    D}(\Theta_*{\mathcal V}_M))$ is the passage from the module $H_0(\overline{\mathfrak
X}_+,\Theta_*{\mathcal V}_M)$ over the non commutative polynomial ring
$k_{\mathcal E}^+[\varphi,\Gamma]$ over $k_{\mathcal E}^+=k[[{\mathfrak
  N}_0]]$($\cong k[[t]]$) to the \'{e}tale
$(\varphi,\Gamma)$-module ${\bf
    D}(\Theta_*{\mathcal V}_M)$ over $k_{\mathcal E}$. To this end, inspired
  by section 5 of \cite{emert} we introduce a notion of
  'standard cyclic modules' over $k_{\mathcal E}^+[\varphi,\Gamma]$; these give
  rise to $(\varphi,\Gamma)$-modules over $k_{\mathcal E}$ whose associated ${\rm
    Gal}_{{\mathbb Q}_p}$-representations are of the form ${\rm
    ind}(\omega_{m+1}^h)\otimes\omega^s\mu_{\beta}$.

The structure of the paper is as follows. We start, section \ref{pmodsl2}, by recalling preliminaries about mod $p$ representations of ${\rm
  SL}_2({\mathbb F}_p)$ as needed in sections \ref{cose}, \ref{iwah},
\ref{gln}. In section \ref{coefsec} we present our main
geometric construction: an isomorphism $\Theta$ between the $N_0$-orbit of a
semiinfinite chamber gallery $C^{({\bullet})}$ in $X$ as above, and the half
tree ${\mathfrak X}_+$; this $\Theta$ can be chosen to be 'equivariant' in a
sense suitable for our purposes. In section \ref{cose} we introduce our notion of equivariant coefficient systems as needed
in this paper, and their (strict)
level $1$-property. We explain the passage from equivariant coefficient systems of
level $1$ on $X$ to such on ${\mathfrak X}_+$ by means of $\Theta$, and we
prove the important Theorem \ref{finite} which says that the ${\mathfrak
  N}_0$-invariants in $H_0(\overline{\mathfrak
X}_+,{\mathcal V})$ for an ${\mathfrak
  N}_0$-equivariant coefficient system ${\mathcal V}$ on ${\mathfrak X}_+$ of strict
level $1$ are just the obvious ones. The functor $M\mapsto{\mathcal V}_M$ from ${\rm
  Mod}^{\rm fin}({\mathcal
  H}(G,I_0))$ to
coefficient systems on $X$ is the subject of section \ref{iwah}. Section \ref{phgase} recalls basic facts on $(\varphi^r,\Gamma)$-modules,
followed by a discussion of standard cyclic modules over $k_{\mathcal
  E}^+[\varphi,\Gamma]$ and their corresponding  ${\rm Gal}_{{\mathbb
    Q}_p}$-representations. Section \ref{secdfu} explains the functors from equivariant coefficient systems
on ${\mathfrak X}_+$ to modules over $k_{\mathcal
  E}^+[\varphi^r,\Gamma]$, to $(\psi^r,\Gamma)$-modules and to
$(\varphi^r,\Gamma)$-modules. This is almost formal and only
rephrases the main construction of
Colmez' functor ${\bf D}$ in terms of the halftree ${\mathfrak X}_+$. We end
 with our computations of the functor $M\mapsto W({\bf
    D}(\Theta_*{\mathcal V}_M))$ in the case $G={\rm
  GL}_{d+1}({\mathbb Q}_p)$ and $r=1$, section \ref{gln}.\\

{\bf Remarks:} (a) Throughout, we could just as well replace $G$ by the group of $F$-rational
points of an $F$-split connected reductive group over $F$, for an arbitrary
local field $F$ with residue field ${\mathbb F}_p$ (but still using the tree
${\mathfrak X}$ of ${\rm GL}_2({\mathbb Q}_p)$ and the associated usual
$(\varphi,\Gamma)$-modules). 

(b) In the sequel \cite{ihg2} to this paper we discuss in more detail the functor ${\bf D}$ introduced here
for split reductive groups over ${\mathbb Q}_p$ other than ${\rm GL}_{d+1}({\mathbb Q}_p)$.\\

{\it Acknowlegdements:} Many thanks go to Peter Schneider who gave helpful remarks on some
technical details concerning Theorem \ref{baumeinbettung}. Moreover, he suggested the description of the coefficient system ${\mathcal V}_M^X$ as given in section \ref{iwah}, significantly simplifying my original approach. I thank the
members of the working group on $p$-adic arithmetic at the University of M\"unster for the
very careful reading of this paper and for numerous remarks. In particular, I thank Marten Bornmann, Jan Kohlhaase and Torsten Schoeneberg for pointing out inaccuracies or clumsy formulations. I thank Laurent Berger for additional comments on \cite{berger}. Thanks also got to the referee for suggesting various improvements in the presentation.

{\it Notations:} We fix a finite extension $K/{\mathbb Q}_p$ with ring of integers ${\mathfrak
  o}$, prime element $\pi_K$ and residue field $k$. We let ${\mathbb Q}_p^{\rm
  alg}$ resp. $k^{\rm
  alg}$ denote an algebraic closure of ${\mathbb Q}_p$ resp. of $k$. For $m\ge1$ we write ${\mathfrak
  o}_m={\mathfrak o}/\pi_K^m$. Recall that Pontryagin duality $V\mapsto V^*={\rm Hom}_{{\mathfrak o}}^{\rm
  ct}(V,K/{\mathfrak o})$ sets up an equivalence between the category of all
torsion ${\mathfrak o}$-modules and the category of all compact
linear-topological ${\mathfrak o}$-modules.

\section{Mod $p$ representations of ${\rm SL}_2({\mathbb F}_p)$}
In this section we recall well known facts on the representation theory of
${\rm SL}_2({\mathbb F}_p)$ and its quotient ${\rm PSL}_2({\mathbb F}_p)$ on $k$-vector spaces. We use usual conventions
concerning induced representations, Hecke algebras and Hecke operators, just
as they
are recalled explicitly in the later section \ref{iwah}. 

\label{pmodsl2}

\begin{lem}\label{einsbcjalg} Let $\overline{\mathcal U}$ be a cyclic group with $p$
  elements. Let $W$ be a $k[\overline{\mathcal U}]$-module, generated
  by a finite dimensional sub $k$-vector space $W'$ of $W$ with ${\rm dim}_k(W')={\rm dim}_k(W^{\overline{\mathcal U}})$. Let $$\eta:k[\overline{\mathcal U}]\otimes_kW'\longrightarrow W$$denote the surjective morphism of $k[\overline{\mathcal U}]$-modules induced from the inclusion $W'\hookrightarrow W$. Then the map $H^1({\mathbb Z}_p,\eta):H^1({\mathbb Z}_p,k[\overline{\mathcal U}]\otimes_kW')\to H^1({\mathbb Z}_p,W)$ (continuous cohomology) is bijective for every surjection ${\mathbb Z}_p\to\overline{\mathcal U}$.   
\end{lem}

{\sc Proof:} The same as for \cite{jalg} Lemma 2.1 (ii) and (iii). For the convenience of the reader, we reproduce the proof.

(i) Any (finitely generated) $k[\overline{\mathcal U}]$-module admits a direct sum decomposition with summands isomorphic to quotients of $k[\overline{\mathcal U}]$: this can be seen  e.g. by applying the structure theorem for modules over the polynomial ring in one variable over $k$, of which $k[\overline{\mathcal U}]$ is a quotient. From this we see that the minimal number of elements needed to generate such a $k[\overline{\mathcal U}]$-module is the same as the dimension of its space of $\overline{\mathcal U}$-invariants. As ${\rm dim}_k(W^{{\overline{\mathcal U}}})=\dim_k(W')$ we therefore see that, as a $k[\overline{\mathcal U}]$-module, $W$ can not be generated by fewer than $\dim_k(W')$ many elements.

(ii) Let $$\epsilon: k[\overline{\mathcal U}]\otimes_k W'\longrightarrow
  W',\quad\quad \bar{u}\otimes w\mapsto w\mbox{ for }\bar{u}\in \overline{\mathcal U}$$ denote the
  augmentation map. We claim that ${\rm ker}(\eta)\subset{\rm ker}(\epsilon)$.

To see this let $x\in\ker(\eta)$. If $x\notin \ker(\epsilon)$ then the class of $x$ in $(k[\overline{\mathcal U}]\otimes_k W')\otimes_{k[\overline{\mathcal U}]}k$ does not vanish (here $k[\overline{\mathcal U}]\to k$ is the augmentation, its kernel is the maximal ideal of the local ring $k[\overline{\mathcal U}]$). Therefore, by Nakayama's Lemma, the quotient $(k[\overline{\mathcal U}]\otimes_k W')/(k[\overline{\mathcal U}].x)$ can be generated by fewer than $\dim_{k}(W')$ many elements. As $x\in \ker(\eta)$ this is a contradiction to what we saw in (i).

(iii) Let $\gamma:{\mathbb Z}_p\to k[\overline{\mathcal U}]\otimes_k
W'$ be a $1$-cocycle such that $\eta\circ\gamma:{\mathbb
  Z}_p\to W$ is a coboundary. As $\eta$ is surjective we may modify $\gamma$
by a coboundary such that now $\eta\circ\gamma=0$. Let $c\in {\mathbb Z}_p$ be
a topological generator. Since $\eta(\gamma(c))=0$ we have $\gamma(c)\in{\rm
  ker}(\epsilon)$ by (ii). But ${\rm
  ker}(\epsilon)=(c-1)k[\overline{\mathcal U}]\otimes_k W'$, so
$\gamma(c)=cf-f$ for some $f\in k[\overline{\mathcal U}]\otimes_k
W'$. Since $c$ generates ${\mathbb Z}_p$ the cocycle
condition on $\gamma$ shows $\gamma(c')=c'f-f$ for any $c'\in {\mathbb Z}_p$,
so $\gamma$ is a coboundary. We have shown injectivity of $H^1({\mathbb
  Z}_p,\eta)$. The surjectivity of $H^1({\mathbb Z}_p,\eta)$ follows from the
surjectivity of $\eta$ and from $H^2({\mathbb Z}_p,?)=0$.\hfill$\Box$\\  

Let either $\overline{\mathcal S}={\rm SL}_2({\mathbb F}_p)$ or $\overline{\mathcal S}={\rm PSL}_2({\mathbb F}_p)$. Let $\overline{\mathcal
  U}$ be the (group of ${\mathbb F}_p$-valued point of the) unipotent radical
of a Borel subgroup in $\overline{\mathcal S}$. For $m\ge1$ we consider the Hecke algebra ${\mathcal
  H}(\overline{\mathcal S},\overline{\mathcal U})_{{\mathfrak
  o}_{m}}={\rm End}_{{{\mathfrak
  o}_{m}}[\overline{\mathcal S}]}({\rm
  ind}_{\overline{\mathcal U}}^{\overline{\mathcal S}}{\bf 1}_{{\mathfrak
  o}_{m}})^{\rm op}$.

\begin{lem}\label{sl2flat} The universal module ${\rm
  ind}_{\overline{\mathcal U}}^{\overline{\mathcal S}}{\bf 1}_{{\mathfrak
  o}_{m}}$ is
  flat over ${\mathcal
  H}(\overline{\mathcal S},\overline{\mathcal U})_{{\mathfrak o}_{m}}$.
\end{lem}

{\sc Proof:} For the course of this proof let us write ${\mathcal H}_{{\mathfrak
  o}_{m}}={\mathcal
  H}(\overline{\mathcal S},\overline{\mathcal U})_{{\mathfrak
  o}_{m}}$ and more specifically ${\mathcal H}_k={\mathcal H}_{{\mathfrak
  o}_{1}}$. The flatness assertion is equivalent with the claim that for any left ideal ${\mathcal I}$ in ${\mathcal H}_{{\mathfrak
  o}_{m}}$ the natural map \begin{gather}{\rm ind}_{\overline{\mathcal U}}^{\overline{\mathcal S}}{\bf
  1}_{{\mathfrak
  o}_{m}}\otimes_{{\mathcal H}_{{\mathfrak
  o}_{m}}}{\mathcal I}\longrightarrow{\rm ind}_{\overline{\mathcal
    U}}^{\overline{\mathcal S}}{\bf
  1}_{{\mathfrak
  o}_{m}}\label{zielinj}\end{gather} is injective. We proceed by induction on $m$. For $m=1$ this is a variant of the proof given for ${\rm
  GL}_2({\mathbb F}_p)$ in Prop. 2.2 of \cite{ollsec}. For the facts on ${\mathcal
  H}_k$ stated below see e.g. \cite{calu}. Let $\overline{\mathcal T}$ be a maximal split torus in $\overline{\mathcal S}$ such that $\overline{\mathcal U}$ is the unipotent radical of the Borel subgroup $\overline{\mathcal B}=\overline{\mathcal U}\overline{\mathcal T}=\overline{\mathcal T}\overline{\mathcal U}$. We have the $\overline{\mathcal S}$-equivariant decomposition$${\rm ind}_{\overline{\mathcal U}}^{\overline{\mathcal S}}{\bf
  1}_{k}\cong\bigoplus_{\beta\in \overline{\mathcal T}^c}{\rm ind}_{\overline{\mathcal B}}^{\overline{\mathcal S}}\beta$$where $\beta$ runs through the set $\overline{\mathcal T}^c$ of $k^{\times}$-valued characters
$\beta$ of $\overline{\mathcal T}$, viewed
as characters of $\overline{\mathcal B}$. For $\beta\in\overline{\mathcal T}^c$ let
$\epsilon_{\beta}\in {\mathcal
  H}_k$ denote the natural projection ${\rm ind}_{\overline{\mathcal U}}^{\overline{\mathcal S}}{\bf
  1}_{k}\to{\rm ind}_{\overline{\mathcal B}}^{\overline{\mathcal S}}\beta$ composed with the inclusion ${\rm ind}_{\overline{\mathcal B}}^{\overline{\mathcal S}}\beta\to{\rm ind}_{\overline{\mathcal U}}^{\overline{\mathcal S}}{\bf
  1}_{k}$ . The $\epsilon_{\beta}$ are pairwise orthogonal idempotents
summing up to the unity element in ${\mathcal H}_k$. Let $T_{n_s}\in {\mathcal
  H}_k$ denote the Hecke operator corresponding to a generator $n_s\in N(\overline{\mathcal T})$ of the Weyl group
$N(\overline{\mathcal T})/\overline{\mathcal T}$. Then
${\mathcal H}_k$ is generated by $T_{n_s}$ together with all the $\epsilon_{\beta}$. For
$\beta\in\overline{\mathcal T}^c$ define $\beta^s\in\overline{\mathcal T}^c$ by $\beta^s(t)=\beta(n_stn_s^{-1})$ for $t\in \overline{\mathcal T}$. We have $T_{n_s}\epsilon_{\beta}=\epsilon_{{\beta}^s}T_{n_s}$.

In view of the above orthogonal decomposition, to
prove injectivity of (\ref{zielinj}) for $m=1$ it is
enough to show the following:

(a) For any $\beta\in\overline{\mathcal T}^c$ with $\beta\ne\beta^s$ and any left ideal
${\mathcal I}$ in ${\mathcal H}_k\epsilon_{\beta}\oplus{\mathcal
  H}_k\epsilon_{{\beta}^s}$ the map $$({\rm ind}_{\overline{\mathcal
      B}}^{\overline{\mathcal S}}{\beta}\oplus {\rm ind}_{\overline{\mathcal B}}^{\overline{\mathcal S}}{\beta}^s)\otimes_{{\mathcal H}_k\epsilon_{\beta}\oplus{\mathcal
  H}_k\epsilon_{{\beta}^s}}{\mathcal I}\longrightarrow{\rm ind}_{\overline{\mathcal B}}^{\overline{\mathcal S}}{\beta}\oplus {\rm ind}_{\overline{\mathcal B}}^{\overline{\mathcal S}}{\beta}^s$$ is injective.

(b) For (the unique) $\beta\in\overline{\mathcal T}^c$ with $\beta=\beta^s$ and any left ideal
${\mathcal I}$ in ${\mathcal H}_k\epsilon_{\beta}$ the map $$({\rm ind}_{\overline{\mathcal B}}^{\overline{\mathcal S}}{\beta})\otimes_{{\mathcal H}_k\epsilon_{\beta}}{\mathcal I}\longrightarrow{\rm ind}_{\overline{\mathcal B}}^{\overline{\mathcal S}}{\beta}$$is injective.

In (a) the image of ${T}_{n_s}$ in
${\mathcal H}_k\epsilon_{\beta}\oplus{\mathcal
  H}_k\epsilon_{{\beta}^s}$ (which we again denote by ${T}_{n_s}$)
satisfies ${T}^2_{n_s}=0$, and besides the relations already stated, there are
no further ones in ${\mathcal H}_k\epsilon_{\beta}\oplus{\mathcal
  H}_k\epsilon_{{\beta}^s}$. As
$\epsilon_{\beta}$ and $\epsilon_{{\beta}^s}$ are orthogonal idempotents, we
may further split up the situation, and the only critical case
to be considered is where ${\mathcal I}={\mathcal
  H}_k\epsilon_{{\beta}^s}{T}_{n_s}$ (or symmetrically ${\mathcal I}={\mathcal
  H}_k\epsilon_{{\beta}}{T}_{n_s}$). We may rewrite the inclusion ${\mathcal I}={\mathcal
  H}_k\epsilon_{{\beta}^s}{T}_{n_s}\to {\mathcal
  H}_k\epsilon_{{\beta}}$ as the exact sequence$${\mathcal H}_k\epsilon_{{\beta}}\stackrel{{T}_{n_s}}{\longrightarrow}{\mathcal H}_k\epsilon_{{\beta}^s}\stackrel{{T}_{n_s}}{\longrightarrow}{\mathcal
  H}_k\epsilon_{{\beta}}$$and need to verify that the resulting sequence$${\rm ind}_{\overline{\mathcal B}}^{\overline{\mathcal S}}{\beta}\stackrel{{T}_{n_s}}{\longrightarrow}{\rm ind}_{\overline{\mathcal B}}^{\overline{\mathcal S}}{\beta}^s\stackrel{{T}_{n_s}}{\longrightarrow}{\rm ind}_{\overline{\mathcal B}}^{\overline{\mathcal S}}{\beta}$$is exact. But this is well known. (This exactness is specific
to our working with the prime field ${\mathbb F}_p$; for non-prime finite fields it fails.) In (b) the argument is similar (but easier,
and valid for any finite field); namely, it boils down to the exactness of$${\rm
  ind}_{\overline{\mathcal B}}^{\overline{\mathcal
    S}}{\beta}\stackrel{{T}_{n_s}}{\longrightarrow}{\rm
  ind}_{\overline{\mathcal B}}^{\overline{\mathcal
    S}}{\beta}\stackrel{{T}_{n_s}+1}{\longrightarrow}{\rm
  ind}_{\overline{\mathcal B}}^{\overline{\mathcal
    S}}{\beta}\stackrel{{T}_{n_s}}{\longrightarrow}{\rm
  ind}_{\overline{\mathcal B}}^{\overline{\mathcal S}}{\beta}$$(for the
trivial $\beta\in\overline{\mathcal T}^c$).

Now let $m>1$. We apply ${\rm ind}_{\overline{\mathcal U}}^{\overline{\mathcal S}}{\bf
  1}_{{\mathfrak
  o}_{m}}\otimes_{{\mathcal H}_{{\mathfrak
  o}_{m}}}(.)$ to the commutative diagram with exact rows$$\xymatrix{0\ar[r]& {\mathcal I}\cap\pi_K^{m-1}{\mathcal H}_{{\mathfrak
  o}_{m}} \ar[d]\ar[r]& {\mathcal I} \ar[d]\ar[r]& {\mathcal I}/{\mathcal I}\cap\pi_K^{m-1}{\mathcal H}_{{\mathfrak
  o}_{m}} \ar[d]\ar[r]&0\\0\ar[r]& \pi_K^{m-1}{\mathcal H}_{{\mathfrak
  o}_{m}}  \ar[r]&{\mathcal H}_{{\mathfrak
  o}_{m}}\ar[r]&{\mathcal H}_{{\mathfrak
  o}_{m-1}} 
       \ar[r]&0}.$$Observe that ${\mathcal H}_{{\mathfrak
  o}_{1}}\cong\pi_K^{m-1}{\mathcal H}_{{\mathfrak
  o}_{m}}$. The bottom row then becomes the exact sequence $0\to{\rm ind}_{\overline{\mathcal U}}^{\overline{\mathcal S}}{\bf
  1}_{k}\to{\rm ind}_{\overline{\mathcal U}}^{\overline{\mathcal S}}{\bf
  1}_{{\mathfrak
  o}_{m}}\to{\rm ind}_{\overline{\mathcal U}}^{\overline{\mathcal S}}{\bf
  1}_{{\mathfrak
  o}_{m-1}}\to0$. The top row remains exact in the middle. By induction
hypothesis, the outer vertical arrows remain injective. Therefore the middle vertical arrow remains injective.\hfill$\Box$\\  

Let $\chi_{\overline{\mathcal U}}$ denote the unique element in ${\rm
  ind}_{\overline{\mathcal U}}^{\overline{\mathcal S}}{\bf 1}_{{\mathfrak
  o}_{m}}$ supported on
$\overline{\mathcal U}$ and taking constant value $1\in{{\mathfrak
  o}_{m}}$ there. For a ${\mathcal
  H}(\overline{\mathcal S},\overline{\mathcal U})_{{\mathfrak
  o}_{m}}$-(left)module $M$ consider the
natural map\begin{gather}M\longrightarrow {\rm
  ind}_{\overline{\mathcal U}}^{\overline{\mathcal S}}{\bf 1}_{{\mathfrak
  o}_{m}}\otimes_{{\mathcal
  H}(\overline{\mathcal S},\overline{\mathcal U})_{{\mathfrak
  o}_{m}}}M,\quad m\mapsto\chi_{\overline{\mathcal U}}\otimes m.\label{hecgleichinv}\end{gather}

\begin{lem}\label{sl2einbett} The map (\ref{hecgleichinv}) is an isomorphism
  from $M$ onto $({\rm
  ind}_{\overline{\mathcal U}}^{\overline{\mathcal S}}{\bf 1}_{{\mathfrak
  o}_{m}}\otimes_{{\mathcal
  H}(\overline{\mathcal S},\overline{\mathcal U})_{{\mathfrak
  o}_{m}}}M)^{\overline{\mathcal U}}$.
\end{lem}

{\sc Proof:} The same as in the
proof of \cite{jalg} Proposition 4.1 (where ${\rm
  GL}_2({\mathbb F}_p)$ instead of ${\rm
  SL}_2({\mathbb F}_p)$ or ${\rm PSL}_2({\mathbb F}_p)$ was considered). Specifically, using Lemma
\ref{sl2flat} this is reduced to the case where $M$ is an irreducible ${\mathcal
  H}(\overline{\mathcal S},\overline{\mathcal U})_{k}$-module.\hfill$\Box$\\  

\begin{lem}\label{einsajalg} Let $\overline{\mathcal U}\ne \overline{\mathcal U}'$ be the unipotent radicals of two opposite Borel subgroups in $\overline{\mathcal S}$. Let $W$ be a $k[\overline{\mathcal S}]$-module which is generated by $W^{\overline{\mathcal U}'}$. Then $W$ is generated by $W^{\overline{\mathcal U}'}$ even as a $k[\overline{\mathcal U}]$-module. 
\end{lem}

{\sc Proof:} The same as in \cite{jalg} Lemma 2.1 (i). (We remark that the analogous
statement is true more generally for ${\rm SL}_2({\mathbb F}_q)$ or ${\rm PSL}_2({\mathbb F}_q)$, for any finite field ${\mathbb F}_q$.) \hfill$\Box$\\  

For concreteness, let now more specifically $\overline{\mathcal U}$ denote the subgroup of unipotent upper triangular matrices
in ${\rm
  SL}_2({\mathbb F}_p)$. We
define $$\nu=\left(\begin{array}{cc}1&1\\0&1\end{array}\right),\quad\quad\quad
n_s=\left(\begin{array}{cc}0&1\\-1&0\end{array}\right),\quad\quad\quad h_s(x)=\left(\begin{array}{cc}x&0\\0&x^{-1}\end{array}\right)$$for
$x\in{\mathbb F}_p^{\times}$. In case $\overline{\mathcal S}={\rm
  PSL}_2({\mathbb F}_p)$ we use the same symbols $\overline{\mathcal U}$, $\nu$, $n_s$, $h_s(x)$ for the respective images in $\overline{\mathcal S}$. Let $t=[\nu]-1$ in $k[\overline{\mathcal U}]$;
this is a generator of the
maximal ideal in the local ring $k[\overline{\mathcal U}]\cong k[t]/(t^p)$. The $k$-algebra ${\mathcal
  H}(\overline{\mathcal S},\overline{\mathcal U})_k$ is generated by the Hecke
operators $T_{h_s(x)}$ for $x\in{\mathbb F}_p^{\times}$ together with
$T_{n_s}$. 

Let $0\le r\le p-1$; in case $\overline{\mathcal S}={\rm
  PSL}_2({\mathbb F}_p)$ suppose in addition that $r$ is even. We define a character $\chi_r:{\mathcal H}(\overline{\mathcal S},\overline{\mathcal
  U})_k\to k$ by requiring\begin{gather}\chi_r(T_{{h_s}(x)})=x^{-r}\quad\quad\mbox{
    for all }x\in{\mathbb F}_p^{\times}.\label{kchd}\end{gather}as well as $\chi_r(T_{n_s})=-1$ if $r=p-1$, and $\chi_r(T_{n_s})=0$ if $0\le r< p-1$.  

 We denote the one-dimensional $k$-vector space underlying $\chi_r$ again by
 $\chi_r$, and we let $e$ denote a basis element of it. By Lemma
\ref{sl2einbett} we may regard $e$ as an element of $({\rm
  ind}_{\overline{\mathcal U}}^{\overline{\mathcal S}}{\bf 1}_k)\otimes_{{\mathcal H}(\overline{\mathcal
    S},\overline{\mathcal U})_k}\chi_r$.   

\begin{lem}\label{symcla} The subspace $k.e=\chi_r$
  of $({\rm
  ind}_{\overline{\mathcal U}}^{\overline{\mathcal S}}{\bf 1}_k)\otimes_{{\mathcal H}(\overline{\mathcal
    S},\overline{\mathcal U})_k}\chi_r$ is preserved by the element
  $t^{r}n_s^{-1}$ of $k[\overline{\mathcal S}]$. We have $t^{r}n_s^{-1}e=r!e$.
\end{lem}

{\sc Proof:} It is well known that $({\rm
  ind}_{\overline{\mathcal U}}^{\overline{\mathcal S}}{\bf 1}_k)\otimes_{{\mathcal H}(\overline{\mathcal
    S},\overline{\mathcal U})_k}\chi_r$ is isomorphic with ${\rm
    Sym}^{r}(k^2)$ as a $k[\overline{\mathcal S}]$-module. By Lemma \ref{sl2einbett} the space
  of $\overline{\mathcal U}$-invariants in
  $({\rm
  ind}_{\overline{\mathcal U}}^{\overline{\mathcal S}}{\bf 1}_k)\otimes_{{\mathcal H}(\overline{\mathcal
    S},\overline{\mathcal U})_k}\chi_r$ is $\chi_r$. Therefore the claims follow by a
  straightforward computation in ${\rm
    Sym}^{r}(k^2)$.\hfill$\Box$\\

Let $1\le r\le p-1$; in case $\overline{\mathcal S}={\rm
  PSL}_2({\mathbb F}_p)$ suppose in addition that $r$ is even. We define an ${\mathcal H}(\overline{\mathcal S},\overline{\mathcal
  U})_k$-module $M_r$ with $k$-basis $e$, $f$ by requiring
$$\quad\quad  T_{n_s}(e)=f\quad\quad\quad\mbox{ and
}\quad\quad\quad\quad T_{n_s}(f)=\left\{\begin{array}{l@{\quad:\quad}l}-f&\quad
      \mbox{  }r=p-1\\0&\quad \mbox{
}r<p-1\end{array}\right.,$$$$ T_{h_s(x)}(e)=x^re\quad\quad\quad\mbox{ and
}\quad\quad\quad T_{h_s(x)}(f)=x^{-r}f\quad\quad\quad\quad\quad\quad\quad\quad\quad$$ for
all $x\in{\mathbb F}_p^{\times}$. By Lemma \ref{sl2einbett} we may regard $M_r$ as a subspace of $({\rm ind}^{\overline{\mathcal
      S}}_{\overline{\mathcal U}}{\bf 1}_k)\otimes_{{\mathcal H}(\overline{\mathcal S},\overline{\mathcal
  U})_k}M_r$.

\begin{lem}\label{sl2fpind} In $({\rm ind}^{\overline{\mathcal
      S}}_{\overline{\mathcal U}}{\bf 1}_k)\otimes_{{\mathcal H}(\overline{\mathcal S},\overline{\mathcal
  U})_k}M_r$ we
have \begin{align}
  t^{p-1}n_s^{-1}e&=f,\label{homcaa}&{}\\t^rn_s^{-1}f&=r!f,\label{homcab}&{}\\n_sf&-e&\in\sum_{i\ge0}k.t^in_se,\label{homcac}\\t^{p-1-r}n_s^{-1}e&-(p-1-r)!e&\in
  \sum_{i\ge0}k.t^in_sf.\label{homcad}\end{align}If $r=p-1$ we have
\begin{align}n_s^{-1}(f+e)&=f+e,&{}\label{homcae}\\t^{p-1}n_s^{-1}e&+e&\in\sum_{i\ge0}k.t^in_s(f+e).\label{homcaf}\end{align}
\end{lem}

{\sc Proof:} In ${\rm ind}^{\overline{\mathcal
      S}}_{\overline{\mathcal U}}{\bf 1}_k$ we
  compute$$t^{p-1}n_s^{-1}\chi_{\overline{\mathcal
      U}}=([\nu]-1)^{p-1}n_s^{-1}\chi_{\overline{\mathcal
      U}}=\sum_{i=0}^{p-1}[\nu]^in_s^{-1}\chi_{\overline{\mathcal
      U}}=T_{n_s}\chi_{\overline{\mathcal U}},$$for the last equality observe
  $\overline{\mathcal U}=\{[\nu]^i\,;\,0\le i\le p-1\}$. This proves formula (\ref{homcaa}). Let $\overline{\mathcal T}=\{h_s(x)\,|\,x\in{\mathbb F}_p^{\times}\}$, the torus of diagonal matrices
in $\overline{\mathcal S}$. Let $\theta:\overline{\mathcal T}\to k^{\times}$ be the
character defined by $\theta(h_s(x)^{-1})=x^r$. Then we have an isomorphism of $k[\overline{\mathcal S}]$-modules$$({\rm ind}^{\overline{\mathcal
      S}}_{\overline{\mathcal U}}{\bf 1}_k)\otimes_{{\mathcal H}(\overline{\mathcal S},\overline{\mathcal
  U})_k}M_r\cong{\rm ind}^{\overline{\mathcal
      S}}_{\overline{\mathcal T}\overline{\mathcal U}}\theta$$sending $e$ (resp. $f$) to the element ${\mathfrak e}$
  (resp. ${\mathfrak f}$) of $({\rm ind}^{\overline{\mathcal
      S}}_{\overline{\mathcal T}\overline{\mathcal U}}\theta)^{\overline{\mathcal U}}$ supported on
  $\overline{\mathcal T}\overline{\mathcal U}$ (resp. supported on
  $\overline{\mathcal S}-\overline{\mathcal T}\overline{\mathcal U}=\overline{\mathcal T}\overline{\mathcal U}n_s\overline{\mathcal U}$) and with ${\mathfrak e}(1)=1$ (resp. with
  ${\mathfrak f}(n_s)=1$). Indeed, in $({\rm ind}^{\overline{\mathcal
      S}}_{\overline{\mathcal U}}{\bf 1}_k)\otimes_{{\mathcal H}(\overline{\mathcal S},\overline{\mathcal
  U})_k}M_r$ we find $$x^r{e}=T_{h_s(x)}({e})=h_s(x)^{-1}{e},$$$$x^r{f}=T_{h_s(x)^{-1}}({f})=h_s(x){f},$$whereas in ${\rm ind}^{\overline{\mathcal
      S}}_{\overline{\mathcal T}\overline{\mathcal U}}\theta$ we find$$\theta(h_s(x)^{-1}){\mathfrak e}(.)={\mathfrak
  e}(h_s(x)^{-1}.)={\mathfrak
  e}(.h_s(x)^{-1})=(h_s(x)^{-1}{\mathfrak
  e})(.),$$$$\theta(h_s(x)^{-1}){\mathfrak f}(.)={\mathfrak
  f}(h_s(x)^{-1}.)={\mathfrak f}(.h_s(x))=(h_s(x){\mathfrak f})(.),$$and to
see that the scaling factor is correct we compare formula (\ref{homcaa}) with $(t^{p-1}n_s^{-1}{\mathfrak e})(n_s)={\mathfrak
  e}(1)=1$. Now $\{{\mathfrak e}\}\cup\{t^in_s{\mathfrak e}\,|\,0\le i\le p-1\}$ is a
$k$-basis of ${\rm ind}^{\overline{\mathcal
      S}}_{\overline{\mathcal T}\overline{\mathcal U}}\theta$. Since we have
  $(n_s{\mathfrak f})(1)={\mathfrak f}(n_s)=1={\mathfrak e}(1)$ but
  $(t^in_s{\mathfrak e})(1)=0$ for all $i\ge0$, formula (\ref{homcac}) follows. Next, $f$ generates a $k[\overline{\mathcal
      S}]$-submodule isomorphic with ${\rm Sym}^rk^2$: we get formula (\ref{homcab}). Finally, this $k[\overline{\mathcal
      S}]$-submodule is in fact generated by $n_sf$ even as a $k[\overline{\mathcal
      U}]$-module (Lemma \ref{einsajalg}), i.e. it coincides with $\sum_{i\ge0}k.t^in_sf$, and the quotient of $({\rm ind}^{\overline{\mathcal
      S}}_{\overline{\mathcal U}}{\bf 1}_k)\otimes_{{\mathcal H}(\overline{\mathcal S},\overline{\mathcal
  U})_k}M_r$ by this submodule is isomorphic with ${\rm Sym}^{p-1-r}k^2$ as a $k[\overline{\mathcal
      S}]$-module. This shows formula (\ref{homcad}). For the stated formulae in
    case $r=p-1$ observe that $f+e$ generates the trivial
    one dimensional $k[\overline{\mathcal
      S}]$-submodule, and dividing it out we are left with a copy of ${\rm Sym}^{p-1}k^2$.\hfill$\Box$\\

\section{Half trees in Bruhat Tits buildings}

\label{coefsec}

Let $G$ be the group of ${\mathbb Q}_p$-rational points of a ${\mathbb
  Q}_p$-split connected reductive group over ${\mathbb Q}_p$. Let $Z$ denote the center of $G$. Fix a maximal
split torus $T$ in $G$, let $N(T)$ be its normalizer in $G$. Let $\Phi$ denote the set of roots of $T$. For
$\alpha\in\Phi$ let $N_{\alpha}$ be the corresponding root subgroup in $G$. Choose a positive system $\Phi^+$ in $\Phi$. Let $N=\prod_{\alpha\in\Phi^+}N_{\alpha}$. 

Let $X$ denote the semi simple Bruhat-Tits building of $G$, let $A$ denote
 its apartment corresponding to $T$. Our notational and terminological convention is that the facets of $A$ or
 $X$ are {\it closed} in $X$ (i.e. {\it contain} all their faces (the lower dimensional facets
 at their boundary)). A chamber is a facet of codimension $0$. For a chamber $C$ in $A$ let $I_C$ be the Iwahori subgroup in $G$ fixing $C$.
Suppose we are given a semiinfinite chamber gallery\begin{gather}C^{(0)},
  C^{(1)},C^{(2)},C^{(3)},\ldots\label{cgall}\end{gather}in $A$ such that, setting $$N_0^{(i)}=I_{C^{(i)}}\cap N=\prod_{\alpha\in\Phi^+}I_{C^{(i)}}\cap N_{\alpha},$$we have\begin{gather}N_0^{(0)}\supset N_0^{(1)}\supset N_0^{(2)}\supset N_0^{(3)}\supset\ldots\quad\mbox{ with }[N_0^{(i)}:N_0^{(i+1)}]=p\mbox{ for all }i\ge0.\label{nseq}\end{gather}We write $$N_0=N_0^{(0)},\quad \quad C=C^{(0)},\quad\quad I=I_C=I_{C^{(0)}}.$$There is a unique sequence $\alpha^{(0)},\alpha^{(1)},\alpha^{(2)},\alpha^{(3)},\ldots$ in $\Phi^+$ such that, setting$$e[i,{\alpha}]=|\{0\le j\le i-1\,|\,\alpha=\alpha^{(j)}\}|$$for $i\ge0$ and $\alpha\in\Phi^+$, we have\begin{gather}N_0^{(i)}=\prod_{\alpha\in\Phi^+}(N_0\cap N_{\alpha})^{p^{e[i,\alpha]}}.\label{n0ineu}\end{gather}Geometrically, $C^{(i+1)}$ and $C^{(i)}$ share a common facet of codimension $1$ contained in a wall which belongs to the translation class of walls corresponding to $\alpha^{(i)}$.\\  

\begin{lem}\label{petlem} $N_0^{(i+1)}$ is a normal subgroup in $N_0^{(i)}$, for any $i\ge0$.
\end{lem}

{\sc Proof:} We claim that, more generally, for any inclusion 
$U_1\subset U_2$ of open subgroups of $N$ with $[U_2:U_1]=p$ and with
$U_i=\prod_{\alpha\in\Phi^+}(U_i\cap N_{\alpha})$ for $i=1,2$ we have: $U_1$
is normal in $U_2$. Indeed, there is a unique $\alpha\in \Phi^+$ with $U_2\cap N_{\beta}=U_1\cap
N_{\beta}$ for all $\beta\ne\alpha$. By standard facts on reductive groups, for $u\in N_{\alpha}$ and $v\in
N_{\gamma}$ (any $\gamma$) the commutator $uvu^{-1}v^{-1}$ is a product
of elements of the $N_{\beta}$ with $\beta\ne\alpha$. The claim follows. \hfill$\Box$\\

{\bf Definition:} We define the half tree $Y$, endowed with an action by the
group $N_0$, as follows. Both its set $Y^0$ of vertices and its set $Y^1$ of edges
are identified with the set of chambers of $X$ of the form $n\cdot C^{(i)}$ for $n\in N_0$ and $i\ge0$. The $N_0$-action is the obvious one. Whenever we view the chamber $C^{(i)}$ of $X$ as a vertex of $Y$ we
denote it by ${\bf v}_{i}$; whenever we view $C^{(i)}$ as an edge of $Y$ we
denote it by ${\bf e}_{i}$. The simplicial structure on $Y$ is given by declaring that for $n\in N_0$ and $i\ge0$ the edge $n\cdot{\bf e}_{i+1}$ contains
the two vertices $n\cdot{\bf v}_{i}$ and $n\cdot{\bf v}_{i+1}$. An element in
$\overline{Y}^1=Y^1-\{{\bf e}_{0}\}$ we also write as the unordered pair of
vertices it contains, e.g. $n\cdot{\bf e}_{i+1}=\{n\cdot{\bf
  v}_{i},n\cdot{\bf v}_{i+1}\}$. In addition we have the edge ${\bf e}_{0}$ of $Y$: it
  contains only the vertex ${\bf v}_{0}$ (hence is a 'loose' end of $Y$) and
  is fixed by $N_0$.

The half tree $\overline{Y}$ is obtained from $Y$ by removing the edge ${\bf e}_{0}$. Thus, its set of vertices is $\overline{Y}^0=Y^0$, its set of edges is $\overline{Y}^1=Y^1-\{{\bf e}_{0}\}$, and the simplicial structure is the one induced from $Y$.\\

{\bf Remark:} By definition, we have a natural bijection $Y^1\cong
  Y^0$: to ${\bf e}_{i}$ it assigns ${\bf v}_{i}$, and more
  generally, to any edge its 'outward pointing' vertex. Thus, instead of
  identifying $Y^0$ with the set of chambers $n\cdot C^{(i)}$ (with $n\in N_0$ and $i\ge0$), one might just
  as well (or perhaps more
  appropriately) identify $Y^0$ with the set of facets of codimension $1$ in
  $X$ shared by two chambers of the form $n\cdot C^{(i)}$.\\

For ${\bf v}\in Y^0$ let \begin{gather}]{\bf v}[\quad=\quad\{{\bf w}\in Y^0\,|\,\{{\bf v},{\bf
  w}\}\in Y^1\mbox{ and }{\bf v}\in [{\bf w},{\bf v}_0]\}\label{deftub}\end{gather}where $[{\bf
  w},{\bf v}_0]\subset Y^0$ denotes the set of vertices the unique geodesic
from ${\bf w}$ to ${\bf v}_0$ in $Y$ is passing through. \\

Let ${\mathfrak X}$ denote the Bruhat-Tits tree of ${\rm PGL}_2({\mathbb
  Q}_p)$. Let ${\mathfrak X}^0$ resp. ${\mathfrak X}^1$ denote the set of
vertices, resp. the set of edges of ${\mathfrak X}$. Let ${\mathfrak
  v}_0\in{\mathfrak X}^0$ denote the vertex fixed by ${\rm PGL}_2({\mathbb
  Z}_p)$. In ${\rm GL}_2({\mathbb
  Q}_p)$ we define the elements$$\varphi=\left( \begin{array}{cc} p & 0 \\ 0 & 1\end{array}\right),\quad\quad\quad \gamma_{{{0}}}=\left( \begin{array}{cc} 1+p & 0 \\ 0 & 1\end{array}\right),\quad\quad\quad \nu=\left( \begin{array}{cc} 1 & 1 \\ 0 & 1\end{array}\right)$$

and the subgroups$$\Gamma=\left( \begin{array}{cc}
    {\mathbb Z}_p^{\times} & 0 \\
    0 & 1\end{array}\right),\quad\quad\quad\Gamma_{{{0}}}=\left( \begin{array}{cc}
    1+p{\mathbb Z}_p & 0 \\
    0 & 1\end{array}\right),\quad\quad\quad{\mathfrak N}_0=\left( \begin{array}{cc} 1 & {\mathbb Z}_p \\
    0 & 1\end{array}\right).$$For $a\in{\mathbb Z}_p^{\times}$ let us write $\gamma(a)=\left( \begin{array}{cc} a & 0 \\ 0 &
    1\end{array}\right)\in\Gamma$. Notice that $\gamma_{{{0}}}$, resp. $\nu$, is a
topological generator of $\Gamma_{{{0}}}$, resp. of ${\mathfrak N}_0$. For $r\in{\mathbb N}$ let $\lfloor{\mathfrak N}_0,\varphi^r\rfloor$
(resp. $\lfloor{\mathfrak N}_0,\varphi^r,\Gamma_{{{0}}}\rfloor$,
resp. $\lfloor{\mathfrak N}_0,\varphi^r,\Gamma\rfloor$) denote the submonoid
({\it not}: subgroup) of ${\rm GL}_2({\mathbb Q}_p)$ generated by the subgroup
${\mathfrak N}_0$ and the element $\varphi^r$ (resp. by the subgroups
${\mathfrak N}_0$ and $\Gamma_{{{0}}}$ and the element $\varphi^r$, resp. by the subgroups
${\mathfrak N}_0$ and $\Gamma$ and the element $\varphi^r$).\\

For $i\in{\mathbb Z}$ we put ${\mathfrak
  v}_i=\varphi^i({\mathfrak
  v}_0)\in{\mathfrak X}^0$ and ${\mathfrak
  e}_i=\{{\mathfrak
  v}_i,{\mathfrak
  v}_{i-1}\}\in{\mathfrak X}^1$. We denote by ${\mathfrak X}_+$ the half tree with
sets of vertices, resp. edges,$${\mathfrak X}_+^0=\coprod_{i\ge0}{\mathfrak
  N}_0\cdot{\mathfrak v}_i,\quad\quad\quad{\mathfrak
  X}_+^1=\coprod_{i\ge0}{\mathfrak N}_0\cdot{\mathfrak
  e}_i.$$Thus ${\mathfrak X}_+$ is one of the halfs obtained by cutting ${\mathfrak X}$ at the
edge ${\mathfrak e}_0$ into two pieces, where by definition ${\mathfrak e}_0$
belongs to ${\mathfrak X}_+$ as its only loose end (i.e. an edge containing only one vertex).      

The half tree $\overline{\mathfrak X}_+$ is obtained from ${\mathfrak X}_+$ by removing the edge ${\mathfrak e}_{0}$. Thus, its set of vertices is $\overline{\mathfrak X}_+^0={\mathfrak X}_+^0$, its set of edges is $\overline{\mathfrak X}_+^1={\mathfrak X}_+^1-\{{\mathfrak e}_{0}\}$, and the simplicial structure is the one induced from ${\mathfrak X}_+$.\\

For ${\mathfrak v}\in{\mathfrak X}_+^0$ we put, similarly to definition (\ref{deftub}),$$]{\mathfrak v}[\quad=\quad\{{\mathfrak w}\in {\mathfrak X}_+^0\,|\,\{{\mathfrak v},{\mathfrak
  w}\}\in {\mathfrak X}_+^1\mbox{ and }{\mathfrak v}\in [{\mathfrak w},{\mathfrak v}_0]\}.$$

\begin{satz}\label{baumeinbettung} (a) There exists an isomorphism of simplicial
  complexes$$\Theta:Y\stackrel{\cong}{\longrightarrow}{\mathfrak X}_+$$such
  that for all ${\bf v}\in Y^0$ and all $b\in{\mathbb Z}_{\ge0}$ the composition of bijections\begin{gather}]{\bf
    v}[\quad\stackrel{\Theta}{\longrightarrow}\quad]\Theta({\bf
    v})[\quad\stackrel{{\nu}^b}{\longrightarrow}\quad]\nu^b(\Theta({\bf
    v}))[\quad\stackrel{\Theta^{-1}}{\longrightarrow}\quad]\Theta^{-1}(\nu^b(\Theta({\bf
    v})))[\label{difficult}\end{gather}is induced by the action of some element $g({\bf v},b)$ of $N_0$.

(b) Let $\phi\in N(T)$ and $r\in{\mathbb
  N}$ such that $\phi(C^{(i)})=C^{(i+r)}$
for all $i\ge0$. Then $\phi(Y^1)\subset Y^1$, and the isomorphism $\Theta$ and the $g({\bf v},b)$ in (a) can be chosen in a way such that\begin{gather}\varphi^r\circ\Theta=\Theta\circ\phi\label{phivarphicom}\end{gather}and
such that in
$G$ we have \begin{gather}g(\phi({\bf
    v}),p^r)\cdot \phi=\phi\cdot g({\bf v},1).\label{secas}\end{gather}

(c) Let $\tau:{\mathbb Z}_p^{\times}\to T$ be
a homomorphism such that $\alpha^{(i)}\circ\tau={\rm id}_{{\mathbb
    Z}_p^{\times}}$ for all $i\ge0$. Then there exists for
each ${\bf v}\in Y^0$ and each $a\in{\mathbb Z}_p^{\times}$ some element $h({\bf
  v},a)$ in $N_0\cdot\tau(a)$ which induces the composition of bijections\begin{gather}]{\bf
    v}[\quad\stackrel{\Theta}{\longrightarrow}\quad]\Theta({\bf
    v})[\quad\stackrel{\gamma(a)}{\longrightarrow}\quad]\gamma(a)(\Theta({\bf
    v}))[\quad\stackrel{\Theta^{-1}}{\longrightarrow}\quad]\Theta^{-1}(\gamma(a)(\Theta({\bf
    v})))[.\label{diamond}\end{gather}For $a\in{\mathbb N}\cap{\mathbb
  Z}_p^{\times}$ these $h({\bf v},a)$ can be chosen in a way such that in
$G$ we have\begin{gather}h(g({\bf v},1){\bf v},a)\cdot g({\bf v},1)=g(h({\bf
    v},a){\bf v},a)\cdot h({\bf v},a).\label{1painch}\end{gather} 

(d) Let $\phi\in N(T)$ be as in (b), let $\tau:{\mathbb Z}_p^{\times}\to T$ be
as in (c) such that for all $a\in{\mathbb Z}_p^{\times}$ we have
$\tau(a)\phi=\phi\tau(a)$ in $G$. The isomorphism $\Theta$, the $g({\bf v},b)$ and the $h({\bf v},a)$
can be chosen in a way such that in $G$ we have 
\begin{gather}h(\phi({\bf v}),a)\cdot \phi=\phi\cdot h({\bf
    v},a).\label{2painch}\end{gather} 
\end{satz}

{\sc Proof:} (a) For any $j\ge0$ we choose a topological generator $\nu_j$ of
$N_{\alpha^{(j)}}\cap N_0\cong{\mathbb Z}_p$. For $m\in{\mathbb Z}_{\ge0}$ we
define the elements $m_0,m_1,m_2,\ldots$ of $\{0,\ldots,p-1\}$ to be the digits
of $m$ in its $p$-adic expansion, i.e. $m=\sum_{j\ge 0}m_jp^j$. We then define the
element$$\nabla(m)=\nu_0^{m_0}\cdots\nu_j^{m_jp^{e[j,\alpha^{(j)}]}}\cdot\nu_{j+1}^{m_{j+1}p^{e[j+1,\alpha^{(j+1)}]}}\cdots$$
of $N_0$. We write $[0,p^i[=\{m\in {\mathbb Z}_{\ge0}\,|\,0\le m<p^i\}$ for $i\ge0$. An easy induction on $i$ shows that
$\{\nabla(m)\,|\,m\in[0,p^i[\}$
is a set of representatives for the set of cosets $N_0/N_0^{(i)}$. Since
$N_0^{(i)}$ is the stabilizer of ${\bf
  v}_i$ in $N_0$ it follows
that the map$$\beta_{Y,i}:[0,p^i[\longrightarrow N_0\cdot{\bf
  v}_i,\quad m\mapsto \nabla(m)\cdot{\bf
  v}_i$$is bijective. On the other hand we have the bijection$$\beta_{{\mathfrak X}_+,i}:[0,p^i[\longrightarrow {\mathfrak N}_0\cdot{\mathfrak v}_i,\quad m\mapsto\nu^m\cdot{\mathfrak v}_i.$$We put $$\Theta_i=\beta_{{\mathfrak
    X}_+,i}\circ\beta_{Y,i}^{-1}\,\,:\,\,N_0\cdot{\bf
  v}_i\stackrel{\cong}{\longrightarrow}{\mathfrak N}_0\cdot{\mathfrak v}_i.$$Taken for all $i\ge0$ this is a bijection
$\Theta:Y^0\cong{\mathfrak X}_+^0$ and is easily seen to define an isomorphism of
simplicial complexes $Y\cong{\mathfrak X}_+$. In the following, for $i\ge0$
and $m\in{\mathbb Z}$ we write
$\beta_{{\mathfrak X}_+,i}(m)=\beta_{{\mathfrak X}_+,i}(m')$ and
$\beta_{Y,i}(m)=\beta_{Y,i}(m')$ where $m'\in [0,p^i[$ is such that
$m-m'\in p^i{\mathbb Z}$. For $m\in[0,p^i[$ we
find\begin{align}]\beta_{Y,i}(m)[\quad&=\quad\{\beta_{Y,i+1}(m+p^it)\,|\,0\le t\le
p-1\},\label{tubform1}\\]\beta_{{\mathfrak X}_+,i}(m)[\quad&=\quad\{\beta_{{\mathfrak
    X}_+,i+1}(m+p^it)\,|\,0\le t\le p-1\}.\label{tubform2}\end{align}For all $b\in{\mathbb Z}_{\ge0}$ we have $\nu^b(\beta_{{\mathfrak
    X}_+,i}(m))=\beta_{{\mathfrak X}_+,i}(b+m)$ in ${\mathfrak
  N}_0\cdot{\mathfrak v}_{i}$ and $$\nu^b(\beta_{{\mathfrak
    X}_+,i+1}(m+p^it))=\beta_{{\mathfrak
    X}_+,i+1}(b+m+p^it)\quad\mbox{ in }{\mathfrak N}_0\cdot{\mathfrak
  v}_{i+1}$$for all $0\le t\le p-1$. On the other hand, if for ${\bf
  v}=\beta_{Y,i}(m)$ with $m\in[0,p^i[$ we put\begin{gather}g({\bf v},b)=\nabla(b+m)\cdot\nabla(m)^{-1}\in
N_0\label{gvmfo}\end{gather}then $g({\bf v},b)\cdot
\beta_{Y,i}(m)=\beta_{Y,i}(b+m)$ in ${N}_0\cdot{\bf v}_{i}$ and$$g({\bf v},b)\cdot
\beta_{Y,i+1}(m+p^it)=\beta_{Y,i+1}(b+m+p^it)$$in
${N}_0\cdot{\bf v}_{i+1}$. To see this last equation we use Lemma \ref{petlem}. Together we obtain$$\nu^b\cdot\Theta_{i+1}({\bf
  w})=\Theta_{i+1}(g({\bf v},b)\cdot {\bf w})$$ for all ${\bf w}\in
]{\bf v}[$, as desired.

(b) For a chamber $D$ contained in $A$ let $\pi_D:X\to A$ denote the
retraction from $X$ to $A$ centered at $D$, i.e. the unique polysimplicial map
restricting to the identity on $A$ and with $\pi_D^{-1}(D)=D$. We have
$Y^1=N_0\{C^{(i)}\}_{i\ge0}=I\{C^{(i)}\}_{i\ge0}=\pi_C^{-1}(\{C^{(i)}\}_{i\ge0})$
and hence
$$\phi(Y^1)=\phi(\pi_C^{-1}(\{C^{(i)}\}_{i\ge0}))=\pi_{C^{(r)}}^{-1}(\{C^{(i)}\}_{i\ge
  r})\quad\subset\quad \pi_C^{-1}(\{C^{(i)}\}_{i\ge0})=Y^1.$$We also see that $\pi_{C^{(j+r)}}^{-1}(\{C^{(i)}\}_{i\ge
j+r})=\phi(\pi_{C^{(j)}}^{-1}(\{C^{(i)}\}_{i\ge j}))$ is stable under $\phi
(N_{\alpha^{(j)}}\cap N_0)\phi^{-1}$, for any $j\ge0$. Now $\phi$ moves the wall containing
  $C^{(j)}\cap C^{(j+1)}$ to the one containing $C^{(j+r)}\cap C^{(j+r+1)}$, therefore $\phi
(N_{\alpha^{(j)}}\cap N_0)\phi^{-1}$ must be the stabilizer of $C^{(j+r)}\cap
  C^{(j+r+1)}$ in either $N_{\alpha^{(j+r)}}$ or in
  $N_{-\alpha^{(j+r)}}$. But the stabilizer of $C^{(j+r)}\cap
  C^{(j+r+1)}$ in $N_{-\alpha^{(j+r)}}$ does not stabilize $\pi_{C^{(j+r)}}^{-1}(\{C^{(i)}\}_{i\ge
j+r})$. Thus\begin{gather}\phi(N_{\alpha^{(j)}}\cap N_0)\phi^{-1}=(N_{\alpha^{(j+r)}}\cap N_0)^{p^{e[r,\alpha^{(j+r)}]}}\label{petin}.\end{gather}Now we refine the construction in (a) by choosing the $\nu_j$ more specifically. First we claim\begin{gather}e[r,{\alpha}^{(j+r)}]+e[i,{\alpha}^{(j)}]=e[i+r,{\alpha}^{(j+r)}]\label{nepin}\end{gather}for
all $i,j\ge0$. Indeed, the above discussion shows in particular that for all
$j,t\ge0$ we have $\alpha^{(j)}=\alpha^{(t)}$ if and only if
$\alpha^{(j+r)}=\alpha^{(t+r)}$, or in other words$$\{0\le t\le i-1\,|\,\alpha^{(j+r)}=\alpha^{(t+r)}\}=\{0\le t\le
i-1\,|\,\alpha^{(j)}=\alpha^{(t)}\}.$$The set on the left hand side contains
$e[i+r,{\alpha}^{(j+r)}]-e[r,{\alpha}^{(j+r)}]$ elements, as can be seen by
performing the bijective shift
$t\mapsto t+r$. The set on the right
hand side contains $e[i,{\alpha}^{(j)}]$ elements. The claim follows. By
formula (\ref{petin}) we may
choose the topological generators $\nu_j$ of the
$N_{\alpha^{(j)}}\cap N_0$ in such a way
that\begin{gather}\phi\nu_j\phi^{-1}=\nu_{j+r}^{p^{e[r,\alpha^{(j+r)}]}}\label{epin}\end{gather}for
all $j\ge0$. An induction using formulae (\ref{nepin}) and (\ref{epin}) then
shows\begin{gather}\phi\nabla(m)\phi^{-1}=\nabla(mp^r)\label{phinab}\end{gather}for
all $m\in{\mathbb Z}_{\ge0}$. We claim
that for all $i\ge0$ we have\begin{align}\phi\circ\beta_{Y,i}&=\beta_{Y,{i+r}}\cdot p^r,\label{betaycom}\\\varphi^r\circ\beta_{{\mathfrak
    X}_+,i}&=\beta_{{\mathfrak
    X}_+,i+r}\cdot p^r\label{betaxcom}\end{align}(equality of maps
$[0,p^i[\to N_0\cdot{\bf v}_{i+r}$ resp. $[0,p^i[\to {\mathfrak
  N}_0\cdot{\mathfrak v}_{i+r}$). Indeed, formula (\ref{betaycom}) follows from formula (\ref{phinab}) and $\phi({\bf v}_{i})={\bf
  v}_{i+r}$, while formula
(\ref{betaxcom}) follows from $\varphi^r\nu\varphi^{-r}=\nu^{p^r}$ and $\varphi^r({\mathfrak v}_{i})={\mathfrak
  v}_{i+r}$. Combining formulae (\ref{betaycom}) and (\ref{betaxcom}) gives
formula (\ref{phivarphicom}), as desired.

To prove formula (\ref{secas}) we write ${\bf
  v}=\beta_{Y,i}(m)$ with $m\in[0,p^i[$ as before. By formula (\ref{betaycom}) we have $\phi({\bf
  v})=\beta_{Y,i+r}(p^rm)$ and hence $g(\phi({\bf
    v}),p^r)=\nabla(p^r(m+1))\cdot\nabla(p^rm)^{-1}$, while on the other hand
  we have $g({\bf
  v},1)=\nabla(m+1)\cdot\nabla(m)^{-1}$. Therefore it will be enough to prove$$\phi
\nabla(m+1)\cdot\nabla(m)^{-1}\phi^{-1}=\nabla(p^r(m+1))\cdot\nabla(p^rm)^{-1}$$in
$N_0$. But this follows from formula (\ref{phinab}) applied to both $m$ and
$m+1$.

(c) Let us first remark that, since $\tau(a)\in T$ in fact belongs to the
maximal compact subgroup of $T$ fixing $A$ pointwise, we have
$\tau(a)C^{(i)}=C^{(i)}$, i.e. $\tau(a){\bf v}_{i}={\bf v}_{i}$ for all $i\ge0$.

We may and do assume that $a\in{\mathbb Z}\cap{\mathbb Z}_p^{\times}$. We write ${\bf
  v}=\beta_{Y,i}(m)$ with $m\in[0,p^i[$ as before. We first claim
that\begin{gather}\Theta^{-1}\circ\gamma(a)\circ\Theta=\tau(a)\quad\quad\quad\mbox{
    in }{\rm Aut}(]{\bf v}_i[).\label{permugleich}\end{gather}To see this
observe (cf. formulae (\ref{tubform1}),
(\ref{tubform2})) that\begin{align}]{\bf v}_i[\quad&=\quad\{(\nu_i^{p^{e[i,\alpha^{(i)}]}})^b{\bf
  v}_{i+1}\,|\,0\le b\le p-1\},\notag\\]{\mathfrak v}_i[\quad&=\quad\{(\nu^{p^i})^b{\mathfrak
  v}_{i+1}\,|\,0\le b\le p-1\},\notag\end{align}and that the bijection
$\Theta:\,\,]{\bf v}_i[\to ]{\mathfrak v}_i[$ sends $(\nu_i^{p^{e[i,\alpha^{(i)}]}})^b{\bf
  v}_{i+1}$ to $(\nu^{p^i})^b{\mathfrak
  v}_{i+1}$. The hypothesis $\alpha^{(i)}\circ\tau={\rm id}_{{\mathbb
    Z}_p^{\times}}$ implies $\tau(a)n=n^{a}\tau(a)$ in $G$ 
for all $n\in N_{\alpha^{(i)}}\cap N_0$. Specifically,$$\tau(a)(\nu_i^{p^{e[i,\alpha^{(i)}]}})^b=((\nu_i^{p^{e[i,\alpha^{(i)}]}})^b)^{a}\tau(a).$$We
thus get$$\tau(a)(\nu_i^{p^{e[i,\alpha^{(i)}]}})^b{\bf
  v}_{i+1}=((\nu_i^{p^{e[i,\alpha^{(i)}]}})^b)^{a}\tau(a){\bf
  v}_{i+1}=(\nu_i^{p^{e[i,\alpha^{(i)}]}})^{ab}{\bf
  v}_{i+1}.$$On
the other hand we have$$\gamma(a)(\nu^{p^i})^b{\mathfrak
  v}_{i+1}=((\nu^{p^i})^b)^{a}\gamma(a){\mathfrak
  v}_{i+1}=(\nu^{p^i})^{ab}{\mathfrak
  v}_{i+1}.$$Comparing these formulae we get our claim (\ref{permugleich}). Now $g({\bf v}_i,am)$ resp. $g({\bf
  v}_i,m)$ induces \begin{align}\Theta^{-1}\circ\nu^{am}\circ\Theta\quad:\quad]{\bf v}_i[\quad\longrightarrow
&\quad]\Theta^{-1}(\nu^{am}(\Theta({\bf
  v}_i)))[\notag\\\mbox{
  resp. }\quad\quad\Theta^{-1}\circ\nu^{m}\circ\Theta\quad:\quad]{\bf
  v}_i[\quad\longrightarrow &\quad]\Theta^{-1}(\nu^{m}(\Theta({\bf
  v}_i)))[\quad=\quad]{\bf v}[.\notag\end{align}We may therefore rewrite the arrow
(\ref{diamond})
as\begin{align}\Theta^{-1}\circ\gamma(a)\circ\Theta&=\Theta^{-1}\circ\nu^{am}\circ\gamma(a)\circ\nu^{-m}\circ\Theta\notag\\{}&=\Theta^{-1}\circ\nu^{am}\circ\gamma(a)\circ\Theta\circ
  g({\bf
    v}_i,m)^{-1}\notag\\{}&\stackrel{(i)}{=}\Theta^{-1}\circ\nu^{am}\circ\Theta\circ\tau(a)\circ
  g({\bf
    v}_i,m)^{-1}\notag\\{}&=g({\bf
    v}_i,am)\circ\tau(a)\circ g({\bf
    v}_i,m)^{-1}\notag\end{align}where $(i)$ uses formula
(\ref{permugleich}). Inserting formula (\ref{gvmfo}) for $g({\bf v}_i,am)$ and
$g({\bf v}_i,m)$ we find that the element\begin{gather}h({\bf
    v},a)=\nabla(am)\cdot\tau(a)\cdot
  \nabla(m)^{-1}\label{hvmfo}\end{gather}of $G$ induces the arrow
(\ref{diamond}). It belongs to $N_0\cdot\tau(a)$ because $\nabla(am)$ and
$\nabla(m)^{-1}$ belong to $N_0$, and $\tau(a)$ normalizes $N_0$.

To show formula (\ref{1painch}) we write ${\bf
  v}=\beta_{Y,i}(m)=\nabla(m){\bf
  v}_i$ with $m\in[0,p^i[$ as before. From $g({\bf
  v},1)=\nabla(m+1)\nabla(m)^{-1}$ we get $g({\bf v},1){\bf
  v}=\beta_{Y,i}(m+1)$, from formula (\ref{hvmfo}) and $\tau(a){\bf v}_i={\bf v}_i$ we get $h({\bf v},a){\bf v}=\beta_{Y,i}(am)$. Thus, inserting formula (\ref{hvmfo}) again, formula (\ref{1painch})
reads\begin{align}{}&\nabla(a(m+1))\cdot\tau(a)\cdot\nabla(m+1)^{-1}\cdot\nabla(m+1)\cdot\nabla(m)^{-1}\notag\\{}=&\nabla(am+a)\cdot\nabla(am)^{-1}\cdot\nabla(am)\cdot\tau(a)\cdot\nabla(m)^{-1}\notag\end{align}which is
obviously correct. 

(d) To show formula (\ref{2painch}) for ${\bf v}=\beta_{Y,i}(m)$ we notice that $\phi({\bf
  v})=\beta_{Y,i}(p^rm)$ by formula (\ref{betaycom}). Inserting
formula (\ref{hvmfo}) for $h(\phi({\bf v}),a)$ and $h({\bf v},a)$ we see that formula
(\ref{2painch})
becomes$$\nabla(ap^rm)\cdot\tau(a)\cdot\nabla(p^rm)^{-1}\cdot\phi=\phi\cdot\nabla(am)\cdot\tau(a)\cdot\nabla(m)^{-1}.$$That
this is correct follows from formula (\ref{phinab}), applied both on the left
hand side
and on the right hand side, and our assumption $\tau(a)\phi=\phi\tau(a)$.\hfill$\Box$\\    

{\bf Remarks:} Let $N_0'$ denote the subgroup of $G$ generated by all
the $N_{\alpha^{(j)}}\cap N_0$ for $j\ge0$. The proof of Theorem \ref{baumeinbettung} shows that the elements $g({\bf v},b)$ (resp. $h({\bf
  v},a)$) of $G$ in fact belong to $N_0'$ (resp. to the subgroup of $G$ generated by $N_0'$ and $\tau(a)$).

\section{Coefficient systems on half trees}

\label{cose}

Let ${\mathfrak T}$ be a half tree with set of vertices ${\mathfrak T}^0$ and
set of edges ${\mathfrak T}^1$, as considered in section \ref{coefsec}.\\ 

{\bf Definition:} (a) A (homological) coefficient system ${\mathcal V}$ in ${\mathfrak
  o}$-modules on ${\mathfrak T}$ is a collection of ${\mathfrak
  o}$-modules ${\mathcal V}(\tau)$ for each simplex $\tau$ of ${\mathfrak T}$, and a collection
of ${\mathfrak
  o}$-linear transition maps $r_y^{\tau}:{\mathcal V}(\tau)\to {\mathcal V}(y)$ for
each $y\in {\mathfrak T}^0$, $\tau\in {\mathfrak T}^1$ with $y\in\tau$. The ${\mathfrak
  o}$-modules $H_0({\mathfrak T},{\mathcal V})$ and $H_1({\mathfrak T},{\mathcal V})$ are defined by the exact
sequence$$0\longrightarrow H_1({\mathfrak T},{\mathcal V})\longrightarrow\bigoplus_{\tau\in {\mathfrak T}^1}{\mathcal
  V}(\tau)\longrightarrow\bigoplus_{y\in {\mathfrak T}^0}{\mathcal V}(y)\longrightarrow
H_0({\mathfrak T},{\mathcal V})\longrightarrow0\label{ybarh0}$$where $v\in{\mathcal
  V}(\tau)$ is sent to $\sum_{y\in {\mathfrak T}^0\atop
  y\in\tau}r_y^{\tau}(v)$.

Morphisms of coefficient systems are defined in the obvious way. A sequence of coefficient systems
$0\to{\mathcal V}_1\to{\mathcal V}_2\to{\mathcal V}_3\to0$ is called exact if for any
simplex $\tau$ of ${\mathfrak T}$ the sequence $ 0\to{\mathcal V}_1(\tau)\to{\mathcal V}_2(\tau)\to{\mathcal V}_3(\tau)\to0$ is exact.  

(b) Let $H$ be a monoid (with neutral element $1$) acting on ${\mathfrak T}$. A coefficient
system ${\mathcal V}$ on ${\mathfrak T}$ is called $H$-equivariant if in addition we are
given an ${\mathfrak
  o}$-linear map $g_{\tau}:{\mathcal V}(\tau)\to{\mathcal V}(g\tau)$ for each
simplex $\tau$ and each $g\in H$, subject to the following conditions:

(a) $g_{h\tau}\circ h_{\tau}=(gh)_{\tau}$ for simplices $\tau$ and $g,h\in H$,

(b) $1_{\tau}={\id}_{{\mathcal V}(\tau)}$ for simplices $\tau$,

(c) $r_{gy}^{g\tau}\circ g_y=g_{\tau}\circ r_y^{\tau}$ for $y\in {\mathfrak T}^0$ and $\tau
\in {\mathfrak T}^1$ with $y\in\tau$, and $g\in H$.\\

We then we have a natural action of $H$ on $H_0({\mathfrak T},{\mathcal V})$ if (at least)
one of the following conditions is satisfied:

--- ${\mathfrak T}$ has no loose ends (to each edge two vertices are assigned), or

--- $H$ acts by automorphisms of ${\mathfrak T}$.\\

\begin{lem}\label{standtree} (a) An exact sequence of coefficient systems in ${\mathfrak
  o}$-modules $0\to{\mathcal
    V}_1\to{\mathcal V}_2\to{\mathcal V}_3\to0$ induces an exact sequence of ${\mathfrak
  o}$-modules$$0\longrightarrow H_1({\mathfrak T},{\mathcal V}_1)\longrightarrow
H_1({\mathfrak T},{\mathcal V}_2)\longrightarrow H_1({\mathfrak T},{\mathcal
  V}_3)\longrightarrow H_0({\mathfrak T},{\mathcal V}_1)\longrightarrow
H_0({\mathfrak T},{\mathcal V}_2)\longrightarrow H_0({\mathfrak T},{\mathcal
  V}_3)\longrightarrow0.$$(b) If the coefficient system ${\mathcal
    V}$ has injective transition maps $r_y^{\tau}$, then $H_1({\mathfrak T},{\mathcal
  V})=0$.
\end{lem}

{\sc Proof:} Both statements are very easy to prove.\hfill$\Box$\\  

{\bf Definition:} We say that a $N_0$-equivariant (resp. ${\mathfrak N}_0$-equivariant)
  coefficient system ${\mathcal V}$ on $Y$ (resp. on ${\mathfrak X}_+$) is of level $1$ if for any $i\ge0$ the
  action of $N_0^{(i+1)}$ on ${\mathcal V}({\bf v}_{i})$ and the
  action of $N_0^{(i)}$ on ${\mathcal V}({\bf e}_{i})$ (resp. the
  action of ${\mathfrak N}_0^{p^{i+1}}$ on ${\mathcal V}({\mathfrak v}_{i})$ and the
  action of ${\mathfrak N}_0^{p^{i}}$ on ${\mathcal V}({\mathfrak e}_{i})$) are trivial.\\

We fix a gallery (\ref{cgall}) and choose an isomorphism $\Theta:Y\stackrel{\cong}{\to}{\mathfrak
  X}_+$ as in Theorem \ref{baumeinbettung}.

\begin{satz}\label{pufo} Let ${\mathcal V}$ be a $N_0$-equivariant
  coefficient system of level $1$ on $Y$. 

(a) The push forward $\Theta_*{\mathcal V}$ of ${\mathcal V}$
  to ${\mathfrak
  X}_+$ is in a natural
  way a ${\mathfrak N}_0$-equivariant
  coefficient system of level $1$ on ${\mathfrak X}_+$.

(b) Let $\phi\in N(T)$ and $r\in{\mathbb
  N}$ be as in Theorem \ref{baumeinbettung}(b). If the action of $N_0$ on ${\mathcal V}$ extends to an action of the
submonoid of $G$ generated by $N_0$ and the element $\phi$, then $\Theta_*{\mathcal V}$ is in a natural
  way $\lfloor{\mathfrak N}_0,\varphi^r\rfloor$-equivariant. 
   
(c) Let $\phi\in N(T)$, $r\in{\mathbb
  N}$ and $\tau:{\mathbb Z}_p^{\times}\to T$ be as in Theorem \ref{baumeinbettung}(d). If the action of $N_0$ on ${\mathcal V}$ extends to an action of the
submonoid of $G$ generated by $N_0$, by the image of $\tau$ and by the element
$\phi$, then $\Theta_*{\mathcal V}$ is in a natural
  way $\lfloor{\mathfrak N}_0,\varphi^r,\Gamma\rfloor$-equivariant.

(d) In (a), resp. (b), resp. (c), the isomorphism class of $H_0(\overline{\mathfrak
  X}_+,\Theta_*{\mathcal V})$, as an ${\mathfrak o}$-module acted on by the respective
monoid, depends on the choice of (\ref{cgall}) alone, resp. of
(\ref{cgall}) and $\phi$ alone, resp. of (\ref{cgall}) and $\phi$ and $\tau$
alone, but not on the choice of $\Theta$. 
\end{satz}
 
{\sc Proof:} (a) Let $g\in{\mathfrak N}_0$ and ${\mathfrak v}\in {\mathfrak
  X}_+^0$. As ${\mathfrak N}_0$ is topologically generated by $\nu$ we find, by Theorem \ref{baumeinbettung}(a), some $g'\in {N}_0$ which induces the bijection\begin{gather}]\Theta^{-1}({\mathfrak
    v})[\quad\stackrel{\Theta}{\longrightarrow}\quad]{\mathfrak
    v}[\quad\stackrel{g}{\longrightarrow}\quad]g{\mathfrak
    v}[\quad\stackrel{\Theta^{-1}}{\longrightarrow}\quad]\Theta^{-1}(g{\mathfrak
    v})[.\label{transg}\end{gather}We define $g_{\mathfrak v}:\Theta_*{\mathcal
    V}({\mathfrak v})\to\Theta_*{\mathcal V}(g{\mathfrak v})$ to be the map $$g_{\mathfrak v}=\Theta\circ g'_{\Theta^{-1}{\mathfrak v}}\circ\Theta^{-1}\quad:\quad\Theta_*{\mathcal V}({\mathfrak v})\longrightarrow\Theta_*{\mathcal
    V}(\Theta(g'\Theta^{-1}({\mathfrak v})))=\Theta_*{\mathcal V}(g{\mathfrak v}).$$This definition is independent on the
  choice of $g'$. Indeed, let also $g''\in {N}_0$ induce the bijection
  (\ref{transg}). Then $g'^{-1}g''$ belongs to the pointwise stabilizer $H$ of $]\Theta^{-1}({\mathfrak
    v})[$ in $N_0$. But the action of $H$ on ${\mathcal V}(\Theta^{-1}{\mathfrak v})$
  is trivial. To see this we may assume, by $N_0$-equivariance of ${\mathcal
    V}$, that $\Theta^{-1}{\mathfrak v}={\bf v}_{i}$ for some $i$. Then
  $H=N_0^{(i+1)}$ and this indeed acts trivially on ${\mathcal
    V}(\Theta^{-1}{\mathfrak v})={\mathcal V}({\bf v}_{i})$ by the level $1$ assumption. Therefore $g'$ and $g''=g'(g'^{-1}g'')$ give the same map ${\mathcal V}(\Theta^{-1}{\mathfrak v})\to{\mathcal
    V}(g\Theta^{-1}{\mathfrak v})$. It follows that $g_{\mathfrak v}$ as
  defined above is well defined. 

Now let $g\in{\mathfrak N}_0$ and ${\mathfrak e}\in {\mathfrak
  X}_+^1$. Choose some $g'\in N_0$ with $g'(\Theta^{-1}({\mathfrak
  e}))=\Theta^{-1}(g({\mathfrak e}))$ and define $g_{\mathfrak
    e}:\Theta_*{\mathcal V}({\mathfrak e})\to\Theta_*{\mathcal V}(g{\mathfrak
    e})$ to be the map$$g_{\mathfrak e}=\Theta\circ g'_{\Theta^{-1}{\mathfrak e}}\circ\Theta^{-1}\quad:\quad\Theta_*{\mathcal V}({\mathfrak e})\longrightarrow\Theta_*{\mathcal
    V}(\Theta(g'\Theta^{-1}({\mathfrak e})))=\Theta_*{\mathcal V}(g{\mathfrak
    e}).$$Again this definition is independent on the
  choice of $g'$. The definitions immediately show that we have
  defined an ${\mathfrak N}_0$-action on $\Theta_*{\mathcal V}$ which is
  again of level $1$.

(b) By formula (\ref{phivarphicom}), the $\phi$-action on ${\mathcal V}$
induces a $\varphi^r$-action on $\Theta_*{\mathcal V}$. It follows from formula (\ref{secas}) --- which corresponds to the formula
$\nu^{p^r}\cdot\varphi^r=\varphi^r\cdot\nu$ in $\lfloor{\mathfrak
  N}_0,\varphi^r\rfloor$ --- that the actions of ${\mathfrak
  N}_0$ and $\varphi^r$ merge as desired.

(c) Let $a\in {\mathbb Z}_p^{\times}$
and ${\mathfrak v}\in {\mathfrak
  X}_+^0$. We define $\gamma(a)_{\mathfrak v}:\Theta_*{\mathcal
    V}({\mathfrak v})\to\Theta_*{\mathcal V}(\gamma(a){\mathfrak v})$ to be
  the map $\gamma(a)_{\mathfrak v}=\Theta\circ
  g'_{\Theta^{-1}{\mathfrak v}}\circ\Theta^{-1}$ where $g'$ is
  some element of $N_0\cdot\tau(a)$ which induces the arrow\begin{gather}]\Theta^{-1}({\mathfrak
    v})[\quad\stackrel{\Theta}{\longrightarrow}\quad]{\mathfrak
    v}[\quad\stackrel{\gamma(a)}{\longrightarrow}\quad]\gamma(a){\mathfrak
    v}[\quad\stackrel{\Theta^{-1}}{\longrightarrow}\quad]\Theta^{-1}(\gamma(a){\mathfrak
    v})[.\label{nochmal}\end{gather}Such a $g'$ does exist, by Theorem \ref{baumeinbettung}(c). As in the proof of (a) we see that
  $\gamma(a)_{\mathfrak v}$ does not depend on the choice of $g'\in
  N_0\cdot\tau(a)$ (i.e. on its left hand factor $g'\cdot\tau(a)^{-1}\in
  N_0$, as long as the arrow (\ref{nochmal}) is induced as indicated). 

Similarly, for ${\mathfrak e}\in \overline{\mathfrak X}^1_+$ and $a\in {\mathbb Z}_p^{\times}$ choose some $g'\in N_0\cdot\tau(a)$ with $g'(\Theta^{-1}({\mathfrak
  e}))=\Theta^{-1}(\gamma(a)({\mathfrak e}))$ and define $\gamma(a)_{\mathfrak
    e}:\Theta_*{\mathcal V}({\mathfrak e})\to\Theta_*{\mathcal V}(\gamma(a){\mathfrak
    e})$ to be the map $\gamma(a)_{\mathfrak e}=\Theta\circ
  g'_{\Theta^{-1}{\mathfrak e}}\circ\Theta^{-1}$. Again this is independent on
  the choice of $g'$. 

We have defined an action of $\Gamma$ on $\Theta_*{\mathcal V}$ (notice that $N_0\cdot\tau(a)\cdot
  N_0\cdot\tau(a')=N_0\cdot\tau(aa')$ for
  $a,a'\in {\mathbb Z}_p^{\times}$). It follows from formula (\ref{2painch})
  --- which corresponds to the formula
$\gamma(a)\cdot\varphi^r=\varphi^r\cdot\gamma(a)$ in $\lfloor{\mathfrak
  N}_0,\varphi^r,\Gamma\rfloor$ --- and from formula (\ref{1painch}) --- which corresponds to the formula $\gamma(a)\cdot\nu=\nu^a\cdot\gamma(a)$ in $\lfloor{\mathfrak
  N}_0,\varphi^r,\Gamma\rfloor$ --- that the actions of $\lfloor{\mathfrak
  N}_0,\varphi^r\rfloor$ and of $\Gamma$ merge as desired.

(d) Let $\Xi:Y\stackrel{\cong}{\to}{\mathfrak X}_+$ be another choice. The
automorphism $\Delta=\Theta\circ\Xi^{-1}$ of ${\mathfrak X}_+$ is
covered by the isomorphism $\Xi_*{\mathcal V}\to\Theta_*{\mathcal V}$ which
on facets $\sigma$ of ${\mathfrak X}_+$ is given by the identity maps $$(\Xi_*{\mathcal
  V})(\sigma)={\mathcal V}(\Xi^{-1}\sigma)={\mathcal V}(\Theta^{-1}\Delta\sigma)=(\Theta_*{\mathcal
  V})(\Delta\sigma).$$It induces a canonical isomorphism $H_0(\overline{\mathfrak
  X}_+,\Xi_*{\mathcal V})\to H_0(\overline{\mathfrak
  X}_+,\Theta_*{\mathcal V})$. That it commutes with the actions of
${\mathfrak N}_0$ (resp. of $\varphi^r$, resp. of $\Gamma$) is a {\it tautological} consequence of the definition of these actions.\hfill$\Box$\\  

{\bf Remark:} Theorem \ref{baumeinbettung}, and then Theorem \ref{pufo}
correspondingly, hold true if the topological generator $\nu$ of ${\mathfrak
  N}_0$ is replaced by any other topological generator of ${\mathfrak
  N}_0$. In general, this choice does affect the isomorphism class of $H_0(\overline{\mathfrak
  X}_+,\Theta_*{\mathcal V})$ as a $\lfloor{\mathfrak N}_0,\varphi^r\rfloor$-
resp. $\lfloor{\mathfrak N}_0,\varphi^r,\Gamma\rfloor$-representation (but not
as a ${\mathfrak N}_0$-representation).\\  

{\bf Remark:} Let $N_0'$ denote the subgroup of $G$ generated by all
the $N_{\alpha^{(j)}}\cap N_0$ for $j\ge0$. By the remark following Theorem
\ref{baumeinbettung} we see that Theorem
\ref{pufo} may be sharpened as follows: to obtain a ${\mathfrak N}_0$-action on $\Theta_*{\mathcal V}$ we do not need a level $1$-action of the full group $N_0$ on ${\mathcal V}$ but
only a level $1$-action (same definition) of $N_0'$; similarly for $\lfloor{\mathfrak N}_0,\varphi^r\rfloor$-, resp. $\lfloor{\mathfrak
  N}_0,\varphi^r,\Gamma\rfloor$-actions.\\

{\bf Remark:} It is straightforward to generalize Theorem \ref{baumeinbettung} (c) and (d) as
follows. In the setting of Theorem \ref{baumeinbettung} (a) let $\tau:{\mathbb Z}_p^{\times}\to T$ be
a homomorphism such that there is some
$w\in{\mathbb N}$ with $\alpha^{(i)}\circ\tau={\rm id}^w_{{\mathbb
    Z}_p^{\times}}$ for all $i\ge0$. Then there exists for
each ${\bf v}\in Y^0$ and each $a\in{\mathbb Z}_p^{\times}$ some element $h({\bf
  v},a)$ in $N_0\cdot\tau(a)$ which induces the composition of bijections$$]{\bf
    v}[\quad\stackrel{\Theta}{\longrightarrow}\quad]\Theta({\bf
    v})[\quad\stackrel{\gamma(a)^w}{\longrightarrow}\quad]\gamma(a)^w(\Theta({\bf
    v}))[\quad\stackrel{\Theta^{-1}}{\longrightarrow}\quad]\Theta^{-1}(\gamma(a)^w(\Theta({\bf
    v})))[.$$For $a\in{\mathbb N}\cap{\mathbb
  Z}_p^{\times}$ these $h({\bf v},a)$ can be chosen in a way such that in
$G$ we have$$h(g({\bf v},1){\bf v},a)\cdot g({\bf v},1)=g(h({\bf
    v},a){\bf v},a^w)\cdot h({\bf v},a).$$Moreover, if $\phi\in
N(T)$ is as in Theorem \ref{baumeinbettung} (b) such that for all $a\in{\mathbb Z}_p^{\times}$ we have
$\tau(a)\phi=\phi\tau(a)$ in $G$, then, as before, we can achieve $h(\phi({\bf v}),a)\cdot \phi=\phi\cdot h({\bf
    v},a)$. Given this setting, Theorem \ref{pufo} (c)
can be generalized as follows. If the action of $N_0$ on ${\mathcal V}$ extends to an action of the
submonoid of $G$ generated by $N_0$, by the image of $\tau$ and by
$\phi$, then $\Theta_*{\mathcal V}$ is in a natural way both $\lfloor{\mathfrak
  N}_0,\varphi^r\rfloor$-equivariant and $\Gamma$-equivariant,
such that$$\gamma(a)\cdot\varphi^r=\varphi^r\cdot\gamma(a)\quad\quad\quad\mbox{ and
}\quad\quad\quad\gamma(a)\cdot\nu=\nu^{a^w}\cdot\gamma(a)$$as endomorphisms of
$\Theta_*{\mathcal V}$. Therefore, if in addition $\tau(a)$ acts
trivially on ${\mathcal V}$ for all $a\in{\mathbb Z}_p^{\times}$ with $a^w=1$,
then by extracting $w$-th roots of the above $\gamma(a)$-operators we obtain an action of the submonoid of $\lfloor{\mathfrak
  N}_0,\varphi^r,\Gamma\rfloor$ generated by ${\mathfrak N}_0$, by
$\varphi^r$ and by $\Gamma^w=\{\gamma^w\,|\,\gamma\in\Gamma\}$.\\

{\bf Definition:} Let $m\ge1$, let ${\mathcal V}$ be a ${\mathfrak N}_0$-equivariant coefficient system in ${\mathfrak
  o}_m$-modules on ${\mathfrak X}_+$. 

(1) We say that ${\mathcal
  V}$ is strictly of level $1$ over $k$ if the ${\mathfrak
  o}_m$-action factors through $k$ and if the following
conditions (a), (b) and (c) are satisfied:   

(a) All transition maps $r_y^{\tau}:{\mathcal V}(\tau)\to {\mathcal V}(y)$ for
$y\in {\mathfrak X}_+^0$, $\tau\in {\mathfrak X}_+^1$ with $y\in\tau$ are injective; in the following we view them as
inclusions.

(b) For all $i\ge0$ we have \begin{gather}{\mathcal V}({\mathfrak
    v}_{i})={\mathfrak N}_0^{p^{i}}\cdot{\mathcal V}({\mathfrak
    e}_{i+1}),\label{striclevgen}\end{gather}i.e. as a
${\mathfrak N}_0^{p^{i}}$-representation, ${\mathcal V}({\mathfrak v}_{i})$ is
generated by ${\mathcal V}({\mathfrak e}_{i+1})$, and\begin{gather}{\mathcal V}({\mathfrak e}_{i})={\mathcal V}({\mathfrak v}_{i})^{{\mathfrak N}_0^{p^{i}}},\label{striclevinv}\end{gather}i.e. ${\mathcal V}({\mathfrak e}_{i})$ is the submodule of ${\mathfrak N}_0^{p^{i}}$-invariants in ${\mathcal V}({\mathfrak v}_{i})$.

(c) The number ${\rm dim}_k({\mathcal V}({\mathfrak
  e}_{i}))$
is finite and independent of $i\ge0$.

(2) We say that ${\mathcal V}$ is strictly of level
$1$ if it admits a finite filtration such that all the subquotiens are strictly of level
$1$ over $k$.

\begin{satz}\label{finite} Let ${\mathcal V}$ be an ${\mathfrak N}_0$-equivariant coefficient system in ${\mathfrak
  o}_m$-modules on
${\mathfrak X}_+$ which is strictly of level $1$. We have a natural isomorphism$${\mathcal V}({\mathfrak
  e}_{0})\stackrel{\cong}{\longrightarrow} H_0(\overline{\mathfrak X}_+,{\mathcal V})^{{\mathfrak N}_0}.$$
\end{satz}

{\sc Proof:} (i) Suppose first that ${\mathcal V}$ is strictly of level
$1$ over $k$. We have canonical isomorphisms of ${\mathfrak
  N}_0$-representations\begin{gather}\bigoplus_{n\in {\mathfrak
      N}_0/{\mathfrak N}_0^{p^i}}{\mathcal V}(n\cdot{\mathfrak
    v}_{i})\cong{\rm ind}_{{\mathfrak N}_0^{p^{i}}}^{{\mathfrak N}_0}
  {\mathcal V}({\mathfrak
    v}_{i}),\label{preshave}\end{gather}\begin{gather}\bigoplus_{n\in
    {\mathfrak N}_0/{\mathfrak N}_0^{p^i}}{\mathcal V}(n\cdot{\mathfrak
    e}_{i})\cong{\rm ind}_{{\mathfrak N}_0^{p^{i}}}^{{\mathfrak N}_0}{\mathcal
    V}({\mathfrak e}_{i}).\label{preshaed}\end{gather}By formula (\ref{striclevinv}) the transition map ${\mathcal V}({\mathfrak
  e}_{0})\to{\mathcal V}({\mathfrak
  v}_{0})$ is an isomorphism between ${\mathcal V}({\mathfrak
  e}_{0})$ and ${\mathcal V}({\mathfrak
  v}_{0})^{{\mathfrak N}_0}$. Hence we need to show that the natural map ${\mathcal V}({\mathfrak
  v}_{0})^{{\mathfrak N}_0}\to H_0(\overline{\mathfrak X}_+,{\mathcal
  V})^{{\mathfrak N}_0}$ is an isomorphism. The
injectivity even of ${\mathcal V}({\mathfrak
  v}_{0})\to H_0(\overline{\mathfrak X}_+,{\mathcal
  V})$ follows from the
injectivity of all the transition maps of ${\mathcal V}$. Next, Lemma
\ref{standtree} (which of course relies on the same argument) shows the exactness of the
sequence$$0\longrightarrow\bigoplus_{\tau\in \overline{\mathfrak X}_+^1}{\mathcal
  V}(\tau)\longrightarrow\bigoplus_{y\in \overline{\mathfrak X}_+^0}{\mathcal V}(y)\longrightarrow
H_0(\overline{\mathfrak X}_+,{\mathcal V})\longrightarrow0.$$Looking at the long
exact cohomology sequence obtained by applying the group cohomology functor
$H^*({\mathfrak N}_0,\cdot)$ we see
that it is now enough to prove that the natural map\begin{gather}(\bigoplus_{\tau\in \overline{\mathfrak X}_+^1}{\mathcal V}(\tau))^{{\mathfrak N}_0}\bigoplus {\mathcal V}({\mathfrak
  v}_{0})^{{\mathfrak N}_0}\longrightarrow(\bigoplus_{y\in \overline{\mathfrak X}_+^0}{\mathcal
  V}(y))^{{\mathfrak N}_0}\label{h0h0yb}\end{gather}is bijective and that the
natural map\begin{gather}H^1({\mathfrak N}_0,\bigoplus_{\tau\in \overline{\mathfrak X}_+^1}{\mathcal
  V}(\tau))\longrightarrow H^1({\mathfrak N}_0,\bigoplus_{y\in
  \overline{\mathfrak X}_+^0}{\mathcal V}(y))\label{h1h0yb}\end{gather}is injective. We recognize the map (\ref{h0h0yb}) as the natural map\begin{gather}(\bigoplus_{\tau\in {\mathfrak X}_+^1}{\mathcal V}(\tau))^{{\mathfrak N}_0}\longrightarrow(\bigoplus_{y\in {\mathfrak X}_+^0}{\mathcal
  V}(y))^{{\mathfrak N}_0}.\label{newh0h0yb}\end{gather}In view of the
isomorphisms (\ref{preshave}), (\ref{preshaed}) and Shapiro's Lemma, we may
rewrite it as$$\bigoplus_{i\ge 0}{\mathcal V}({\mathfrak
  e}_{i})^{{{\mathfrak N}_0^{p^{i}}}}\longrightarrow\bigoplus_{i\ge0}{\mathcal V}({\mathfrak
  v}_{i})^{{{\mathfrak N}_0^{p^{i}}}}$$(with ${\mathcal V}({\mathfrak
  e}_{0})^{{{\mathfrak N}_0^{p^{0}}}}$ mapping to ${\mathcal V}({\mathfrak
  v}_{0})^{{{\mathfrak N}_0^{p^{0}}}}$ and with ${\mathcal V}({\mathfrak
  e}_{i})^{{{\mathfrak N}_0^{p^{i}}}}$ for $i\ge1$ mapping to both ${\mathcal V}({\mathfrak
  v}_{i})^{{{\mathfrak N}_0^{p^{i}}}}$ and ${\mathcal V}({\mathfrak
  v}_{i-1})^{{{\mathfrak N}_0^{p^{i-1}}}}$). Since by hypothesis all the
maps $${\mathcal V}({\mathfrak
  e}_{i})^{{{\mathfrak N}_0^{p^{i}}}}={\mathcal V}({\mathfrak
  e}_{i})\longrightarrow{\mathcal V}({\mathfrak
  v}_{i})^{{{\mathfrak N}_0^{p^{i}}}}$$ are bijective we have proven
the bijectivity of the map (\ref{h0h0yb}). Similarly as with the map
(\ref{newh0h0yb}) we proceed with the map (\ref{h1h0yb}): we use Shapiro's Lemma to rewrite it
as$$\bigoplus_{i\ge0}H^1({\mathfrak N}_0^{p^{i}},k[\frac{{\mathfrak N}_0^{p^{i}}}{{\mathfrak N}_0^{p^{i+1}}}]\otimes_k{\mathcal V}({\mathfrak
  e}_{i+1}))\longrightarrow\bigoplus_{i\ge0}H^1({\mathfrak N}_0^{p^{i}},{\mathcal V}({\mathfrak
  v}_{i})).$$To see its injectivity it is enough to
show that the maps$$H^1({\mathfrak N}_0^{p^{i}},k[\frac{{\mathfrak N}_0^{p^{i}}}{{\mathfrak N}_0^{p^{i+1}}}]\otimes_k{\mathcal V}({\mathfrak
  e}_{i+1}))\longrightarrow H^1({\mathfrak N}_0^{p^{i}},{\mathcal V}({\mathfrak
  v}_{i}))$$are injective for all $i\ge0$. But in view of our
hypotheses on ${\mathcal V}$ this follows from  Lemma \ref{einsbcjalg},
applied to ${\mathfrak N}_0^{p^{i}}\cong{\mathbb Z}_p$ acting through
its quotient
${\mathfrak N}_0^{p^{i}}/{\mathfrak N}_0^{p^{i+1}}\cong {\mathbb F}_p$ on ${\mathcal V}({\mathfrak
  v}_{i})$. 

(ii) For general ${\mathcal V}$ strictly of level
$1$ we argue by induction on the minimal length of a filtration with subquotients
strictly of level $1$ over $k$. Consider an exact sequence $0\to{\mathcal
  V}_1\to{\mathcal V}\to{\mathcal V}_2\to0$ where ${\mathcal
  V}_1$ and ${\mathcal V}_2$ are strictly of level $1$ and different from ${\mathcal V}$. The strict level $1$ property of ${\mathcal V}_2$ implies
that ${\mathcal V}_2$ has injective transition maps. By Lemma
\ref{standtree} we obtain the exactness of\begin{gather}0\longrightarrow H_0(\overline{\mathfrak
    X}_+,{\mathcal V}_1)\longrightarrow H_0(\overline{\mathfrak
    X}_+,{\mathcal V})\longrightarrow H_0(\overline{\mathfrak
    X}_+,{\mathcal V}_2)\longrightarrow0.\label{h0ex}\end{gather}It follows that in the commutative diagram$$\xymatrix{0\ar[r]& {\mathcal V}_1({\mathfrak
  e}_{0})               \ar[d]\ar[r]& {\mathcal V}({\mathfrak
  e}_{0})     \ar[d]\ar[r]&{\mathcal V}_2({\mathfrak
  e}_{0})        \ar[d]\ar[r]&0\\0\ar[r]& H_0(\overline{\mathfrak
    X}_+,{\mathcal V}_1)^{{\mathfrak N}_0}   \ar[r]&H_0(\overline{\mathfrak
    X}_+,{\mathcal V})^{{\mathfrak N}_0}\ar[r]&H_0(\overline{\mathfrak X}_+,{\mathcal
    V}_2)^{{\mathfrak N}_0} 
       &}$$the bottom sequence is exact. For the top sequence this is
     obvious. The outer vertical arrows are bijective by induction hypothesis. Hence the middle vertical arrow is
     bijective.\hfill$\Box$\\     

\begin{pro}\label{natgam} Let $r\in{\mathbb N}$, let ${\mathcal V}$ be a
  $\lfloor{\mathfrak N}_0,\varphi^r\rfloor$-equivariant coefficient system on
  ${\mathfrak X}_+$ whose ${\mathfrak N}_0$-action is of level $1$. The
  $\lfloor{\mathfrak N}_0,\varphi^r\rfloor$-action on ${\mathcal V}$ naturally
  extends to a $\lfloor{\mathfrak N}_0,\varphi^r,\Gamma_{{{0}}}\rfloor$-action such that $\Gamma_0$ acts trivially on ${\mathcal V}({\mathfrak e}_0)$. 
\end{pro}

{\sc Proof:} The action of $\Gamma_{{{0}}}$ on ${\mathfrak X}_+^0$ respects the orbits
${\mathfrak N}_0\cdot {\mathfrak v}_i$ for all $i\ge0$. On such an orbit
${\mathfrak N}_0\cdot {\mathfrak v}_i$ it is described by the
formula\begin{gather}\gamma_{{{0}}}(\nu^n({\mathfrak v}_i))=\nu^{(p+1)n}({\mathfrak
    v}_i)\label{gamfor}\end{gather}for all $n\ge0$. Applied to the elements
of$$]\nu^n({\mathfrak v}_i)[\quad=\quad\{\nu^{n+tp^i}({\mathfrak v}_{i+1})\,|\,0\le
t<p\},$$formula (\ref{gamfor}) (with $i+1$ instead of $i$) shows that the
restrictions of $\gamma_{{{0}}}$ and $\nu^{np}$ to the subset $]\nu^n({\mathfrak
  v}_i)[$ of ${\mathfrak N}_0\cdot {\mathfrak v}_{i+1}$ coincide, because $(p+1)(n+tp^i)\equiv pn+n+tp^i$ modulo $(p^{i+1})$. Thus we
obtain that for all $j\in{\mathbb Z}$ and all ${\mathfrak v}\in {\mathfrak
  X}_+^0$ there is some $\beta\in{\mathfrak N}_0$ such that the restrictions
of $\gamma_{{{0}}}^j$ and $\beta$ to the subset $]{\mathfrak v}[$ of ${\mathfrak
  X}_+^0$ coincide. Therefore, if we define the map $$\gamma^j_{{{0}},{\mathfrak
    v}}:{\mathcal V}({\mathfrak v})\longrightarrow{\mathcal
  V}(\gamma_{{{0}}}^j({\mathfrak v}))$$ as the map $\beta_{{\mathfrak v}}:{\mathcal
  V}({\mathfrak v})\to{\mathcal V}(\beta({\mathfrak v}))={\mathcal
  V}(\gamma_{{{0}}}^j({\mathfrak v}))$ then, arguing as in the proof of Theorem
\ref{pufo}, we see that by the level $1$ property of ${\mathcal V}$ this definition is indepedent on the choice of $\beta$. In the same way we define maps $\gamma^j_{{0},{\mathfrak e}}:{\mathcal V}({\mathfrak e})\longrightarrow{\mathcal V}(\gamma_{{{0}}}^j({\mathfrak e}))$ for ${\mathfrak e}\in{\mathfrak X}_+^1$. We have defined an action of $\Gamma_{{{0}}}$ on ${\mathcal V}$. 

Now let again $\nu^n({\mathfrak v}_i)$ (some $n\ge0$, some $i\ge0$) be an
arbitrary element of ${\mathfrak X}_+^0$. Let $x\in {\mathcal
  V}(\nu^n({\mathfrak v}_i))$. We compute (dropping the names of vertices in subscripts)$$\gamma_{{{0}}}(\nu(x))\stackrel{(i)}{=}\nu^{p(n+1)}(\nu(x))=\nu^{p+1}(\nu^{pn}(x))\stackrel{(ii)}{=}\nu^{p+1}(\gamma_{{{0}}}(x))$$where
in $(i)$ we used that $\gamma_{{{0}}}$ and $\nu^{p(n+1)}$ coincide on
$]\nu^{n+1}({\mathfrak v}_i)[$ whereas in $(ii)$ we used that $\gamma_{{{0}}}$ and
$\nu^{pn}$ coincide on $]\nu^{n}({\mathfrak v}_i)[$. Similarly we
compute$$\gamma_{{{0}}}(\varphi^r(x))\stackrel{(i)}{=}\nu^{p^{r+1}n}(\varphi^r(x))=\varphi^r(\nu^{pn}(x))\stackrel{(ii)}{=}\varphi^r(\gamma_{{{0}}}(x))$$where
in $(i)$ we used that $\gamma_{{{0}}}$ and $\nu^{p^{r+1}n}$ coincide on
$]\nu^{p^rn}({\mathfrak v}_{i+r})[$ whereas in $(ii)$ we used that $\gamma_{{{0}}}$
and $\nu^{pn}$ coincide on $]\nu^{n}({\mathfrak v}_i)[$. A similar computation
can be done on the values of ${\mathcal V}$ at all ${\mathfrak e}\in{\mathfrak
  X}_+^1$. We have shown\begin{gather}\gamma_{{{0}}}\circ\varphi^r=\varphi^r\circ\gamma_{{{0}}}\quad\quad\quad\mbox{
  and }\quad\quad\quad\gamma_{{{0}}}\circ\nu=\nu^{p+1}\circ\gamma_{{{0}}}\label{ganuph}\end{gather}as endomorphisms of
${\mathcal V}$. This means that the actions of $\lfloor{\mathfrak
  N}_0,\varphi^r\rfloor$ and of $\Gamma_{{{0}}}$ merge to an action of
$\lfloor{\mathfrak N}_0,\varphi^r,\Gamma_{{{0}}}\rfloor$ on ${\mathcal V}$.  \hfill$\Box$\\    

\section{Pro-$p$ Iwahori-Hecke modules and coefficient systems}

\label{iwah}

Let $I_0$ denote the maximal pro-$p$ subgroup in $I$. Let ${\rm
  ind}_{I_0}^{G}{\bf 1}_{{\mathfrak o}}$ denote the ${\mathfrak o}$-module of
${\mathfrak o}$-valued compactly supported functions $f$ on $G$ such that $f(ig)=f(g)$ for all
$g\in G$, all $i\in I_0$. It is a $G$-representation by means of
$(g'f)(g)=f(gg')$ for $g,g'\in G$. Let $${\mathcal
  H}(G,I_0)={\rm End}_{{\mathfrak
    o}[G]}({\rm ind}_{I_0}^{G}{\bf 1}_{{\mathfrak o}})^{\rm op}$$denote the corresponding
pro-$p$-Iwahori Hecke algebra with coefficients in ${\mathfrak o}$. Then ${\rm
  ind}_{I_0}^{G}{\bf 1}_{{\mathfrak o}}$ is naturally a right ${\mathcal
  H}(G,I_0)$-module. For a subset $H$ of $G$ let $\chi_H$ denote the
characteristic function of $H$. For $g\in G$ let $T_g\in {\mathcal
  H}(G,I_0)$ denote the Hecke operator corresponding to the double coset
$I_0gI_0$. It sends $f:G\to{\mathfrak o}$
to $$T_g(f):G\longrightarrow{\mathfrak o},\quad\quad h\mapsto\sum_{x\in
  I_0\backslash G}\chi_{I_0gI_0}(hx^{-1})f(x).$$In particular we have\begin{gather}T_g(\chi_{I_0})=\chi_{I_0g}=g^{-1}\chi_{I_0}\quad\quad\mbox{ if
}gI_0=I_0g.\label{hecnor}\end{gather}

For any facet $F$ of $X$ let $I_0^F$ denote the 'pro-unipotent radical' of the stabilizer of $F$. More precisely, following \cite{os} section 3 (where instead the notation $I_F$ is used), $I_0^F$ consists of all $g\in{\bf G}_F^0({\mathbb Z}_p)$ mapping to the unipotent radical in $\overline{\bf G}_F^0$. Here ${\bf G}_F$ is the smooth affine ${\mathbb Z}_p$-group scheme whose general fibre is (the reductive ${\mathbb Q}_p$-group scheme underlying) $G$ and such that ${\bf G}_F({\mathbb Z}_p)$ is the pointwise stabilizer of the preimage of $F$ in the enlarged building of $G$. For example, $I_0=I_0^C$. Since $F_1\subset F_2$ implies $I_0^{F_1}\subset I_0^{F_2}$ the assignment $F\mapsto ({\rm ind}_{I_0}^{G}{\bf 1}_{{\mathfrak o}})^{I_0^F}$ is a $G$-equivariant coefficient system on $X$. Since the right action of ${\mathcal
  H}(G,I_0)$ on ${\rm ind}_{I_0}^{G}{\bf 1}_{{\mathfrak o}}$ commutes with the left $G$-action this is, in fact, also a coefficient system of right ${\mathcal
  H}(G,I_0)$-modules. Given a left ${\mathcal
  H}(G,I_0)$-module $M$ we therefore obtain a new $G$-equivariant coefficient system ${\mathcal V}_M^X$ on $X$ by putting$$ {\mathcal V}_M^X(F)=({\rm ind}_{I_0}^{G}{\bf 1}_{{\mathfrak o}})^{I_0^F}\otimes_{{\mathcal
  H}(G,I_0)}M.$$In fact we will only be interested in the restriction of ${\mathcal V}_M^X$ to facets of codimension $0$ and $1$. We may regard the restriction of ${\mathcal
  V}^X_M$ to $Y$ as a coefficient system ${\mathcal
  V}_M$ on $Y$. To be explicit, for $n\in N_0$, $i\ge0$, to the vertex
$n\cdot{\bf v}_{i}$, resp. the edge $n \cdot{\bf e}_{i}$, of $Y$, we assign the codimension-$1$-facet, resp. codimension-$0$-facet, $$\eta^0(n\cdot{\bf v}_{i})=n C^{(i)}\cap n C^{(i+1)},\quad\quad\quad\mbox{ resp. }\quad\quad\quad \eta^1(n\cdot{\bf e}_{i})=n C^{(i)},$$of $X$. Both $\eta^0$ and $\eta^1$ are injective mappings, equivariant under all sub monoids of $G$
which respect $Y$, and $(\eta^0,\eta^1)$ respects facet inclusions. We now put$${\mathcal
  V}_{M}({\bf v})={\mathcal
  V}^X_{M}(\eta^0({\bf v})),\quad\quad\quad\mbox{ resp. }\quad\quad\quad{\mathcal
  V}_{M}({\bf e})={\mathcal
  V}^X_{M}(\eta^1({\bf e})),$$for ${\bf v}\in Y^0$, resp. ${\bf e}\in Y^1$. This defines a coefficient system ${\mathcal
  V}_{M}$ on $Y$, equivariant under all sub monoids of $G$
which respect $Y$. 

It is clear that all these constructions are covariantly functorial in $M$.\\

For $m\ge1$ let us write ${\mathcal H}(G,I_0)_{{\mathfrak
  o}_{m}}={\mathcal
  H}(G,I_0)\otimes_{\mathfrak o}{\mathfrak o}_m$. We denote by ${\rm Mod}^{\rm
  fin}({\mathcal
  H}(G,I_0)_{{\mathfrak o}_m})$ the category of ${\mathcal
  H}(G,I_0)_{{\mathfrak o}_m}$-modules which are finitely generated as ${\mathfrak o}_m$-modules.

We fix a gallery
(\ref{cgall}) and choose an isomorphism $\Theta:Y\stackrel{\cong}{\to}{\mathfrak
  X}_+$ as in Theorem \ref{baumeinbettung}.\\

Let $F$ be a codimension-$1$-face of $C$. There is a unique quotient $\overline{\mathcal S}$ of $I_0^F$ which is isomorphic with either ${\rm SL}_2({\mathbb
  F}_p)$ or ${\rm PSL}_2({\mathbb
  F}_p)$. The image of $I_0$ in $\overline{\mathcal S}$ is the unipotent
radical $\overline{\mathcal U}$ of a Borel subgroup in $\overline{\mathcal S}$.

\begin{pro}\label{jalgheck} (a) $\Theta_*{\mathcal V}_M$ is in a natural way a ${\mathfrak N}_0$-equivariant coefficient
  system of level $1$ on ${\mathfrak X}_+$.

(b)  For any left-${\mathcal
  H}(G,I_0)_{{\mathfrak
  o}_{m}}$-module $M$ the Hecke algebra ${\mathcal H}(\overline{\mathcal
  S},\overline{\mathcal U})={\rm End}_{{\mathfrak o}[\overline{\mathcal S}]}({\rm
  ind}_{\overline{\mathcal U}}^{\overline{\mathcal S}}{\bf 1}_{\mathfrak o})^{\rm op}$ naturally acts on $M$, and we have natural isomorphisms $M\cong{\mathcal V}_M^X(C)$ and \begin{gather}{\rm ind}_{\overline{\mathcal U}}^{\overline{\mathcal S}}{\bf 1}_{\mathfrak
  o}\otimes_{{\mathcal H}(\overline{\mathcal
  S},\overline{\mathcal U})}M\cong{\mathcal
  V}^X_M(F)\label{rapecomp}\end{gather}such that the transition map ${\mathcal V}_M^X(C)\to{\mathcal V}_M^X(F)$ gets identified with the natural map \begin{gather}M\longrightarrow {\rm ind}_{\overline{\mathcal U}}^{\overline{\mathcal S}}{\bf 1}_{\mathfrak
  o}\otimes_{{\mathcal H}(\overline{\mathcal
  S},\overline{\mathcal U})}M.\label{ufer}\end{gather}

(c) If $M\in {\rm Mod}^{\rm fin}({\mathcal
  H}(G,I_0)_{{\mathfrak o}_m})$ then $\Theta_*{\mathcal V}_M$ is strictly of level $1$.

(d) If $0\to M_1\to M\to M_2\to0$ is an exact sequence of ${\mathcal H}(G,I_0)_{{\mathfrak
  o}_{m}}$-modules for some
$m\in{\mathbb N}$ then the induced sequence $0\to \Theta_*{\mathcal
  V}_{M_1}\to\Theta_*{\mathcal V}_{M}\to\Theta_*{\mathcal V}_{M_2}\to0$ is
exact.
\end{pro}

{\sc Proof:} (a) This follows from Theorem \ref{pufo}.

(b) For a facet $D$ of $X$ let $J_D$ denote the stabilizer of $D$ in $G$. Put $${\mathcal H}_{F}={\rm End}_{{\mathfrak o}[J_F]}({\rm ind}_{I_0}^{J_F}{\bf 1}_{{\mathfrak o}})^{\rm op}.$$As $I_0^F\subset J_F$ we have natural embeddings$${\rm ind}_{\overline{\mathcal U}}^{\overline{\mathcal S}}{\bf 1}_{\mathfrak
  o}\cong {\rm ind}_{I_0}^{I_0^F}{\bf 1}_{\mathfrak
  o}\hookrightarrow{\rm ind}_{I_0}^{J_F}{\bf 1}_{{\mathfrak o}}\hookrightarrow{\rm ind}_{I_0}^{G}{\bf 1}_{{\mathfrak o}}.$$Thus Frobenius reciprocity provides inclusions of Hecke algebras$${\mathcal H}(\overline{\mathcal
  S},\overline{\mathcal U})\subset {\mathcal H}_{F}\subset{\mathcal
  H}(G,I_0).$$We obtain a natural map\begin{gather}{\rm ind}_{\overline{\mathcal U}}^{\overline{\mathcal S}}{\bf 1}_{\mathfrak
  o}\otimes_{{\mathcal H}(\overline{\mathcal
  S},\overline{\mathcal U})}M\longrightarrow {\rm ind}_{I_0}^{J_F}{\bf 1}_{{\mathfrak o}}\otimes_{{\mathcal H}_{F}}M.\label{kleiner}\end{gather}The elements of $J_C$ normalize $I_0$,
 hence we have the group
 homomorphism $J_C\longrightarrow {\mathcal
  H}(G,I_0)^{\times}$, $g\mapsto T_{g^{-1}}$, into the group ${\mathcal
  H}(G,I_0)^{\times}$ of
invertible elements of ${\mathcal
  H}(G,I_0)$. [Indeed, for $g_1, g_2\in J_C$ we form the composition $T_{g_1}\circ T_{g_2}$ in
${\rm End}_{{\mathfrak
    o}[G]}({\rm ind}_{I_0}^{G}{\bf 1}_{{\mathfrak o}})$ and compute
$(T_{g_1}\circ
T_{g_2})(\chi_{I_0})=T_{g_1}(\chi_{I_0g_2})=T_{g_1}(g_2^{-1}\chi_{I_0})=g_2^{-1}T_{g_1}(\chi_{I_0})=g_2^{-1}\chi_{I_0g_1}=g_2^{-1}g_1^{-1}\chi_{I_0}=\chi_{I_0g_1g_2}$.] Similarly we have a homomorphism $J_C\cap
J_F\longrightarrow {\mathcal
  H}_F^{\times}$. As $J_C\cap
J_F$ and $I_0^F$ together generate $J_F$ and as their intersection is $I_0$ we deduce that the map
(\ref{kleiner}) is an isomorphism. But also the map\begin{gather}({\rm ind}_{I_0}^{J_F}{\bf 1}_{{\mathfrak o}})\otimes_{{\mathcal H}_{F}}{\mathcal
  H}(G,I_0)_{{\mathfrak
  o}_{m}}\longrightarrow ({\rm ind}_{I_0}^{G}{\bf 1}_{{\mathfrak
  o}_m})^{I_0^F},\quad\quad f\otimes h\mapsto h(f)\label{pera}\end{gather}is
an isomorphism. This is proven in Lemma 3.11.i and section 4.9 of \cite{os} in
the setting where the coefficient ring is a field (instead of our ${\mathfrak
  o}_m$). However, the proof given in loc. cit. for coefficient fields of
characteristic $p$ also applies to the coefficient ring ${\mathfrak o}_m$. Composing the isomorphism (\ref{kleiner}) with the isomorphism obtained from applying $(.)\otimes_{{\mathcal
  H}(G,I_0)_{{\mathfrak
  o}}} M$ to (\ref{pera}) we obtain the desired isomorphism (\ref{rapecomp}). On the other hand, as $I_0=I_0^C$ we have ${\mathcal
  H}(G,I_0)\cong  ({\rm ind}_{I_0}^{G}{\bf 1}_{{\mathfrak o}})^{I_0^C}$ and hence $M\cong{\mathcal V}_M^X(C)$ naturally. By construction, the transition map ${\mathcal V}_M^X(C)\to{\mathcal V}_M^X(F)$ gets identified with the map (\ref{ufer}).

(c) By Lemma \ref{sl2einbett} the map (\ref{ufer}) is an isomorphism between $M$ and $({\rm ind}_{\overline{\mathcal U}}^{\overline{\mathcal S}}{\bf 1}_{\mathfrak
  o}\otimes_{{\mathcal H}(\overline{\mathcal
  S},\overline{\mathcal U})}M)^{\overline{\mathcal U}}$. In other words, $M$ is the submodule of $I_0$-invariants in the $I_0^F$-representation ${\mathcal
  V}^X_M(F)$. By $G$-equivariance of ${\mathcal
  V}^X_M$ we deduce property (\ref{striclevinv}) for
$\Theta_*{\mathcal V}_M$. Similarly, Lemma \ref{einsajalg} ensures property
(\ref{striclevgen}) for $\Theta_*{\mathcal V}_M$. All this applies of course similarly
to the subquotients of $M$ with respect to the filtration $\{p^iM\}_{i\ge0}$. The other properties
required for being strictly of level $1$ are clear. 

(d) The exactness at any
edge of ${\mathfrak X}_+$ (which corresponds to a chamber of $Y$) follows
immediately from the definitions. To prove exactness at any vertex of
${\mathfrak X}_+$ it is enough, by $G$-equivariance, to prove exactness of $0\to {\mathcal
  V}^X_{M_1}(F)\to{\mathcal V}^X_{M}(F)\to{\mathcal V}^X_{M_2}(F)\to0$ for all codimension-$1$-faces $F$ of $C$. From what we
learned in the proof of (b) and (c) we see that here
we need to prove exactness of$$0\longrightarrow{\rm ind}_{\overline{\mathcal U}}^{\overline{\mathcal
  S}}{\bf 1}_{{\mathfrak
  o}_{m}}\otimes M_1\stackrel{}{\longrightarrow}{\rm ind}_{\overline{\mathcal U}}^{\overline{\mathcal
  S}}{\bf 1}_{{\mathfrak
  o}_{m}}\otimes M
     {\longrightarrow}{\rm ind}_{\overline{\mathcal U}}^{\overline{\mathcal
  S}}{\bf 1}_{{\mathfrak
  o}_{m}}\otimes M_2{\longrightarrow}0$$where all tensor products are taken
over ${\mathcal H}(\overline{\mathcal S},\overline{\mathcal U})_{{\mathfrak
  o}_{m}}={\rm End}_{{\mathfrak
  o}_{m}[\overline{\mathcal
  S}]}({\rm ind}_{\overline{\mathcal U}}^{\overline{\mathcal
  S}}{\bf 1}_{{\mathfrak
  o}_{m}})^{\rm op}$. This exactness follows from the flatness of ${\rm ind}_{\overline{\mathcal U}}^{\overline{\mathcal
  S}}{\bf 1}_{{\mathfrak
  o}_{m}}$ over ${\mathcal H}(\overline{\mathcal S},\overline{\mathcal U})_{{\mathfrak
  o}_{m}}$, Lemma \ref{sl2flat}.\hfill$\Box$\\

{\bf Remark:} An alternative description of ${\mathcal V}^X_M$ (which we do not need), at least when restricted to facets of codimension $0$ and $1$, and hence of ${\mathcal V}_M$ can be given if
$M$ is realized inside a smooth $G$-representation $V$. Recall that for such $V$ the submodule $V^{I_0}$ of $I_0$-invariants is in a natural way a (left) module over ${\mathcal H}(G,I_0)$. Let $M$ be an ${\mathcal
  H}(G,I_0)$-sub module of $V^{I_0}$. For any facet $D$ of $X$ we let$${\mathcal
  V}^X_{(V,M)}(D)=\sum_{g\in G\atop gC\supset D}gM$$(sum inside $V$). This defines a $G$-equivariant coefficient system ${\mathcal
  V}^X_{(V,M)}$ on $X$. If $M\in {\rm Mod}^{\rm fin}({\mathcal
  H}(G,I_0)_{{\mathfrak o}_m})$ for some $m\ge1$ then there is a natural $G$-equivariant morphism ${\mathcal
  V}^X_{M}\to{\mathcal
  V}^X_{(V,M)}$ which, at least when restricted to facets of codimension $0$ and $1$, is an isomorphim. To see this one can use arguments from the proof of Propostition
\ref{jalgheck} (the starting point is to see, using Lemmata \ref{sl2flat} and \ref{sl2einbett}, that for any codimension-$1$-face $F$ of $C$ with corresponding subgroup
$I_0^F$ of $G$ we have $(I_0^F.M)^{I_0}=M$).

\section{$(\varphi^r,\Gamma)$-modules }

\label{phgase}

\subsection{$(\psi^r,\Gamma)$-modules and $(\varphi^r,\Gamma)$-modules}

Let ${\mathcal O}_{\mathcal E}^+={\mathfrak o}[[{\mathfrak N}_0]]$ denote the completed group ring of
${\mathfrak N}_0$ over ${\mathfrak o}$. Let ${\mathcal O}_{\mathcal E}$ denote the $p$-adic completion of the
localization of ${{{\mathcal O}_{\mathcal E}^+}}$ with respect to the
complement of $\pi_K{{{\mathcal O}_{\mathcal E}^+}}$.

Let $\varphi_{{{{\mathcal O}_{\mathcal E}^+}}}$ denote the endomorphism of
${{{\mathcal O}_{\mathcal E}^+}}$ induced by the endomorphism
$n\mapsto\varphi n\varphi^{-1}$ of ${\mathfrak N}_0$. Let $\varphi_{{\mathcal
    O}_{\mathcal E}}$ denote the endomorphism of ${\mathcal O}_{\mathcal E}$
induced from $\varphi_{{{{\mathcal O}_{\mathcal E}^+}}}$ by
functoriality. As an ${\mathfrak o}$-module, ${\mathcal O}_{\mathcal E}$ decomposes
as ${\mathcal O}_{\mathcal E}={\rm im}(\varphi_{{\mathcal O}_{\mathcal E}})\oplus({\mathfrak N}_0-{\mathfrak N}_0^p){\rm im}(\varphi_{{\mathcal O}_{\mathcal E}})$ (notice that $\varphi{\mathfrak N}_0\varphi^{-1}={\mathfrak
  N}_0^p$). We denote by $\psi_{{\mathcal O}_{\mathcal E}}$ the ${\mathfrak
  o}$-linear endomorphism of ${\mathcal
  O}_{\mathcal E}$ with $\psi_{{\mathcal O}_{\mathcal
    E}}\circ\varphi_{{\mathcal O}_{\mathcal E}}={\rm id}$ and with kernel
${\rm ker}(\psi_{{\mathcal O}_{\mathcal E}})=({\mathfrak N}_0-{\mathfrak
  N}_0^p){\rm im}(\varphi_{{\mathcal O}_{\mathcal E}})$. Similarly (or simply
by restriction) we define the ${\mathfrak
  o}$-linear endomorphism  $\psi_{{{{\mathcal O}_{\mathcal E}^+}}}$ of ${{{\mathcal O}_{\mathcal E}^+}}$. The conjugation action $(\gamma,n)\mapsto \gamma n\gamma^{-1}$ of $\Gamma$ on ${\mathfrak N}_0$ induces an
action $(\gamma,a)\mapsto \gamma\cdot a$ of $\Gamma$ on ${{{\mathcal
      O}_{\mathcal E}^+}}$ and on ${\mathcal O}_{\mathcal E}$.

On ${{k_{\mathcal E}^+}}=k[[{\mathfrak
  N}_0]]={{{\mathcal O}_{\mathcal E}^+}}\otimes_{\mathfrak o}
k$ we have the $\varphi$-operator $\varphi_{{{k_{\mathcal
        E}^+}}}=\varphi_{{{\mathcal O}_{\mathcal E}^+}}\otimes_{\mathfrak
  o}k$, the $\psi$-operator $\psi_{{{k_{\mathcal E}^+}}}=\psi_{{{\mathcal
      O}_{\mathcal E}^+}}\otimes_{\mathfrak o}k$ and the induced action of $\Gamma$, and similarly on ${{k_{\mathcal E}}}={\mathcal O}_{\mathcal E}\otimes_{{\mathfrak o}}k$.\\

Let $r\in{\mathbb N}$. We need the non-commutative
polynomial ring ${\mathcal O}_{\mathcal
  E}^+[\varphi^r_{{\mathcal O}_{\mathcal
  E}^+}]$ over ${\mathcal O}_{\mathcal
  E}^+$ in which relations are as imposed by multiplication in ${\rm
GL}_2({\mathbb Q}_p)$ (i.e. $\varphi^r_{{\mathcal O}_{\mathcal
  E}^+}\cdot[n]=[n]^{p^r}\cdot\varphi^r_{{\mathcal O}_{\mathcal
  E}^+}$ for $n\in {\mathfrak N}_0$) and the twisted group ring
${\mathcal O}_{\mathcal
  E}^+[\varphi^r_{{\mathcal O}_{\mathcal
  E}^+},\Gamma]$ over ${\mathcal O}_{\mathcal
  E}^+[\varphi^r_{{\mathcal O}_{\mathcal
  E}^+}]$, again with relations as imposed by
multiplication in ${\rm GL}_2({\mathbb Q}_p)$. (I.e. for
$\gamma\in\Gamma$ and $n\in {\mathfrak N}_0$ we have $\varphi^r_{{\mathcal O}_{\mathcal
  E}^+}\cdot\gamma=\gamma\cdot\varphi^r_{{\mathcal O}_{\mathcal
  E}^+}$ and $\gamma\cdot[n]=[\gamma n\gamma^{-1}]\cdot\gamma$. Specifically,
for $a\in{\mathbb Z}_p^{\times}$, setting $t=[\nu]-1$ we have $\gamma(a)\cdot
t=((t+1)^a-1)\cdot\gamma(a)$.) Similarly we define ${\mathcal O}_{\mathcal
  E}[\varphi^r_{{\mathcal O}_{\mathcal
  E}}]$ and ${\mathcal O}_{\mathcal
  E}[\varphi^r_{{\mathcal O}_{\mathcal
  E}},\Gamma]$, as well as ${{k_{\mathcal E}^+}}[\varphi^r_{{{k_{\mathcal
        E}^+}}}]$ and ${{k_{\mathcal E}^+}}[\varphi^r_{{{k_{\mathcal
        E}^+}}},\Gamma]$ (equivalently, ${\mathcal O}_{\mathcal
  E}[\varphi^r_{{\mathcal O}_{\mathcal
  E}}]={\mathcal O}_{\mathcal
  E}^+[\varphi^r_{{\mathcal O}_{\mathcal
  E}^+}]\otimes_{{\mathcal O}_{\mathcal
  E}^+}   {\mathcal O}_{\mathcal
  E}$ and ${\mathcal O}_{\mathcal
  E}[\varphi^r_{{\mathcal O}_{\mathcal
  E}},\Gamma]={\mathcal O}_{\mathcal
  E}^+[\varphi^r_{{\mathcal O}_{\mathcal
  E}^+},\Gamma]\otimes_{{\mathcal O}_{\mathcal
  E}^+}   {\mathcal O}_{\mathcal
  E}$ as well as ${{k_{\mathcal E}^+}}[\varphi^r_{{{k_{\mathcal
        E}^+}}}]={\mathcal O}_{\mathcal
  E}^+[\varphi^r_{{\mathcal O}_{\mathcal
  E}^+}]\otimes_{{\mathfrak o}}k$ and ${{k_{\mathcal E}^+}}[\varphi^r_{{{k_{\mathcal
        E}^+}}},\Gamma]={\mathcal O}_{\mathcal
  E}^+[\varphi^r_{{\mathcal O}_{\mathcal
  E}^+},\Gamma]\otimes_{{\mathfrak o}}k$). We will often just write
${{k_{\mathcal E}^+}}[\varphi^r]$ and ${{k_{\mathcal
      E}^+}}[\varphi^r,\Gamma]$, and also drop the exponent $r$ in case $r=1$.\\

{\bf Definition:} (a) An \'{e}tale $(\varphi^r,\Gamma)$-module over
${\mathcal O}_{\mathcal E}$ is a module ${\bf D}$ over ${\mathcal O}_{\mathcal
  E}[\varphi^r_{{\mathcal O}_{\mathcal
  E}},\Gamma]$, finitely generated as a module over ${\mathcal O}_{\mathcal
  E}$ and such that the structure map $\varphi^r_{{\bf D}}$ is \'{e}tale,
i.e. $\varphi^r_{{\bf D}}$ is injective and satisfies ${\bf D}=\oplus_{n\in{\mathfrak
  N}_0/{\mathfrak N}_0^{p^r}}n\varphi^r_{{\bf D}}({\bf D})$. In this
situation there is a unique ${\mathfrak o}$-linear endomorphism $\psi^r_{{\bf
    D}}$ of ${\bf D}$ with $\psi^r_{{\bf
    D}}\circ\varphi^r_{{\bf
    D}}={\rm id}_{\bf D}$ and with kernel $\oplus_{n\in({\mathfrak
  N}_0/{\mathfrak N}_0^{p^r})-\{{\mathfrak N}_0^{p^r}\}}n\varphi^r_{{\bf
  D}}({\bf D})$. A $(\varphi^r,\Gamma)$-module over ${{k_{\mathcal E}}}$ is an \'{e}tale
$(\varphi^r,\Gamma)$-module over ${\mathcal O}_{{\mathcal E}}$ whose ${\mathcal
  O}_{{\mathcal E}}$-action factors through the quotient ${{k_{\mathcal E}}}$ of ${\mathcal
  O}_{{\mathcal E}}$.

(b) Similarly we define \'{e}tale $(\varphi^r,\Gamma_0)$-modules over
${\mathcal O}_{\mathcal E}$ and over ${{k_{\mathcal E}}}$.\\

{\bf Definition:} A $(\psi^r,\Gamma)$-module over ${{k_{\mathcal E}^+}}$ is a
finitely generated free
  ${{k_{\mathcal E}^+}}$-module $D^{\sharp}$, together with a $k$-linear
  endomorphism $\psi_{D^{\sharp}}^r$ satisfying
  $\psi_{D^{\sharp}}^r(\varphi^r_{{{k_{\mathcal
          E}^+}}}(\alpha)x)=\alpha\psi_{D^{\sharp}}^r(x)$ and a continuous
  semilinear action of $\Gamma$ that commutes with $\psi_{D^{\sharp}}^r$.

It is called non degenerate if ${\rm ker}(\psi_{D^{\sharp}}^r)$ contains no non zero
${k_{\mathcal E}^+}$-sub module.\\

$\bullet$ Let $D^{\sharp}$ be a $(\psi^r,\Gamma)$-module over ${{k_{\mathcal
      E}^+}}$. Viewing $D^{\sharp}$ as a linearly compact $k$-vector space we
endow the topological dual $(D^{\sharp})^*={\rm Hom}_{k}^{\rm
  ct}(D^{\sharp},k)$ with the structure of a ${{k_{\mathcal E}^+}}$-module by setting
$(\alpha\cdot\ell)(x)=\ell({\alpha}\cdot x)$ for $\ell\in (D^{\sharp})^*$, $x\in D^{\sharp}$
and $\alpha\in {{k_{\mathcal E}^+}}$.\footnote{We choose not to use the involution on ${{k_{\mathcal E}^+}}$ induced by inversion in ${\mathfrak
    N}_0$ when defining the ${{k_{\mathcal E}^+}}$-action on the dual, taking advantage of the commutativity of ${\mathfrak
    N}_0$.} Notice that $(D^{\sharp})^*$ is a torsion ${{k_{\mathcal E}^+}}$-module. We define a $k$-linear endomorphism
  $\varphi^r_{(D^{\sharp})^*}$ on $(D^{\sharp})^*$ through
  $(\varphi^r_{(D^{\sharp})^*}(\ell))(x)=\ell(\psi^r_{D^{\sharp}}(x))$ and an action by $\Gamma$ through
  $(\gamma(\ell))(x)=\ell(\gamma^{-1}(x))$. One checks that this defines on
  $(D^{\sharp})^*$ the structure of a module over ${{k_{\mathcal
        E}^+}}[\varphi^r,\Gamma]$.

$\bullet$ Let $D^{\sharp}$ be a $(\psi^r,\Gamma)$-module over ${{k_{\mathcal
      E}^+}}$ which is non degenerate and such that $\psi^r_{D^{\sharp}}$ is
surjective. Then ${\bf D}=D^{\sharp}\otimes_{k_{\mathcal E}^+}{k_{\mathcal E}}$
carries a unique structure of a $(\varphi^r,\Gamma)$-module over
${{k_{\mathcal E}}}$ compatible with the $(\psi^r,\Gamma)$-structure on
$D^{\sharp}$ (i.e. the $\Gamma$-actions coincide, and $\psi^r_{\bf D}$ extends $\psi^r_{D^{\sharp}}$). For $r=1$ this is well known [for the uniqueness one may use
Proposition II.3.4(ii) of \cite{col}; see also \cite{vien} Proposition
3.3.24.]. For general $r\ge1$ the proofs are the same. (For the $(\psi^r,\Gamma)$-modules met later on in this paper, the corresponding $(\varphi^r,\Gamma)$-modules over
${{k_{\mathcal E}}}$ come along simultaneously and explicitly, see section \ref{secdfu}.)\\

Let ${\rm Gal}_{{\mathbb Q}_p}={\rm Gal}({\mathbb Q}_p^{\rm alg}/{\mathbb Q}_p)$ denote the absolute Galois group of ${\mathbb Q}_p$.

\begin{satz}\label{fofu} (Fontaine) There is an equivalence functor ${\bf D}\mapsto W({\bf D})$ from the category of $(\varphi,\Gamma)$-modules over ${{k_{\mathcal E}}}$ to the category of finite dimensional (smooth) $k$-representations of ${\rm Gal}_{{\mathbb Q}_p}$.
\end{satz}

{\sc Proof:} This is shown in \cite{fon}, for a brief account see \cite{berger} Theorem 2.1.2.\hfill$\Box$\\

Let ${\mathcal
  I}_{{\mathbb Q}_p}$ denote the inertia subgroup of ${\rm Gal}_{{\mathbb Q}_p}$. For
$m\ge0$ let
$\omega_{m+1}:{\mathcal I}_{{\mathbb Q}_p}\to\overline{\mathbb F}_p^{\times}$
denote the fundamental character of level $m+1$. It is given by the formula
$\omega_{m+1}(g)=\overline{g(\pi_{m+1})/\pi_{m+1}}$ where
$\pi_{m+1}\in{\mathbb Q}_p^{\rm alg}$ is chosen such that
$\pi_{m+1}^{p^{m+1}-1}=-p$. We denote by $\omega$ the cyclotomic character modulo $p$ of ${\rm Gal}_{{\mathbb Q}_p}$. 

For $0\le h\le p^{m+1}-1$ let ${\rm ind}(\omega_{m+1}^h)$ be the $(m+1)$-dimensional ${\rm Gal}_{{\mathbb Q}_p}$-representation over $k$ with ${\rm det}({\rm ind}(\omega_{m+1}^h))=\omega^h$ and ${\rm ind}(\omega_{m+1}^h)|_{{\mathcal
  I}_{{\mathbb Q}_p}}=\oplus_{j=0}^m\omega_{m+1}^{p^jh}$ described in \cite{berger}. For
$\beta\in({k}^{\rm alg})^{\times}$ let $\mu_{\beta}$ denote the unique unramified (i.e. trivial on ${\mathcal I}_{{\mathbb Q}_p}$) character of ${\rm Gal}_{{\mathbb Q}_p}$ sending the geometric Frobenius (i.e. lifting $x\mapsto x^{-p}$) to $\beta$. If $\beta^{m+1}\in k^{\times}$ then ${\rm
  ind}(\omega_{m+1}^h)\otimes \mu_{\beta}$ is defined over $k$.   

One says that an integer $1\le h\le p^{m+1}-2$ is primitive (with respect to
$m+1$) if there is no $n<m+1$ dividing $m+1$ such that $h$ is a multiple of
$(p^{m+1}-1)/(p^n-1)$. If $h$ is primitive then ${\rm
  ind}(\omega_{m+1}^h)$ is absolutely irreducible. For all this
see \cite{berger} section 2.1.

\subsection{Standard cyclic $k_{\mathcal E}^+[\varphi^r,\Gamma]$-modules}

We put $t=[\nu]-1\in {{k_{\mathcal E}^+}}=k[[{\mathfrak
  N}_0]]$. This is a uniformizer in the complete discrete valuation ring ${{k_{\mathcal E}^+}}$, thus ${{k_{\mathcal E}^+}}=k[[t]]$ and ${{k_{\mathcal E}}}=k((t))$.

We often identify elements of ${\mathbb F}_p^{\times}$ with their Teichm\"uller lifting in ${\mathbb Z}_p^{\times}$. In particular, for $x\in{\mathbb F}_p^{\times}$ we may consider the element $\gamma(x)\in\Gamma$.\\

{\bf Definition:} (a) We say that a $k_{\mathcal E}^+[\varphi^r]$-module $H$
is standard cyclic of perimeter $m+1\in{\mathbb N}$ if it is generated by ${\rm
  ker}(t|_H)=H^{{\mathfrak N}_0}$, if it is a torsion
${{k_{\mathcal E}^+}}$-module and if there are a $k$-basis
$e_0,\ldots,e_m$ of ${\rm ker}(t|_H)$, integers
$0\le k_0,\ldots,k_m\le p^r-1$ and units $\varrho_0,\ldots,\varrho_m\in
k^{\times}$ with $$t^{k_i}\varphi^r e_{i-1}=\varrho_ie_i$$ for all $0\le i\le m$. Here and below we extend the indexing of the $k_i$, $e_i$, $\varrho_i$ by
$i\in\{0,\ldots,m\}$ to an indexing by $i\in{\mathbb Z}$ such that
$k_i=k_{i+m+1}$, $e_i=e_{i+m+1}$, $\varrho_i=\varrho_{i+m+1}$.

(b) We say that a $k_{\mathcal E}^+[\varphi^r,\Gamma]$-module $H$
is standard cyclic of perimeter $m+1\in{\mathbb N}$ if it is standard cyclic
of perimeter $m+1$ as a $k_{\mathcal E}^+[\varphi^r]$-module, in such a way that all
the above $e_i$ can be chosen as eigenvectors for the action of $\Gamma$, with
eigenvalues in $k^{\times}$.

(c) Giving $H$ the discrete topology, we endow the topological dual $H^*={\rm Hom}_{k}^{\rm
  ct}(H,k)$ with the structure of a ${{k_{\mathcal E}^+}}$-module by setting
$(\alpha\cdot\ell)(x)=\ell({\alpha}\cdot x)$ for $\ell\in H^*$, $x\in H$
and $\alpha\in {{k_{\mathcal E}^+}}$. We define a $k$-linear endomorphism
  $\psi^r_{H^*}$ on $H^*$ through
  $(\psi^r_{H^*}(\ell))(x)=\ell(\varphi^r(x))$ and an action by $\Gamma$ through
  $(\gamma(\ell))(x)=\ell(\gamma^{-1}(x))$.\\

\begin{pro}\label{stnle} (a) If $k_i>0$ for at least one $0\le i\le m$ then
  $H^*$ is a non degenerate $(\psi^r,\Gamma)$-module over ${{k_{\mathcal
        E}^+}}$ with surjective operator $\psi^r_{H^*}$, free of
  rank $m+1$ as a $k_{\mathcal E}^+$-module. 

(b) Suppose that for any $1\le j\le m$ there is some $0\le i\le
  m$ with $k_i{{{}}}\ne k_{i+j}{{{}}}$. Then $H$ is irreducible as a
  ${{k_{\mathcal E}^+}}[\varphi^r]$-module.

(c) For $0\le i\le m$ let $\eta_i:\Gamma\to k^{\times}$ be the character
with $\gamma e_i=\eta_i(\gamma)e_i$ for all $\gamma\in\Gamma$. Suppose that for any $1\le j\le m$ which satisfies $k_i=k_{i+j}$
for all $0\le i\le m$ there is some $0\le i\le m$ with $\eta_i\ne\eta_{i+j}$. Then $H$ is irreducible as a
  ${{k_{\mathcal E}^+}}[\varphi^r,\Gamma]$-module.   

(d) If $k_i=0$ for all $0\le i\le m$ then $t=0$ on $H$, and $e_0,\ldots,e_m$ is a
  $k$-basis of $H$.
\end{pro}

{\sc Proof:} Statement (d) is clear. We prove statement (a).

{\it Claim:} The action of $t$ on $H$ is surjective.

For any $x\in H$ we need to find some $y\in H$ with $ty=x$. As ${\rm
  ker}(t|_H)$ generates $H$ we may assume, by additivity, that
  $x=t^{n_1}\varphi^{rn_1'}\cdots t^{n_l}\varphi^{rn_l'}e_i$ for some $0\le i\le
  d$, some $l\ge0$, some $n_j, n'_j\ge0$. Using $\varphi^r t=t^{p^r}\varphi^r$ we rewrite this as
  $x=t^{n}\varphi^{rn'}e_i$ (some $n,n'\ge0$). If $n\ge 1$ we are
  done. Otherwise we put
  $w_i=k_i+p^rk_{i-1}+\ldots+p^{ri}k_0+p^{r(i+1)}k_{m}+\ldots+p^{rm}k_{i+1}$ and substitute
  $t^{w_i{{{}}}}\varphi^{r(m+1)} e_{i}$ for $e_i$ --- up to a scalar in
  $k^{\times}$ these are the same, as follows from $\varphi^r t=t^{p^r}\varphi^r$. We use $\varphi^r t=t^{p^r}\varphi^r$ again to rewrite the
  result as $x=t^{n}\varphi^{rn'}e_{i}$ for some new $n,n'\ge0$. Now we have $n\ge 1$
  because $w_i\ge 1$ for any $i$, as follows from the hypothesis that $k_i>0$ for at least one $i$. The claim is proven.

Let $\ell\in H^*$, $\ell\ne0$. Choose $x\in H$ with
$\ell(x)\ne0$. By the above claim (iterated $n$ times) we find for any $n\ge0$ some $y\in H$ with ${t}^ny=x$, hence
  $(t^n\ell)(y)=\ell({t}^ny)=\ell(x)\ne0$, hence $t^n\ell\ne0$. It follows
  that $H^*$ is a torsion free ${{k_{\mathcal E}^+}}$-module. On the other
  hand, we know $H^*\otimes_{{{k_{\mathcal E}^+}}}k$ (with $t\mapsto0$)
  because this is dual to ${\rm ker}(t|_H)$. In view of Nakayama's Lemma we
  obtain that $H^*$ is ${{k_{\mathcal E}^+}}$-free of rank $m+1$.

Now $\varphi^r$ is injective on $H$ [as $\varphi^r t=t^{p^r}\varphi^r$ we have $t\cdot {\rm ker}(\varphi^r)\subset{\rm ker}(\varphi^r)$, thus if we had ${\rm ker}(\varphi^r)\ne0$ then also ${\rm ker}(\varphi^r)\cap{\rm ker}(t)\ne0$, but this is false], hence $\psi^r_{H^*}$ is surjective; similarly, ${\rm ker}(\psi^r_{H^*})$ contains no non zero
${k_{\mathcal E}^+}$-sub module.

(b) Let $0\ne
Z\subset H$ be a non zero
${{k_{\mathcal E}^+}}[\varphi^r]$-sub module. As $H$ is a torsion
${{k_{\mathcal E}^+}}$-module so is $Z$, hence ${\rm ker}(t|_Z)$ is non
zero. Let us put $$\eta(z)=\eta(\sum_{0\le i\le m}x_ie_i)={\rm max}\{k_i{{{}}}\,|\,0\le i\le m,\,\, x_i\ne0\}$$for non zero elements $z=\sum_{0\le i\le
  m}x_ie_i\in {\rm ker}(t|_Z)$ (with $x_i\in k$). For such $0\ne z\in {\rm ker}(t|_Z)$ let us put
$\Lambda(z)=t^{\eta(z)}\varphi^r z$. This is again an element in ${\rm ker}(t|_Z)$ (as $Z$ is stable
under $t$ and ${\varphi}^r$); moreover $\Lambda(z)\ne0$. We apply $\Lambda$ repeatedly: The hypothesis shows that for sufficiently large $n\ge0$ we have $\Lambda^n(z)\in k^{\times}\cdot
  e_i$ for some $0\le i\le m$. But then we further see that
  $\Lambda^{n+j}(z)\in k^{\times}\cdot e_{i+j}$ for all $j\ge0$. It follows
  that $Z$ contains all the $e_i$, hence $Z=H$ as the $e_i$ generate the ${{k_{\mathcal E}^+}}[\varphi^r]$-module $H$.

(c) Let $0\ne
Z\subset H$ be a non zero
${{k_{\mathcal E}^+}}[\varphi^r,\Gamma]$-sub module. If $H$ is not already irreducible as a
  ${{k_{\mathcal E}^+}}[\varphi^r]$-module then, by the proof of (b), there is some $0\ne z\in {\rm ker}(t|_Z)$ such that for all $n\ge0$, if we write $\Lambda^n(z)=\sum_{0\le i\le
  m}x_{i,n}e_i$, then the number $|\{i\,|\,x_{i,n}\ne0\}|$ is larger than $1$
and independent on $n$. Thus $x_{i,0}\ne0$ and $x_{i+j,0}\ne0$ for some $i,j$,
and this $j$ violates the hypothesis in (b). By the hypothesis in (c), replacing $z$ by $\Lambda^n(z)$
for a suitable $n$ (and replacing $i$ by $i+n$) we may assume that
$\eta_i\ne\eta_{i+j}$. This allows us to produce a non-zero element $\sum_{0\le i\le
  m}y_ie_i$ in ${\rm ker}(t|_Z)$ such that $|\{i\,|\,y_{i}\ne0\}|<|\{i\,|\,x_{i,0}\ne0\}|$. Proceeding by induction we obtain that $e_i\in {\rm ker}(t|_Z)$ for some $0\le i\le m$ and hence, as in (b), that $Z=H$.\hfill$\Box$\\

For $0\le j\le m$ let $i_j=p^r-1-k_{m+1-j}$ and $$w_j=k_i+p^rk_{j-1}+\ldots+p^{rj}k_0+p^{r(j+1)}k_{m}+\ldots+p^{rm}k_{j+1}.$$Let $h_0=0$ and for $1\le j\le m+1$ let $h_j=\sum_{i=0}^{j-1}i_{m+i+1-j}p^{ri}$. Let $\varrho=\prod_{i=0}^m\varrho_i$; then
$t^{w_i}\varphi^{r(m+1)}e_i=\varrho e_i$ for all $0\le i\le m$ (use $\varphi^r t=t^{p^r}\varphi^r$). 

\begin{lem}\label{normint} (a) There is some $0\le s\le
p-2$ such that $\gamma(x)e_i=x^{-h_i-s}e_i$ for all $x\in{\mathbb F}_p^{\times}$,
all $0\le i\le m$.  

(b) $h=h_{m+1}/(p-1)$ is an integer.
\end{lem}

{\sc Proof:} (a) For $a\in{\mathbb N}\cap {\mathbb Z}_p^{\times}$ we formally compute
$[\nu]^a-1=\sum_{j=1}^a{a \choose j}([\nu]-1)^j=\sum_{j=1}^a{a \choose j}t^j$,
thus $\gamma(a)t\gamma(a^{-1})-at=([\nu]^a-1)-at\in t^2k_{\mathcal
  E}^+$ and therefore also \begin{gather}\gamma(a)t^k\gamma(a^{-1})-a^kt^k\in t^{k+1}k_{\mathcal
  E}^+\label{naly}\end{gather} for all $k\in{\mathbb Z}_{\ge0}$. For $0\le i\le m$ let
$n_i\in{\mathbb Z}$ such that $\gamma(a)e_i=a^{n_i}e_i$ for all $a\in{\mathbb
  N}\cap {\mathbb Z}_p^{\times}$. Now let $1\le i\le
m$. Then$$\varrho_ia^{n_{i}}e_i=\varrho_i\gamma(a)e_i=\gamma(a)\varrho_ie_i=\gamma(a)t^{k_i}\varphi^re_{i-1}\stackrel{(i)}{=}t^{k_i}\varphi^ra^{k_i}\gamma(a)e_{i-1}=t^{k_i}\varphi^ra^{k_i+n_{i-1}}e_{i-1}$$where
in $(i)$ we used formula (\ref{naly}) and the fact that
$t^{k_i}\varphi^re_{i-1}$ belongs to ${\rm ker}(t|_H)$ (as it equals
$\varrho_ie_i$). Comparing this with
$\varrho_ie_i=t^{k_i}\varphi^re_{i-1}$ we obtain $k_i\equiv n_i-n_{i-1}$ modulo $(p-1)\mathbb Z$. On the
other hand, the definition of $k_i$ shows $k_i\equiv -i_{m+1-i}\equiv
h_{i-1}-h_i$ modulo $(p-1)\mathbb Z$. Together this means $n_i+h_i\equiv
n_{i-1}+h_{i-1}$ modulo $(p-1)\mathbb Z$ and statement (a) follows.

(b) The fraction $h$ is an integer if and only if $p-1$ divides $h_{m+1}$, if
and only if $p-1$  divides $\sum_{j=0}^mi_j$, if and only if $p-1$ divides
$\sum_{j=0}^mk_j$, if and only if $p-1$ divides $w_i$ for any $0\le i\le m$. For any such $i$ and any $\gamma\in
\Gamma$ we have$$\gamma
t^{w_i}\gamma^{-1}\varphi^{r(m+1)}\gamma(e_i)=\gamma(t^{w_i}\varphi^{r(m+1)}e_i)=\varrho\gamma(e_i)=t^{w_i}\varphi^{r(m+1)}\gamma(e_i)$$and
this is a non zero element in $H$. Specifically, taking $\gamma=\gamma(a)$ and
using formula (\ref{naly}) we
obtain $a^{w_i}=1$ in ${\mathbb F}_p^{\times}$ for all $a\in{\mathbb N}\cap{\mathbb Z}_p^{\times}$. This
implies that $p-1$ divides $w_i$.\hfill$\Box$\\ 

\begin{lem}\label{maerz} Suppose that $k_i>0$ for at least one $i$. The $(\varphi^r,\Gamma)$-module ${\bf
    D}$ over $k_{\mathcal E}$ associated with
  the $(\psi^r,\Gamma)$-module $H^*$ admits a $k_{\mathcal E}$-basis
  $g_0,\ldots,g_m$ such that \begin{align}\varphi^r_{{\bf
    D}}(g_j)&=g_{j+1}\quad\quad\quad\quad\mbox{ for }0\le j\le
  m-1,\label{malaur1}\\\varphi^r_{{\bf
    D}}(g_m)&=\varrho^{-1}t^{-h_{m+1}}g_0,\label{malaur2}\\\gamma(x)(g_j)&- x^{s}g_{j}\in
    t\cdot k_{\mathcal E}^+\cdot g_j \quad\quad\mbox{ for }0\le j\le m\mbox{
      and }x\in{\mathbb F}_p^{\times}.\label{malaur3}\end{align}These formulae completely
  characterize the actions of $\varphi^r_{{\bf
    D}}$ and $\Gamma$.
\end{lem} 

{\sc Proof:} Let ${\bf
    D}$ denote the $\varphi^r$-module over $k_{\mathcal E}$ with $k_{\mathcal E}$-basis
  $g_0,\ldots,g_m$ and with $\varphi^r$-operator $\varphi^r_{{\bf
    D}}$ given by formulae (\ref{malaur1}) and (\ref{malaur2}). It is
\'{e}tale and hence admits the usual canonical left inverse $\psi^r_{{\bf
    D}}$. For $0\le j\le m$ let $f_j=t^{h_j}g_j$. For $1\le j\le m$ and $\alpha\in k_{\mathcal E}^+$ we compute\begin{align}\psi_{{\bf
    D}}^r(\alpha f_j)=t^{h_{j-1}}\psi_{{\bf
    D}}^r(\alpha t^{i_{m+1-j}}g_j)=t^{h_{j-1}}\psi^r_{k_{\mathcal E}^+}(\alpha
t^{i_{m+1-j}})g_{j-1}=\psi^r_{k_{\mathcal E}^+}(\alpha
t^{i_{m+1-j}})f_{j-1}\label{laur1}\end{align}where in the first equation we used
$h_j=p^rh_{j-1}+i_{m+1-j}$. We also have $h_{m+1}=p^rh_{m}+i_{0}$ and
therefore we similarly compute\begin{align}\psi_{{\bf
    D}}^r(\alpha f_0)=\psi_{{\bf
    D}}^r(\alpha g_0)=\varrho \psi_{{\bf
    D}}^r(\alpha t^{h_{m+1}}\varphi^r_{{\bf D}}(g_m))&=\varrho t^{h_m}
\psi_{k_{\mathcal E}^+}^r(\alpha t^{i_0})g_m\notag\\{}&=\varrho\psi_{k_{\mathcal E}^+}^r(\alpha t^{i_0})f_m.\label{laur2}\end{align}Let $D^{\sharp}$ be the free $k_{\mathcal E}^+$-submodule of ${\bf
    D}$ with basis $f_0,\ldots,f_m$. Formulae (\ref{laur1}) and (\ref{laur2})
  show that $D^{\sharp}$ is stable under $\psi_{{\bf
    D}}^r$, hence that $D^{\sharp}$ is a $\psi^r$-module over $k_{\mathcal E}$. Its
$\psi^r$-operator $\psi^r_{D^{\sharp}}= \psi^r_{{\bf D}}|_{D^{\sharp}}$ is
surjective, as follows from formulae (\ref{laur1}) and (\ref{laur2}) in view
of\begin{gather}\psi^r_{{{k_{\mathcal
          E}^+}}}(t^{\sum_{i=0}^rn_ip^i})=(-1)^{\sum_{i=0}^{r-1}n_i}t^{n_r}\quad\quad\mbox{ for }0\le n_0,\ldots,n_{r-1}\le p-1\mbox{ and }0\le n_r.\label{epsifo}\end{gather}For $0\le i\le m$ we define $e'_i\in (D^{\sharp})^*$ by $e'_i(f_j)=\delta_{ij}$ and by $e'_i|_{t\cdot
  D^{\sharp}}=0$. The set $\{e'_0,\ldots,e'_m\}$ is a $k$-basis of ${\rm
  ker}(t|_{(D^{\sharp})^*})$, and it generates $(D^{\sharp})^*$ as a ${{k_{\mathcal
      E}^+}}[\varphi^r]$-module as will follow from formulae (\ref{bermonouni1})
and (\ref{bermonouni2}) below. We claim\begin{align}t^{k_i{{{}}}}\varphi^r_{(D^{\sharp})^*}
  e'_{i-1}&=e'_i&\mbox{ for }1\le i\le m,\label{bermonouni1}\\t^{k_0{{{}}}}\varphi^r_{(D^{\sharp})^*}e'_{m}&=\varrho e'_0,&{}\label{bermonouni2}\\t^{w_i}\varphi_{(D^{\sharp})^*}^{r(m+1)}e'_i&=\varrho e'_i&\mbox{ for }0\le i\le m.\label{beridemmonouni}\end{align}For $n\ge0$ we compute\begin{align}(t^{k_i{{{}}}}\varphi^r_{(D^{\sharp})^*}
  e'_{i-1})(t^nf_j)&=(\varphi^r_{(D^{\sharp})^*}e'_{i-1})(t^{k_i+n}f_j)\notag\\{}&=e'_{i-1}(\psi^r_{D^{\sharp}}(t^{k_i+n}f_j)).\notag\end{align}Inserting formula (\ref{laur1}) we thus get \begin{gather}(t^{k_i{{{}}}}\varphi^r_{(D^{\sharp})^*}
  e'_{i-1})(t^nf_j)=e'_{i-1}(\psi^r_{{{k_{\mathcal E}^+}}}(t^{k_i+n+p^r-1-k_j})f_{j-1})\label{endber1}\end{gather}for $1\le
j\le m$ while for $j=0\equiv m+1$, inserting formula (\ref{laur2}) we get\begin{gather}(t^{k_i{{{}}}}\varphi^r_{(D^{\sharp})^*}
  e'_{i-1})(t^nf_0)=e'_{i-1}(\varrho\psi^r_{{{k_{\mathcal E}^+}}}(t^{k_i+n+p^r-1-k_0})f_{m}).\label{endber2}\end{gather}The right hand
side in formula (\ref{endber1}) (resp. in formula (\ref{endber2})) vanishes if
$i\ne j$ (resp. if $i\ne m+1\equiv 0$). To evaluate the right hand
side if $1\le i=j\le m$ (resp. if $i=m+1\equiv 0$) we again use formula (\ref{epsifo}). It
shows that for $n>0$ the right hand side in formula (\ref{endber1}) (resp. in formula (\ref{endber2}))
vanishes and that for $n=0$ and $1\le i=j\le m$ its value is $1$ (resp. for
$n=0$ and $i=m+1\equiv 0$ its value is $\varrho$). We have proven formulae (\ref{bermonouni1}) and (\ref{bermonouni2}). Formula
     (\ref{beridemmonouni}) follows by iteration (as $\varphi^r
     t=t^{p^r}\varphi^r$).

As both $H^*$ and $D^{\sharp}$ are ${{k_{\mathcal
      E}^+}}$-free of the same rank, a comparison of formulae (\ref{bermonouni1}), (\ref{bermonouni2}) with those
describing the action of
$\varphi^r$ on $H\cong (H^*)^*$ shows that there is an isomorphism of $k_{\mathcal E}^+[\varphi^r]$-modules
$H\cong(D^{\sharp})^*$ sending the $e_i$ to suitable $k^{\times}$-rescalings
of the $e'_i$.

We use this isomorphism to transport the $\Gamma$-action on $H$ to a
$\Gamma$-action on $(D^{\sharp})^*$. For $0\le i,j\le m$ and $x\in {\mathbb
  F}_p$ it satisfies$$e'_i(\gamma(x)(f_j)-x^{h_j+s}f_j)=(\gamma(x^{-1})e'_i)(f_j)-x^{h_j+s}e'_i(f_j).$$We
claim that this
vanishes. Indeed, if $i=j$ then this follows from
$\gamma(x^{-1})e'_i=x^{h_i+s} e'_i$ (the definition of $s$), whereas if $i\ne
j$ even both summands vanish individually. The claim proven we
infer $$\gamma(x)(f_j)-x^{h_j+s}f_j\in\bigcap_{i}{\rm ker}(e'_i)=t\cdot
  D^{\sharp}.$$On the other hand, formula (\ref{naly}) says $\gamma(x)(f_j)-x^{h_j}t^{h_j}\gamma(x)(g_j)\in t^{h_j+1}k_{{\mathcal E}}^+\cdot g_j\subset t\cdot D^{\sharp}$. Together we obtain$$\gamma(x)(g_j)-x^{s}g_j\in t^{1-h_j}\cdot
  D^{\sharp}.$$As both $\gamma(x)(g_j)$ and $x^{s}g_j$
  belong to $k_{\mathcal E}^+\cdot g_j$ and as  $k_{\mathcal E}^+\cdot g_j\cap t^{1-h_j}\cdot
  D^{\sharp}=t\cdot k_{\mathcal E}^+\cdot g_j$ we have proven formula (\ref{malaur3}).

For the final statement it is enough to show that formulae (\ref{malaur1}), (\ref{malaur2}) and
(\ref{malaur3}) characterize the $k_{\mathcal E}^+[\varphi^r,\Gamma]$-module
$(D^{\sharp})^*$. We have already seen that they imply formulae
(\ref{bermonouni1}), (\ref{bermonouni2}), (\ref{beridemmonouni}), i.e. they
characterize the action of $\varphi^r_{(D^{\sharp})^*}$. It therefore remains to
see that they characterize the action of $\Gamma$ on the $k_{\mathcal
  E}^+[\varphi^r]$-generators $e'_i$ of $(D^{\sharp})^*$. We claim$$\gamma(x)e'_i=x^{-h_i-s} e'_i.$$Indeed, both sides vanish on
$t\cdot D^{\sharp}$. To compare their values on an argument in $D^{\sharp}$ we
may ignore summands belonging to $t\cdot D^{\sharp}$. Thus formula
(\ref{malaur3}) gives
us\begin{align}(\gamma(x)e'_i)(f_j)&=e'_i(\gamma(x)^{-1}f_j)\notag\\&=e'_i(x^{-h_j-s}f_j)\notag\end{align}for any
$0\le i,j\le m$. We are done.\hfill$\Box$\\

 Let $\beta\in(k^{\rm alg})^{\times}$ such that $(-1)^m\beta^{-m-1}=\varrho=\prod_{i=0}^m\varrho_i$.

\begin{pro}\label{galstacom} (Berger) Suppose that $r=1$ and that $k_i>0$ for at least one $i$. Let ${\bf
    D}$ be the $(\varphi,\Gamma)$-module over $k_{\mathcal E}$ associated with
  the $(\psi,\Gamma)$-module $H^*$. We have an isomorphism of ${\rm Gal}_{{\mathbb
      Q}_p}$-representations $$W({\bf
    D})\cong {\rm ind}(\omega_{m+1}^h)\otimes\omega^s\mu_{\beta}.$$
\end{pro} 

{\sc Proof:} In \cite{berger} section 2.2 it is shown that the
$(\varphi,\Gamma)$-module ${\bf
    D}'$ over $k_{\mathcal E}$ with $W({\bf
    D}')\cong {\rm ind}(\omega_{m+1}^h)\otimes\omega^s\mu_{\beta}$ admits a
  basis in which the actions of $\varphi$ and $\Gamma$ satisfy
  the formulae (\ref{malaur1}), (\ref{malaur2}), (\ref{malaur3}), hence we conclude with Lemma \ref{maerz}.\hfill$\Box$\\

{\bf Remarks:} (a) In \cite{berger} section 2.2 even the precise formula for
$\gamma(g_j)$, for $\gamma\in\Gamma$, is worked out, sharpening formula (\ref{malaur3}).

(b) The results in \cite{berger} are in fact stated there only
under the hypothesis that $h$ be primitive. However, it is easily checked that
those statements of \cite{berger} which we are using in Proposition \ref{galstacom} hold true without that
primitivity assumption.

\subsection{$(\varphi^r,\Gamma)$-modules and $(\varphi,\Gamma)$-modules}

Here we briefly explain the interest in \'{e}tale $(\varphi^r,\Gamma)$-modules for {\it any} $r\in{\mathbb N}$ (we will not need this later on in the present paper): There is an
exact
functor from the category of \'{e}tale $(\varphi^r,\Gamma)$-modules to the category of \'{e}tale
$(\varphi,\Gamma)$-modules (the rank gets multiplied by the factor $r$). To the latter, of course, e.g. Theorem \ref{fofu} applies.\\

Let ${\bf D}=({\bf D},{\varphi}^r_{{\bf D}})$ be an \'{e}tale $\varphi^r$-module over ${\mathcal
  O}_{\mathcal E}$. For $0\le i\le r-1$ let ${\bf D}^{(i)}={\bf D}$ be a copy of ${\bf D}$. For
$1\le i\le r-1$ define
$\varphi_{\widetilde{\bf D}}:{\bf D}^{(i)}\to{\bf D}^{(i-1)}$ to be the
identity map on ${\bf D}$, and define $\varphi_{\widetilde{\bf D}}:{\bf
  D}^{(0)}\to{\bf D}^{(r-1)}$ to be the structure map ${\varphi}^r_{{\bf D}}$
on ${\bf D}$. Together we obtain a ${\mathbb Z}_p$-linear endomorphism
$\varphi_{\widetilde{\bf D}}$ on $$\widetilde{\bf D}=\bigoplus_{i=0}^{r-1}{\bf
  D}^{(i)}.$$Define an ${\mathcal
  O}_{\mathcal E}$-action on $\widetilde{\bf D}$ by the
formula\begin{gather}x\cdot((d_i)_{0\le i\le r-1})=(\varphi^i_{{\mathcal
  O}_{\mathcal E}}(x)d_i)_{0\le i\le r-1}.\label{oede}\end{gather}

\begin{lem}\label{phisemilin} The endomorphism $\varphi_{\widetilde{\bf D}}$ of
  $\widetilde{\bf D}$ is semilinear with respect to the ${\mathcal
  O}_{\mathcal E}$-action (\ref{oede}), hence it defines on $\widetilde{\bf
  D}$ the structure of an \'{e}tale $\varphi$-module over ${\mathcal
  O}_{\mathcal E}$.
\end{lem}

{\sc Proof:} \begin{align}\varphi_{\widetilde{\bf D}}(x\cdot((d_i)_{i}))&= \varphi_{\widetilde{\bf D}}((\varphi^i_{{\mathcal
  O}_{\mathcal E}}(x)d_i)_i)\notag\\{}&=((\varphi^i_{{\mathcal
  O}_{\mathcal E}}(x)d_{i+1})_{0\le i\le r-2},({\varphi}^r_{{\bf D}}(x\cdot d_0))_{r-1})\notag\\{}&=((\varphi^i_{{\mathcal
  O}_{\mathcal E}}(x)d_{i+1})_{0\le i\le r-2},({\varphi}^r_{{\mathcal
  O}_{\mathcal E}}(x){\varphi}^r_{{\bf D}}(d_0))_{r-1})\notag\\{}&={\varphi}_{{\mathcal
  O}_{\mathcal E}}(x)((d_{i+1})_{0\le i\le r-2},({\varphi}^r_{{\bf D}}(d_0))_{r-1})\notag\\{}&={\varphi}_{{\mathcal
  O}_{\mathcal E}}(x)\varphi_{\widetilde{\bf D}}((d_i)_{i})\notag.\end{align}\hfill$\Box$\\

Let $\Gamma'$ be an open subgroup of $\Gamma$, let ${\bf D}$ be an \'{e}tale $(\varphi^r,\Gamma')$-module over ${\mathcal
  O}_{\mathcal E}$. Define an action of $\Gamma'$ on $\widetilde{\bf D}$ by$$\gamma\cdot((d_i)_{0\le i\le r-1})=(\gamma\cdot d_i)_{0\le i\le r-1}.$$

\begin{lem}\label{gamsemilin} The $\Gamma'$-action on $\widetilde{\bf D}$
  commutes with $\varphi_{\widetilde{\bf D}}$ and is semilinear with respect to the ${\mathcal
  O}_{\mathcal E}$-action (\ref{oede}), hence we obtain on $\widetilde{\bf
  D}$ the structure of an \'{e}tale $(\varphi,\Gamma')$-module over ${\mathcal
  O}_{\mathcal E}$. We thus obtain an exact functor from the category of \'{e}tale $(\varphi^r,\Gamma')$-modules to the category of \'{e}tale $(\varphi,\Gamma')$-modules over ${\mathcal
  O}_{\mathcal E}$.
\end{lem}

{\sc Proof:} This is immediate from the respective properties of the
$\Gamma'$-action on ${\bf
  D}$. \hfill$\Box$\\

\section{The functor ${\bf D}$}

\label{secdfu}

 The topological dual $V^*$ of a
smooth ${\mathfrak N}_0$-representation on a torsion ${\mathfrak o}$-module $V$ (endowed with the discrete topology) is a compact left ${{{\mathcal O}_{\mathcal E}^+}}$-module, with ${{{\mathcal O}_{\mathcal E}^+}}$ acting through $(a\cdot f)(v)=f({a}\cdot
v)$ for $a\in {{{\mathcal O}_{\mathcal E}^+}}$ for $f\in V^*$ and $v\in
V$.\footnote{As ${\mathfrak N}_0$ is commutative this formula indeed defines an ${{{\mathcal O}_{\mathcal E}^+}}$-action. One might argue that it would be more natural to endow the dual with the ${{{\mathcal O}_{\mathcal E}^+}}$-action $(a\cdot f)(v)=f(\widetilde{a}\cdot
v)$ where $\widetilde{(.)}:{{{\mathcal O}_{\mathcal E}^+}}\to{{{\mathcal O}_{\mathcal E}^+}}$ denotes the involution induced by the inversion map
on ${\mathfrak N}_0$. Of course, everything we are going to develop holds true with this alternative ${{{\mathcal O}_{\mathcal E}^+}}$-action as well.}\\

Let ${\mathcal V}$ be a ${\mathfrak N}_0$-equivariant coefficient system on
${\mathfrak X}_+$ of level $1$. Clearly the ${\mathfrak N}_0$-action on both
$H_0(\overline{\mathfrak X}_+,{\mathcal V})$ and $H_0({\mathfrak
  X}_+,{\mathcal V})$ is smooth, thus$$D({\mathcal V})=H_0(\overline{\mathfrak
  X}_+,{\mathcal V})^*\quad\quad\quad\mbox{ and }\quad\quad\quad D'({\mathcal
  V})=H_0({\mathfrak X}_+,{\mathcal V})^*$$are compact ${{{\mathcal O}_{\mathcal E}^+}}$-modules. We put $${\bf D}({\mathcal V})={\mathcal O}_{\mathcal E}\otimes_{{{{\mathcal O}_{\mathcal E}^+}}}D({\mathcal V}).$$

\begin{pro}\label{strictiso} Suppose that ${\mathcal V}$ is strictly of level $1$. 

(a) $D({\mathcal V})$ can be generated as an ${{{\mathcal O}_{\mathcal E}^+}}$-module by ${\rm dim}_k({\mathcal V}({\mathfrak e}_0)\otimes_{\mathfrak o}k)$ many
  elements. In particular, ${\bf D}({\mathcal V})$ can be generated as an
  ${\mathcal O}_{\mathcal E}$-module by ${\rm dim}_k({\mathcal V}({\mathfrak e}_0)\otimes_{\mathfrak o}k)$ many
  elements. 

(b) The natural map ${\mathcal O}_{\mathcal E}\otimes_{{{{\mathcal O}_{\mathcal E}^+}}}D'({\mathcal V})\to {\bf
    D}({\mathcal V})$ is bijective.  
\end{pro}

{\sc Proof:} By Theorem \ref{finite} we know ${\mathcal V}({\mathfrak
  e}_0)\cong H_0(\overline{\mathfrak X}_+,{\mathcal V})^{{\mathfrak N}_0}$. It
follows that the $k$-vector space $H_0(\overline{\mathfrak X}_+,{\mathcal
  V})^{{\mathfrak N}_0,\pi_K=0}$ can be generated by ${\rm dim}_k({\mathcal V}({\mathfrak e}_0)\otimes_{\mathfrak o}k)$ many
  elements. By duality this means that the $k={{{\mathcal O}_{\mathcal E}^+}}/{\mathfrak m}$-vector space
$D({\mathcal V})/{\mathfrak m}D({\mathcal V})$ can be generated by ${\rm
  dim}_k({\mathcal V}({\mathfrak e}_0)\otimes_{\mathfrak o}k)$ many elements; here ${\mathfrak m}$
denotes the maximal ideal in the local ring ${{{\mathcal O}_{\mathcal E}^+}}$. Now we conclude with the topological Nakayama Lemma (see
\cite{baho}).

(b) We have an exact sequence $0\to {\mathcal V}({\mathfrak e}_0)\to H_0(\overline{\mathfrak
  X}_+,{\mathcal V})\to H_0({\mathfrak
  X}_+,{\mathcal V})\to0$ giving rise to an exact
sequence of ${{{{\mathcal O}_{\mathcal E}^+}}}$-modules (taking the Pontryagin dual is exact) $$0\longrightarrow D'({\mathcal V})\longrightarrow D({\mathcal
  V})\longrightarrow {\mathcal
  V}({\mathfrak e}_0)^*\longrightarrow0.$$Tensoring its first non trivial arrow with ${\mathcal O}_{\mathcal E}$ over
${{{{\mathcal O}_{\mathcal E}^+}}}$ gives the map
in question. Its surjectivity follows from ${\mathcal O}_{\mathcal E}\otimes_{{{{\mathcal O}_{\mathcal E}^+}} } {\mathcal
  V}({\mathfrak e}_0)^*=0$ which holds true because ${\mathcal
  V}({\mathfrak e}_0)$ and hence ${\mathcal
  V}({\mathfrak e}_0)^*$ is finitely generated over ${\mathfrak
    o}$.\hfill$\Box$\\

Now suppose that ${\mathcal V}$ is $\lfloor{\mathfrak
  N}_0,\varphi^r\rfloor$-equivariant for some $r\in{\mathbb N}$ such that the structure maps
$\varphi^r_x:{\mathcal V}(x)\to {\mathcal V}(\varphi^rx)$ and $\varphi^r_{\tau}:{\mathcal V}(\tau)\to {\mathcal V}(\varphi^r\tau)$ are bijective for
all $x\in{\mathfrak X}_+^0$ and $\tau\in{\mathfrak X}_+^1$. The $\varphi^r$-action provides an endomorphism $\varphi^r_{C_0({\mathfrak
  X}_+,{\mathcal V})}$ of $C_0({\mathfrak
  X}_+,{\mathcal V})=\oplus_{x\in{\mathfrak X}_+^0}{\mathcal V}(x)$ and hence an endomorphism $\varphi^r_{H_0(\overline{\mathfrak
  X}_+,{\mathcal V})}$ of $H_0(\overline{\mathfrak
  X}_+,{\mathcal V})$. We define$$\psi^r_{D({\mathcal V})}:D({\mathcal V})\longrightarrow D({\mathcal
  V}),\quad\quad d\mapsto d\circ \varphi^r_{H_0(\overline{\mathfrak
  X}_+,{\mathcal V})}.$$It is straightforward to check\begin{gather}\psi^r_{D({\mathcal
    V})}(\varphi^r_{{{{{\mathcal O}_{\mathcal E}^+}}}}(a)\cdot
d)=a\cdot(\psi^r_{D({\mathcal V})}(d))\quad\mbox{ for }a\in{{{\mathcal O}_{\mathcal E}^+}}, d\in D({\mathcal V}).\label{easyphi}\end{gather}As ${\mathcal
  O}_{\mathcal E}=\oplus_{n\in{\mathfrak
  N}_0/{\mathfrak N}_0^{p^r}}n\varphi^r_{{\mathcal
  O}_{\mathcal E}}({\mathcal
  O}_{\mathcal E})$, any element in ${\mathcal
  O}_{\mathcal E}$ can be written as a sum of products $\varphi_{{\mathcal
    O}_{\mathcal E}}^{r}(b)\cdot c$ with $b\in {\mathcal
  O}_{\mathcal E}$ and $c\in {{{\mathcal O}_{\mathcal E}^+}}$. Thus any
element in ${\bf
    D}({\mathcal V})$ is a sum of elements of the form $\varphi_{{\mathcal
    O}_{\mathcal E}}^{r}(b)\otimes d$ with $b\in {\mathcal
  O}_{\mathcal E}$ and $d\in D({\mathcal V})$. It therefore follows from
formula (\ref{easyphi}) that there is a well defined ${\mathfrak o}$-linear map $$\psi_{{\bf
    D}({\mathcal V})}^r:{\bf
    D}({\mathcal V})\longrightarrow{\bf
    D}({\mathcal V}),\quad\quad\varphi_{{\mathcal
    O}_{\mathcal E}}^{r}(b)\otimes d\mapsto b\otimes\psi^r_{D({\mathcal V})}(d).$$The map $\varphi^r:{\mathfrak
  X}_+^0\to{\mathfrak X}_+^0$ is injective. Thus the inverse of $\varphi^r$ induces an
isomorphism$$\varphi^{-r}:\bigoplus_{x\in{\varphi^r\mathfrak
      X}_+^0}{\mathcal
    V}(x)\stackrel{\cong}{\longrightarrow}\bigoplus_{x\in{\mathfrak
      X}_+^0}{\mathcal V}(x).$$We extend it to a map$$\psi^r_{C_0({\mathfrak
  X}_+,{\mathcal V})}:C_0({\mathfrak
  X}_+,{\mathcal V})=\bigoplus_{x\in{\mathfrak X}_+^0}{\mathcal
    V}(x)\stackrel{}{\longrightarrow}C_0({\mathfrak
  X}_+,{\mathcal V})=\bigoplus_{x\in{\mathfrak
      X}_+^0}{\mathcal V}(x)$$by requiring that its restriction to $\oplus_{x\in{\mathfrak X}_+^0-{\varphi^r\mathfrak
      X}_+^0}{\mathcal
    V}(x)$ vanishes. The definition implies\begin{gather}\psi^r_{C_0({\mathfrak
  X}_+,{\mathcal V})}(\varphi^r_{{{{\mathcal O}_{\mathcal E}^+}}}(a)\cdot h)=a\cdot\psi^r_{C_0({\mathfrak
  X}_+,{\mathcal V})}(h)\quad \mbox{ for }a\in{{{{\mathcal O}_{\mathcal E}^+}}}, h\in\bigoplus_{x\in{\mathfrak
      X}_+^0}{\mathcal V}(x).\label{psineufor}\end{gather}The map $\psi^r_{C_0({\mathfrak
  X}_+,{\mathcal V})}$ induces an endomorphism$$\psi^r_{H_0({\mathfrak
  X}_+,{\mathcal V})}:H_0({\mathfrak
  X}_+,{\mathcal V})\longrightarrow H_0({\mathfrak X}_+,{\mathcal
  V})$$and formula (\ref{psineufor}) becomes\begin{gather}\psi^r_{H_0({\mathfrak
  X}_+,{\mathcal V})}(\varphi^r_{{{{\mathcal O}_{\mathcal E}^+}}}(a)\cdot h)=a\cdot\psi^r_{H_0({\mathfrak
  X}_+,{\mathcal V})}(h)\label{h0phps}\end{gather}for $a\in{{{{\mathcal O}_{\mathcal E}^+}}}$ and $h\in H_0({\mathfrak
  X}_+,{\mathcal V})$. We define the endomorphism $$\varphi_{D'({\mathcal
    V})}^r:D'({\mathcal V})\longrightarrow D'({\mathcal V}),\quad\quad
d\mapsto d\circ \psi^r_{H_0({\mathfrak
  X}_+,{\mathcal V})}.$$We claim \begin{gather}\varphi_{D'({\mathcal V})}^r(a\cdot d)=\varphi_{{{{\mathcal O}_{\mathcal E}^+}}}^r(a)\cdot\varphi_{D'({\mathcal V})}^r(d)\quad\mbox{ for }a\in{{{\mathcal O}_{\mathcal E}^+}}, d\in  D'({\mathcal V}).\label{phps}\end{gather}Indeed, for $h\in H_0({\mathfrak
  X}_+,{\mathcal V})$ we compute \begin{align}(\varphi_{D'({\mathcal V})}^r(a\cdot
d))(h)&=(a\cdot d)(\psi^r_{H_0({\mathfrak
  X}_+,{\mathcal V})}(h))\notag\\{}&=d({a}\cdot\psi^r_{H_0({\mathfrak
  X}_+,{\mathcal V})}(h))\notag\\{}&\stackrel{(i)}{=}d(\psi^r_{H_0({\mathfrak
  X}_+,{\mathcal V})}(\varphi^r_{{{{\mathcal O}_{\mathcal E}^+}}}({a})\cdot h))\notag\\{}&=(\varphi_{D'({\mathcal V})}^r(d))({\varphi_{{{{\mathcal O}_{\mathcal E}^+}}}^r(a)}\cdot h)\notag\\{}&=(\varphi_{{{{\mathcal O}_{\mathcal E}^+}}}^r(a)\cdot\varphi_{D'({\mathcal V})}^r(d))(h)\notag\end{align}where in $(i)$ we used formula (\ref{h0phps}). Because of formula (\ref{phps}) we may proceed to define$$\varphi^r_{{\bf
    D}({\mathcal V})}=\varphi_{{\mathcal O}_{\mathcal E}}^r\otimes \varphi_{D'({\mathcal V})}^r:{\bf D}({\mathcal V})\longrightarrow{\bf
    D}({\mathcal V})$$where we use the isomorphism ${\mathcal O}_{\mathcal E}\otimes_{{{{\mathcal O}_{\mathcal E}^+}}}D'({\mathcal V})\cong {\bf
    D}({\mathcal V})$ of Proposition \ref{strictiso} as an identification. \\

\begin{pro}\label{etale} We have the formulae\begin{gather}\psi_{{\bf D}({\mathcal
      V})}^r\circ(b\cdot\varphi_{{\bf D}({\mathcal
      V})}^r)=\psi^r_{{\mathcal O}_{\mathcal E}}(b)\cdot{\rm id}_{{\bf D}({\mathcal
      V})}\quad\quad\mbox{ for }b\in {\mathcal O}_{\mathcal E},\label{phipsi1}\end{gather}\begin{gather}\sum_{n\in{\mathfrak
  N}_0/{\mathfrak
  N}_0^{p^r}}n\circ \varphi^r_{{\bf D}({\mathcal
      V})}\circ \psi^r_{{\bf D}({\mathcal
      V})}\circ n^{-1}={\rm id}_{{\bf D}({\mathcal
      V})}.\label{phipsi2}\end{gather}In particular, we have $\psi_{{\bf D}({\mathcal
      V})}^r\circ\varphi_{{\bf D}({\mathcal
      V})}^r={\rm id}_{{\bf D}({\mathcal
      V})}$ and $\varphi_{{\bf D}({\mathcal
      V})}^r$ is an \'{e}tale map.
\end{pro}

{\sc Proof:} To prove formula (\ref{phipsi1}) we first remark that for $b\in{{{\mathcal O}_{\mathcal E}^+}}$ we have $$\psi^r_{H_0({\mathfrak
  X}_+,{\mathcal V})}\circ(b\cdot\varphi^r_{H_0(\overline{\mathfrak
  X}_+,{\mathcal V})})=\psi^r_{{{{\mathcal O}_{\mathcal E}^+}}}({b})\cdot {\rm id}.$$Indeed, this is true already on
$0$-chains. To check this the $0$-chain may be assumed to be supported on a single vertex, and $b$
may be assumed to belong either to ${\rm im}(\varphi^r_{{{{\mathcal O}_{\mathcal E}^+}}})$ or to ${\rm ker}(\psi^r_{{{{\mathcal O}_{\mathcal E}^+}}})$; in either case the claim follows easily from the
definitions. We use this to see\begin{align}\psi_{D({\mathcal
      V})}^r(b\cdot\varphi_{D'({\mathcal V})}^r(d))(h)&=(b\cdot\varphi_{D'({\mathcal V})}^r(d))(\varphi^r_{H_0(\overline{\mathfrak
  X}_+,{\mathcal V})}(h))\notag\\{}&=d(\psi^r_{H_0({\mathfrak
  X}_+,{\mathcal V})}({b}\cdot\varphi^r_{H_0(\overline{\mathfrak
  X}_+,{\mathcal V})}(h)))\notag\\{}&=d(\psi^r_{{{{\mathcal O}_{\mathcal E}^+}}}({b})\cdot h)\notag\\{}&=\psi^r_{{{{\mathcal O}_{\mathcal E}^+}}}({b})\cdot d(h).\notag\end{align}Any element of ${\mathcal O}_{\mathcal E}$ can be written as a sum of products $b=b_1\cdot b_2$ with
$b_1\in{\mathcal O}_{\mathcal E}$ and $b_2\in {{{\mathcal O}_{\mathcal E}^+}}$ and $b_1=\varphi^r_{{\mathcal O}_{\mathcal
    E}}(\psi^r_{{\mathcal O}_{\mathcal E}}(b_1))$. Inserting what we just saw we compute\begin{align}\psi_{{\bf D}({\mathcal V})}^r(b\cdot\varphi_{{\bf
     D}({\mathcal V})}^r)(a\otimes d)&=\psi_{{\bf D}({\mathcal
     V})}^r(b\cdot\varphi_{{\mathcal O}_{\mathcal E}}^r(a)\otimes
 \varphi_{D'({\mathcal V})}^r(d))\notag\\{}&=\psi_{{\bf D}({\mathcal
     V})}^r(\varphi^r_{{\mathcal O}_{\mathcal
    E}}(\psi^r_{{\mathcal O}_{\mathcal E}}(b_1))\cdot\varphi_{{\mathcal O}_{\mathcal E}}^r(a)\otimes
 b_2\cdot\varphi_{D'({\mathcal V})}^r(d))\notag\\{}&=\psi_{{\bf D}({\mathcal
     V})}^r(\varphi^r_{{\mathcal O}_{\mathcal
    E}}(\psi^r_{{\mathcal O}_{\mathcal E}}(b_1)\cdot a)\otimes
 b_2\cdot\varphi_{D'({\mathcal V})}^r(d))\notag\\{}&=\psi^r_{{\mathcal O}_{\mathcal
     E}}(b_1)\cdot a\otimes\psi_{D({\mathcal
     V})}^r(b_2\cdot\varphi_{D'({\mathcal V})}^r(d))\notag\\{}&=\psi^r_{{\mathcal O}_{\mathcal
     E}}(b_1)\cdot a\otimes \psi^r_{{{{\mathcal O}_{\mathcal E}^+}}}(b_2)\cdot d\notag\\{}&=\psi^r_{{\mathcal O}_{\mathcal
     E}}(b_1)\cdot \psi^r_{{\mathcal O}_{\mathcal
     E}}(b_2)\cdot a\otimes  d\notag\\{}&=\psi^r_{{\mathcal O}_{\mathcal
     E}}(b)\cdot a\otimes d.\notag\end{align}We have proven formula (\ref{phipsi1}). To prove formula (\ref{phipsi2}) we view the injective map $D'({\mathcal V})\to D({\mathcal V})$ as an inclusion. We find some $N>0$ such that $\psi^r_{D({\mathcal
      V})}(t^ND({\mathcal V}))\subset D'({\mathcal V})$. The map $t^ND({\mathcal V})\to D({\mathcal V})$ induces an isomorphism ${\mathcal O}_{\mathcal E}\otimes_{{{{\mathcal O}_{\mathcal E}^+}}}t^ND({\mathcal V})\cong {\mathcal O}_{\mathcal E}\otimes_{{{{\mathcal O}_{\mathcal E}^+}}}D({\mathcal V})={\bf
    D}({\mathcal V})$, therefore we may write an element in ${\bf
    D}({\mathcal V})$ as a sum of elements $\varphi_{{\mathcal O}_{\mathcal E}}^r(b)\otimes d$ with $d\in D({\mathcal V})$ such that $\psi^r_{D({\mathcal
      V})}(d)\in D'({\mathcal V})$, and then more generally $\psi^r_{D({\mathcal
      V})}(n^{-1}\cdot d)\in D'({\mathcal V})$ for $n\in {\mathfrak
  N}_0$. We compute\begin{align}\sum_{n\in{\mathfrak
  N}_0/{\mathfrak
  N}_0^{p^r}}n\circ \varphi^r_{{\bf D}({\mathcal
      V})}\circ \psi^r_{{\bf D}({\mathcal
      V})}\circ n^{-1}(\varphi_{{\mathcal O}_{\mathcal E}}^r(b)\otimes d)&=\sum_{n\in{\mathfrak
  N}_0/{\mathfrak
  N}_0^{p^r}}n\circ \varphi^r_{{\bf D}({\mathcal
      V})}\circ \psi^r_{{\bf D}({\mathcal
      V})} (\varphi_{{\mathcal O}_{\mathcal E}}^r(b)\otimes n^{-1}\cdot d)\notag\\{}&=\sum_{n\in{\mathfrak
  N}_0/{\mathfrak
  N}_0^{p^r}}n\circ \varphi^r_{{\bf D}({\mathcal
      V})}(b\otimes \psi^r_{D({\mathcal
      V})}(n^{-1}\cdot d))\notag\\{}&=\sum_{n\in{\mathfrak
  N}_0/{\mathfrak
  N}_0^{p^r}}\varphi^r_{{\mathcal O}_{\mathcal E}}(b)\otimes n\cdot(\varphi^r_{{D}'({\mathcal
      V})}(\psi^r_{D({\mathcal
      V})}(n^{-1}\cdot d))).\notag\end{align}Therefore we need to show the equality\begin{gather}\sum_{n\in{\mathfrak
  N}_0/{\mathfrak
  N}_0^{p^r}}n\cdot(\varphi^r_{{D}'({\mathcal
      V})}(\psi^r_{D({\mathcal
      V})}(n^{-1}\cdot d)))=d\label{torsko}\end{gather}of linear forms on $H_0(\overline{\mathfrak X}_+,{\mathcal
  V})$. An inductive argument using formula
(\ref{striclevgen}) shows that any element of $H_0(\overline{\mathfrak X}_+,{\mathcal
  V})$ can be represented by a $0$-chain supported on$${\mathfrak N}_0\varphi^r{\mathfrak X}_+^0=\coprod_{n_0\in{\mathfrak
  N}_0/{\mathfrak
  N}_0^{p^r}}n_0\varphi^r{\mathfrak X}_+^0.$$Thus, to prove (\ref{torsko}) it is enough, by the definitions of $\varphi^r_{{D}'({\mathcal
      V})}$ and $\psi^r_{D({\mathcal
      V})}$, to prove\begin{gather}\sum_{n\in{\mathfrak N}_0/{\mathfrak N}_0^{p^r}}(n^{-1}\circ \varphi^r_{C_0({\mathfrak X}_+,{\mathcal V})} \circ\psi^r_{C_0({\mathfrak
  X}_+,{\mathcal V})}\circ n)(c)=c\label{beidpa}\end{gather}for all $c\in C_0({\mathfrak X}_+,{\mathcal V})$ supported on $n_0\varphi^r{\mathfrak X}_+^0$
for some $n_0\in{\mathfrak
  N}_0$. For $n\in{\mathfrak N}_0$ with $nn_0\notin {\mathfrak
  N}_0^{p^r}$ we have $nn_0\varphi^r{\mathfrak X}_+^0\cap
\varphi^r{\mathfrak X}_+^0=\emptyset$ and then $c$ is killed by
$\psi^r_{C_0({\mathfrak
  X}_+,{\mathcal V})}\circ n$. On the other hand, if $nn_0\in {\mathfrak
  N}_0^{p^r}$ then $c$ is fixed by $n^{-1}\circ \varphi_{C_0({\mathfrak
  X}_+,{\mathcal V})}^r\circ
\psi_{C_0({\mathfrak
  X}_+,{\mathcal V})}^r\circ n$. Formula (\ref{beidpa}) is proven.
 \hfill$\Box$\\

\begin{kor}\label{torsch} (a) The ${\mathcal O}_{\mathcal E}$-module ${\bf D}({\mathcal
      V})$ can be generated by ${\rm dim}_k({\mathcal V}({\mathfrak e}_0)\otimes_{\mathfrak o}k)$ many elements and carries a natural structure of an \'{e}tale
    $(\varphi^r,\Gamma_{{{0}}})$-module. 

(b) If the $\lfloor{\mathfrak
  N}_0,\varphi^r\rfloor$-action on ${\mathcal V}$ extends to a $\lfloor{\mathfrak
  N}_0,\varphi^r,\Gamma\rfloor$-action on ${\mathcal V}$, then ${\bf D}({\mathcal
      V})$ is an \'{e}tale
    $(\varphi^r,\Gamma)$-module.  
\end{kor}

{\sc Proof:} (a) Except for the $\Gamma_{{{0}}}$-action everything else has already been established in Propositions
    \ref{strictiso} and \ref{etale}. By proposition \ref{natgam} the $\lfloor{\mathfrak
  N}_0,\varphi^r\rfloor$-action on ${\mathcal V}$ naturally extends to a $\lfloor{\mathfrak
  N}_0,\varphi^r,\Gamma_{{{0}}}\rfloor$-action on ${\mathcal V}$. In particular it
induces a $\Gamma_{{{0}}}$-action on $H_0({\mathfrak
  X}_+,{\mathcal V})$. We define a $\Gamma_{{{0}}}$-action on $D'({\mathcal
      V})$ by setting $(\gamma\cdot d)(h)=d(\gamma^{-1}\cdot h)$ for
    $\gamma\in\Gamma_{{{0}}}$, $d\in D'({\mathcal
      V})$, $h\in H_0({\mathfrak
  X}_+,{\mathcal V})$. For $a\in {{{\mathcal O}_{\mathcal E}^+}}$ we compute$$(\gamma\cdot(a\cdot d))(h)=(a\cdot d)(\gamma^{-1}\cdot
h)=d({a}\cdot\gamma^{-1}\cdot
h)=(\gamma\cdot d)(\gamma a\gamma^{-1}\cdot
h)$$$$=((\gamma{a}\gamma^{-1})\cdot (\gamma \cdot d))(h)=((\gamma\cdot{a})\cdot(\gamma\cdot d))(h)$$i.e. the action is
semilinear. Moreover, it follows immediately from the definitions that it commutes with $\varphi^r_{D'({\mathcal
      V})}$. Therefore we obtain a semilinear action of $\Gamma_{{{0}}}$ on ${\bf D}({\mathcal
      V})$, commuting with $\varphi^r_{{\bf D}({\mathcal
      V})}$, by putting $\gamma(a\otimes d)=\gamma\cdot a\otimes \gamma\cdot d$
    for $a\in {\mathcal O}_{\mathcal E}$ and $d\in D'({\mathcal
      V})$.

Of course, the same construction also endows $D({\mathcal
      V})$ with a semilinear $\Gamma_{{{0}}}$-action, commuting with $\psi^r_{D({\mathcal
      V})}$, and extending to the same $\Gamma_{{{0}}}$-action on ${\bf D}({\mathcal
      V})$. 

(b) This is the same argument (without invoking proposition
\ref{natgam}).\hfill$\Box$\\

\begin{kor}\label{reconst} (a) The restriction of $\psi_{{\bf D}({\mathcal
      V})}^r$ to $\cap_{n\in{\mathfrak N}_0-{\mathfrak N}_0^{p^r}}{\rm ker}(\psi_{{\bf D}({\mathcal
      V})}^r(n\cdot .))$ is a bijection onto ${\bf D}({\mathcal
      V})$. Its inverse, composed with the inclusion of $\cap_{n\in{\mathfrak N}_0-{\mathfrak N}_0^{p^r}}{\rm ker}(\psi_{{\bf D}({\mathcal
      V})}^r(n\cdot .))$ into ${\bf D}({\mathcal
      V})$, is the map $\varphi_{{\bf D}({\mathcal
      V})}^r$; in particular, the latter can be reconstructed from $\psi_{{\bf D}({\mathcal
      V})}^r$. 

(b) Conversely, $\psi_{{\bf D}({\mathcal
      V})}^r$ can be reconstructed from $\varphi_{{\bf D}({\mathcal
      V})}^r$. 
\end{kor}

{\sc Proof:} This is a formal consequence of Proposition \ref{etale}.

(a) (See \cite{vien} Proposition 3.3.24 for the abstract argument.) All we need to show is$${\rm im }(\varphi_{{\bf D}({\mathcal
      V})}^r)=\bigcap_{n\in{\mathfrak N}_0-{\mathfrak N}_0^{p^r}}{\rm ker}(\psi_{{\bf D}({\mathcal
      V})}^r(n\cdot .)).$$Let $d\in {\rm im }(\varphi_{{\bf D}({\mathcal
      V})}^r)$, say $d=\varphi_{{\bf D}({\mathcal
      V})}^r(c)$. Let $n\in {\mathfrak N}_0-{\mathfrak N}_0^{p^r}$. Then $\psi_{{\bf D}({\mathcal
      V})}^r(n\cdot d)=\psi_{{\bf D}({\mathcal
      V})}^r(n\cdot \varphi_{{\bf D}({\mathcal
      V})}^r(c))=\psi^r_{{\mathcal O}_{\mathcal E}}(n)\cdot c=0$ where we used
  formula (\ref{phipsi1}) and then $\psi^r_{{\mathcal O}_{\mathcal
      E}}(n)=0$. Conversely, let $d\in\cap_{n\in{\mathfrak N}_0-{\mathfrak N}_0^{p^r}}{\rm ker}(\psi_{{\bf D}({\mathcal
      V})}^r(n\cdot .))$. Formula (\ref{phipsi2}) shows $$d=\sum_{n\in{\mathfrak
  N}_0/{\mathfrak
  N}_0^{p^r}}n\cdot\varphi^r_{{\bf D}({\mathcal
      V})}(\psi^r_{{\bf D}({\mathcal
      V})}(n^{-1}\cdot d)).$$By hypothesis, only the summand for the coset ${\mathfrak
  N}_0^{p^r}$ survives, showing $d=\varphi^r_{{\bf D}({\mathcal
      V})}(\psi^r_{{\bf D}({\mathcal
      V})}(d))$, hence $d\in {\rm im }(\varphi_{{\bf D}({\mathcal
      V})}^r)$.

(b) Proposition \ref{etale} implies ${\bf D}({\mathcal
      V})={\rm im}(\varphi_{{\bf D}({\mathcal
      V})}^r)\oplus\sum_{n\in{\mathfrak N}_0-{\mathfrak N}_0^{p^r}}n\cdot{\rm im}(\varphi_{{\bf D}({\mathcal
      V})}^r)$ and that $\varphi_{{\bf D}({\mathcal
      V})}^r$ is injective. Thus $\psi_{{\bf D}({\mathcal
      V})}^r$ is the projection onto ${\rm im}(\varphi_{{\bf D}({\mathcal
      V})}^r)$ with kernel $\sum_{n\in{\mathfrak N}_0-{\mathfrak N}_0^{p^r}}n\cdot{\rm im}(\varphi_{{\bf D}({\mathcal
      V})}^r)$, composed with the inverse of $\varphi_{{\bf D}({\mathcal
      V})}^r$.\hfill$\Box$\\

We fix a gallery (\ref{cgall}) and choose an isomorphism $\Theta:Y\stackrel{\cong}{\to}{\mathfrak
  X}_+$ as in Theorem \ref{baumeinbettung}. 

\begin{satz}\label{exakt} (a) Given $\phi$ as in Theorem
  \ref{baumeinbettung}(b), the assignment $M\mapsto {\bf D}(\Theta_*{\mathcal
      V}_M)$ is an exact contravariant functor from ${\rm Mod}^{\rm fin}({\mathcal
  H}(G,I_0)_{{\mathfrak o}_m})$ to the category of \'{e}tale
    $(\varphi^r,\Gamma_0)$-modules over
${\mathcal O}_{\mathcal E}$.

(b) Given $\phi$
and $\tau$ as in Theorem \ref{baumeinbettung}(d), the assignment $M\mapsto {\bf D}(\Theta_*{\mathcal
      V}_M)$ is an exact contravariant functor from ${\rm Mod}^{\rm fin}({\mathcal
  H}(G,I_0)_{{\mathfrak o}_m})$ to the category of \'{e}tale
    $(\varphi^r,\Gamma)$-modules over
${\mathcal O}_{\mathcal E}$.

(c) The functors in (a) and (b) depend canonically on the choice of
(\ref{cgall}) and $\phi$ alone, resp. of (\ref{cgall}) and $\phi$ and $\tau$ alone,
not on the choice of $\Theta$. 

(d) For $M\in{\rm Mod}^{\rm fin}({\mathcal
  H}(G,I_0)_{{\mathfrak o}_1})={\rm Mod}^{\rm fin}({\mathcal
  H}(G,I_0)_{k})$ we have$${\rm dim}_{k_{\mathcal E}}{\bf D}(\Theta_*{\mathcal
      V}_M)\le{\rm dim}_kM.$$
\end{satz}

{\sc Proof:} The assignment
$M\mapsto {\mathcal V}_M$ (and hence $M\mapsto \Theta_*{\mathcal V}_M$) is exact by Proposition \ref{jalgheck}. Therefore
the assignment $M\mapsto H_0(\overline{\mathfrak
  X}_+,\Theta_*{\mathcal V}_M)$ is exact, cf. the exact sequence
(\ref{h0ex}). By exactness of taking Pontryagin duals it follows that $M\mapsto D(\Theta_*{\mathcal V}_M)$ is exact. Finally, by the
flatness of ${{{\mathcal O}_{\mathcal E}^+}}\to{\mathcal O}_{\mathcal E}$ we
obtain that $M\mapsto {\bf D}(\Theta_*{\mathcal V}_M)$ is exact.

The independence on the choice of $\Theta$ follows from the corresponding independence
statement in Theorem \ref{pufo}.\hfill$\Box$\\

\section{The case ${\rm GL}_{d+1}({\mathbb Q}_p)$}

\label{gln}

We consider the case $G={\rm GL}_{d+1}({\mathbb Q}_p)$ for some $d\ge1$ and
keep all the previous notations $T, N(T), X, A$
etc.. We fix a chamber $C$ in $A$, and as before we denote by $I$ resp. $I_0$ the corresponding Iwahori subgroup, resp. pro-$p$-Iwahori subgroup of $G$. The (affine) reflections
in the codimension-$1$-faces of $C$ form a set $S$ of Coxeter
generators for the affine Weyl group which we view as a subgroup of the
extended affine Weyl group $N(T)/Z(T\cap I)$. Put $W=N(T)/T$, the finite Weyl group.

We find elements $u,s_d\in N(T)$ such that $uC=C$ (equivalently, $uI=Iu$, or also $uI_0=I_0u$), such that $u^{d+1}\in\{p\cdot {\rm id}, p^{-1}\cdot {\rm
  id}\}$ and such that, setting $$s_i=u^{d-i}s_du^{i-d}\quad\quad\mbox{ for }0\le i\le
d$$the set $\{s_0,s_1,\ldots,s_d\}$
maps bijectively to $S$; we henceforth regard this bijection as an identifications. Let $\iota:{\rm SL}_2({\mathbb Q}_p)\to G$ denote the
embedding corresponding to $s_d$. For $0\le
i\le d$ we
put $$n_{s_i}=u^{d-i}\cdot\iota(\left(\begin{array}{cc}0&1\\-1&0\end{array}\right))\cdot
u^{i-d}\quad\mbox{
  and
}\quad
h_{s_i}(x)=u^{d-i}\cdot\iota(\left(\begin{array}{cc}x&0\\0&x^{-1}\end{array}\right))\cdot
u^{i-d}$$for $x\in{\mathbb
F}_p^{\times}$ where we use the Teichm\"uller character to regard
$\left(\begin{array}{cc}x&0\\0&x^{-1}\end{array}\right)$ as an element of
${\rm SL}_2({\mathbb Q}_p)$. Similarly, by means of the Teichm\"uller
character we regard the group $\overline{T}=I/I_0=(T\cap I)/(T\cap I_0)$ as a subgroup of $T$.

We write ${\mathcal H}(G,I_0)_k={\mathcal H}(G,I_0)\otimes_{\mathfrak
  o}k$. Let ${\mathcal H}(G,I_0)_{{\rm aff},k}$ denote the $k$-subalgebra of
${\mathcal H}(G,I_0)_k$ generated by the $T_{n_s}$ for $s\in{S}$ and the $T_t$ for $t\in \overline{T}$. Let ${\mathcal
  H}(G,I_0)'_{{\rm aff},k}$ denote the $k$-subalgebra of
${\mathcal H}(G,I_0)_k$ generated by ${\mathcal H}(G,I_0)_{{\rm aff},k}$
together with the elements $T_{p\cdot{\rm id}}$ and $T_{p^{-1}\cdot{\rm
    id}}=T^{-1}_{p\cdot{\rm id}}$. We put $$\phi=s_du\in N(T)\quad\quad\quad\mbox{ and }\quad\quad\quad
C^{(i)}=\phi^iC\quad\mbox{ for }i\ge0.$$We have $\phi^d=\xi(p)\in T$ for some
$\xi\in{\rm Hom}_{{\rm alg}}({\mathbb G}_m,T)$. Let $$\tau=(.)^m\cdot\xi|_{{\mathbb
    Z}_p^{\times}}:{\mathbb Z}_p^{\times}\longrightarrow T,\quad
a\mapsto a^m\cdot\xi(a)$$for some $m\in{\mathbb Z}$. For $x\in{\mathbb F}_p^{\times}$ use the Teichm\"uller lifting to define $\tau(x)\in T$.

\begin{lem} (For a suitable choice of $N_0$ we have:) $\{C^{(i)}\}_{i\ge0}$, $\phi$ and
$\tau$ satisfy the assumptions of Theorem \ref{baumeinbettung} (with $r=1$ there).
\end{lem}

{\sc Proof:} We choose a system $\Delta$ of simple roots in such a way that the image of
$S-\{s_0\}$ in $W$ is the set of simple
reflections corresponding to $\Delta$. (The above embedding $\iota:{\rm SL}_2({\mathbb
  Q}_p)\to G$ is then taken, more precisely, to be the one corresponding to
the simple root associated with $s_d$.) We take $N_0$ to be the group of ${\mathbb Z}_p$-valued points
of the unipotent radical of the Borel subgroup containing $T$ corresponding to
$\Delta$. Then $\phi^{d}$ belongs to $T$ and we have
$\phi^{d}N_0\phi^{-d}\subset N_0$ with
$[N_0:\phi^{d}N_0\phi^{-d}]=p^{d}$. From this the claims easily follow. (All this can also be checked be means of the
concrete realizations of $T$, $\phi$, $u$ etc. given below.) \hfill$\Box$\\

We have the isomorphism $\Theta:Y\cong{\mathfrak X}_+$ as constructed
in Theorem \ref{baumeinbettung}, by choosing the element $\nu_0$ in the proof
of Theorem \ref{baumeinbettung} to be $\nu_0=\iota(\nu)$. Recall that the chamber $C$ of $X$ corresponds to the edge ${\bf e}_0$ of $Y$. The codimension-$1$-face
of $C$ corresponding to the vertex ${\bf v}_0$ of $Y$ is $$F=C\cap C^{(1)}=C\cap\phi C=C\cap s_du C=C\cap
s_dC.$$ Therefore $s_d$ is the simple reflection corresponding to the simple
root $\alpha^{(0)}$, and in the present setting we have
$$\varphi=\phi=s_du.$$

Placing ourselves into the setting of sections \ref{pmodsl2} and \ref{iwah} we observe:

\begin{lem}\label{thetaslink} (a) The image of ${\mathfrak N}_0\subset{\rm SL}_2({\mathbb Z}_p)$ in $\overline{\mathcal S}={\rm SL}_2({\mathbb F}_p)$ is ${\mathfrak N}_0/{\mathfrak
  N}_0^p$, and this is the unipotent radical $\overline{\mathcal U}$ of a Borel subgroup in $\overline{\mathcal S}$.

(b) We have an isomorphism between $\overline{\mathcal S}$ and the maximal reductive
(over ${\mathbb F}_p$) quotient of $I_0^F$, inducing an isomorphism between $\overline{\mathcal U}$ and the image of $I_0$, hence an
embedding of $k$-algebras\begin{align}{\mathcal
  H}(\overline{\mathcal S},\overline{\mathcal U})_k={\rm End}_{k[\overline{\mathcal S}]}({\rm
  ind}_{\overline{\mathcal U}}^{\overline{\mathcal S}}{\bf 1}_k)^{\rm
  op}&\cong{\rm End}_{k[I_0^F]}({\rm
  ind}_{I_0}^{I_0^F}{\bf 1}_k)^{\rm op}\notag\\{}&\hookrightarrow{\mathcal H}(G,I_0)_{{\rm
    aff},k}\quad\subset\quad {\mathcal H}(G,I_0)_{k}.\label{heckeinbet}\end{align}It
sends $T_{h_s(x)}$ to $T_{h_{s_d}(x)}$ (for $x\in{\mathbb F}_p^{\times}$) and
$T_{n_s}$ to $T_{n_{s_d}}$.    
\end{lem}

{\sc Proof:} (a) is clear. (b) The subgroup of ${\rm SL}_2({\mathbb Q}_p)$ generated by the stabilizers of
the edges in $\overline{\mathfrak X}_+$ emanating from ${\mathfrak v}_0$ is
${\rm SL}_2({\mathbb Z}_p)$, its maximal reductive (over
${\mathbb F}_p$) quotient is $\overline{\mathcal S}={\rm SL}_2({\mathbb F}_p)$. From this everything follows.\hfill$\Box$\\

For a ${\mathcal H}(G,I_0)_k$-module $M$, Theorem \ref{pufo} and Proposition \ref{jalgheck} tell us that $\Theta_*{\mathcal V}_{M}$ is in a natural way a $\lfloor{\mathfrak N}_0,\varphi,\Gamma\rfloor$-equivariant coefficient system on ${\mathfrak X}_+$, with ${\mathfrak N}_0$-action strictly of level $1$. For later use we remark that $\Theta_*{\mathcal
  V}_M({\mathfrak v}_0)={\mathcal V}_{M}({\bf v}_0)={\mathcal V}_{M}^X(F)$ is
stable under $s_d$ (which acts on the $G$-equivariant coefficient system
${\mathcal V}_{M}^X$).\\

From now on, for concreteness, we specialize our discussion as follows. $T$ is the
subgroup of diagonal matrices in $G$ and $I_0$ is the pro-$p$-Iwahori
subgroup of $G$ consisting of all matrices in ${\rm GL}_{d+1}({\mathbb Z}_p)$ which (minus the
identity) are strictly upper triangular modulo $p$. We further assume that we are in one of the following two 'opposite' cases (which we treat simultaneously).\footnote{Notice that our
analysis requires explixit matrix computations only when we justify certain statements
involving $\tau$. For some of
these statements our specific choice of $\tau$ is important, i.e. some of
them may fail if $\tau$ is replaced by $a\mapsto a^m\cdot\tau(a)$ for some arbitrary $m\in{\mathbb Z}$. In general, the condition
to be imposed on $\tau$ should be that the image of $\tau$ together with ${\rm SL}_{d+1}({\mathbb F}_p)$
generates ${\rm GL}_{d+1}({\mathbb F}_p)$.} The first one is where
  \begin{gather}u=\left(\begin{array}{cc}0&E_d\\p&0\end{array}\right),\quad\quad s_d=\left(\begin{array}{cc}E_{d-1}&0\\0&s\end{array}\right),\quad\quad\tau(a)=\left(\begin{array}{cc}
    E_d&0\\0&a^{-1}\end{array}\right)\label{exam}\end{gather}for $a\in{\mathbb
Z}_p^{\times}$. Here $E_l$ denotes the identity $l\times l$-matrix, and
$s=\left(\begin{array}{cc}0&1\\1&0\end{array}\right)\in{{\rm GL}_2}$. In this
case we have $N_0'=\left(\begin{array}{cc}
    E_d&*\\0&1\end{array}\right)$ (with column $*$ having entries in ${\mathbb
  Z}_p$), where $N_0'$ denotes the subgroup of $G$ generated by all the
$N_{{\alpha}^{(j)}}\cap N_0$ for $j\ge0$ (cf. the remark following Theorem \ref{baumeinbettung}).

The 'opposite' case (corresponding to the other end of the Dynkin diagram) is where\begin{gather}u=\left(\begin{array}{cc}0&p^{-1}\\E_d&0\end{array}\right),\quad\quad s_d=\left(\begin{array}{cc}-s&0\\0&E_{d-1}\end{array}\right),\quad\quad\tau(a)=\left(\begin{array}{cc}a&0\\0&E_d\end{array}\right).\label{exam1}\end{gather}In this
case we have $N_0'=\left(\begin{array}{cc}1&*\\0&E_d\end{array}\right)$. 

\subsection{Supersingular ${\mathcal H}(G,I_0)_k$-modules}

\label{supsing}

 For a character $\lambda:\overline{T}\to
k^{\times}$ we denote by $S_{\lambda}$ the set of all $s\in{S}$ with
$\lambda(h_s(x))=1$ for all $x\in{\mathbb
F}_p^{\times}$. Suppose we are given a character $\lambda:\overline{T}\to
k^{\times}$ and a subset ${\mathcal
  J}$ of $S_{\lambda}$. There
is a uniquely determined character $$\chi_{\lambda,{\mathcal
  J}}:{\mathcal H}(G,I_0)_{{\rm
    aff},k}\longrightarrow k$$which sends $T_t$ to $\lambda(t^{-1})$ for $t\in
\overline{T}$, which sends $T_{n_s}$ to $0$ for $s\in {S}-{\mathcal
  J}$ and which
sends $T_{n_s}$ to $-1$ for $s\in {\mathcal
  J}$ (see \cite{vigneras} Proposition 2). Let $b\in k^{\times}$. The character $\chi_{\lambda,{\mathcal
  J}}$ extends uniquely to a character$$\chi_{\lambda,{\mathcal
  J},b}:{\mathcal H}(G,I_0)'_{{\rm
    aff},k}\longrightarrow k$$which sends $T_{u^{d+1}}$ to $b$ (see the
proof of \cite{vigneras} Proposition 3). We
define the ${\mathcal H}(G,I_0)_k$-module$$M[\lambda,{\mathcal
  J},b]={\mathcal H}(G,I_0)_k\otimes_{{\mathcal H}(G,I_0)'_{{\rm
    aff},k}}k.e$$where $k.e$ denotes the one dimensional $k$-vector space on
the basis element $e$, endowed with the action of ${\mathcal H}(G,I_0)'_{{\rm
    aff},k}$ by the character $\chi_{\lambda,{\mathcal
  J},b}$. 

For $0\le i\le d$ we let $e_i=T_{u^{-i}}\otimes e\in M[\lambda,{\mathcal
  J},b]$. The pair $(\lambda^{[i]},{\mathcal
  J}^{[i]})$ defined by $\lambda^{[i]}(t)=\lambda(u^{-i}tu^{i})$ and ${\mathcal
  J}^{[i]}=u^i{\mathcal
  J}u^{-i}$ satisfies the same assumptions as $(\lambda,{\mathcal
  J})$, hence gives rise to a corresponding
character $\chi_{\lambda^{[i]},{\mathcal
  J}^{[i]},b}$ of ${\mathcal H}(G,I_0)'_{{\rm
    aff},k}$. 

For $m\in{\mathbb Z}$ and $0\le i\le d$ with $m-i\in(d+1){\mathbb Z}$ we set
$e_{m}=e_i$, ${\mathcal
  J}^{[m]}={\mathcal
  J}^{[i]}$ and $\lambda^{[m]}=\lambda^{[i]}$, and
similarly for all other objects defined below which are indexed by $0\le i\le d$. Consider the condition\begin{gather}\mbox{ The pairs }(\lambda^{[0]},{\mathcal
  J}^{[0]}),\ldots,(\lambda^{[d]},{\mathcal
  J}^{[d]})\mbox{ are pairwise distinct.}\label{diff}\end{gather}

\begin{pro}\label{marfrarac} (a) If $(\lambda,{\mathcal
  J})$ satisfies (\ref{diff}) then $M[\lambda,{\mathcal
  J},b]$ is an absolutely simple supersingular ${\mathcal H}(G,I_0)_k$-module of
$k$-dimension $d+1$. For another pair $(\lambda',{\mathcal
  J}')$ satisfying (\ref{diff}) the ${\mathcal H}(G,I_0)_k$-modules $M[\lambda,{\mathcal
  J},b]$ and $M[\lambda',{\mathcal
  J}',b]$ are isomorphic if and only if $(\lambda',{\mathcal
  J}')$ and $(\lambda,{\mathcal
  J})$ are conjugate by some power of $u$. 

(b) Any absolutely simple supersingular
  ${\mathcal H}(G,I_0)_k$-module of $k$-dimension $d+1$ is isomorphic with $M[\lambda,{\mathcal
  J},b]$ for
  suitable $\lambda$, ${\mathcal
  J}$, $b$ satisfying (\ref{diff}).
 
(c) As a ${\mathcal H}(G,I_0)'_{{\rm
    aff},k}$-module, $M[\lambda,{\mathcal
  J},b]$ decomposes as\begin{gather}M[\lambda,{\mathcal
  J},b]\cong\bigoplus_{0\le i\le
  d}k.e_i\label{modgplude}\end{gather}with ${\mathcal H}(G,I_0)'_{{\rm
    aff},k}$ acting on $k.e_i$ by the character $\chi_{\lambda^{[i]},{\mathcal
  J}^{[i]},b}$.
\end{pro}

{\sc Proof:} For (a) see \cite{vigneras} Proposition 3 and Theorem 5. Notice
that conjugating the pair $(\lambda,{\mathcal
  J})$ by powers of $u$ is equivalent with cyclically permuting the set of
pairs $\{(\lambda^{[i]},{\mathcal
  J}^{[i]})\}_{i}$. For (b) see \cite{vigneras} Theorem 5 together with \cite{oll} Theorem 7.3. For (c) see the proof of \cite{vigneras} Proposition
3. To see e.g. that $T_t$ for $t\in\overline{T}$ acts by $\lambda^{[i]}(t^{-1})$ on
$k.e_i$ we compute$$T_te_i=T_tT_{u^{-i}}\otimes e=T_{u^{-i}t}\otimes
e=T_{u^{-i}}T_{u^{-i}tu^i}\otimes e$$$$=T_{u^{-i}}\otimes T_{u^{-i}tu^i}e=T_{u^{-i}}\otimes\lambda(u^{-i}t^{-1}u^i)e=T_{u^{-i}}\otimes\lambda^{[i]}(t^{-1}) e=\lambda^{[i]}(t^{-1})e_i.$$\hfill$\Box$\\

For $0\le i\le d$ we define a number $0\le k_i=k_i(\lambda,{\mathcal
  J})\le p-1$ such that\begin{gather}\lambda^{[i]}(h_{s_d}(x))=\lambda(h_{s_{i-1}}(x))=x^{k_i}\quad\quad\mbox{ for all
}x\in{\mathbb F}_p^{\times},\label{kljfml}\end{gather}as follows. If
$\lambda^{[i]}\circ h_{s_d}$ is not the constant character ${\bf 1}$ then $k_i$ is already uniquely determined by formula (\ref{kljfml}). Next notice that $\lambda^{[i]}\circ
h_{s_d}={\bf 1}$ is equivalent with $s_{i-1}\in S_{\lambda}$. If
$\lambda^{[i]}\circ h_{s_d}={\bf 1}$ and $s_{i-1}\in {\mathcal
  J}$ we put $k_i=p-1$, if
$\lambda^{[i]}\circ h_{s_d}={\bf 1}$ and $s_{i-1}\notin {\mathcal
  J}$ we put $k_i=0$. We put \begin{gather}w_i=w_i(\lambda,{\mathcal
  J})=k_{i}{{{}}}+pk_{i-1}{{{}}}+p^2k_{i-2}{{{}}}+\ldots+p^ik_0{{{}}}+p^{i+1}k_d{{{}}}+\ldots+p^dk_{i+1}{{{}}},\label{wiede}\end{gather} $$\delta(\lambda,{\mathcal J})=(-1)^{d}\lambda(-{\rm id})\prod_{i=0}^dk_i!.$$Put ${\mathcal V}=\Theta_*{\mathcal V}_{M[\lambda,{\mathcal
  J},b]}$. The image of $e_i\in M[\lambda,{\mathcal
  J},b]$ in $H_0(\overline{\mathfrak X}_+,{\mathcal V})$ we denote again by
$e_i$.

\begin{pro}\label{abslang} (a) $H_0(\overline{\mathfrak X}_+,{\mathcal V})$ is a standard cyclic $k_{\mathcal E}^+[\varphi,\Gamma]$-module of perimeter $d+1$. The set $\{e_0,\ldots,e_d\}$ is a
$k$-basis of the kernel ${\rm ker}(t)$ of $t$, and we have the following formulae:\begin{align}t^{k_i{{{}}}}\varphi
  e_{i-1}&=k_i!\lambda^{[i]}(\tau(-1)) e_i&\mbox{ for }1\le i\le d,\label{monouni1}\\t^{k_0{{{}}}}\varphi
  e_{d}&=k_0!\lambda^{[0]}(\tau(-1))b^{-1}
  e_0,\label{monouni2}&{}\\t^{w_i}\varphi^{d+1}e_i&=(-1)^{d}\delta(\lambda,{\mathcal J})b^{-1}
  e_i&\mbox{ for }0\le i\le d,\label{idemmonouni}\\\gamma(x) e_i&=\lambda^{[i]}(\tau(x)) e_i&\mbox{ for }0\le i\le d\mbox{ and }x\in{\mathbb F}_p^{\times}.\label{auchdisplay}\end{align}
\end{pro}

{\sc Proof:} As $H_0(\overline{\mathfrak X}_+,{\mathcal V})$
is an inductive limit of ${{k_{\mathcal E}^+}}$-modules which are finite dimensional over $k$, it is a torsion ${{k_{\mathcal E}^+}}$-module. From Theorem \ref{finite} we obtain $M[\lambda,{\mathcal
  J},b]={\mathcal V}({\mathfrak
  e}_0)=H_0({\mathfrak
  X}_+,{\mathcal V})^{{\mathfrak N}_0}={\rm ker}(t)$. From this and the strict
level $1$ property (namely formula (\ref{striclevgen})) together with the fact
that for any ${\bf e}\in Y^1$ we find some $n\in {N}_0$ and $m\ge0$ with $n\phi^m{\bf e}_0={\bf e}$ it follows that ${\rm ker}(t)$ generates $H_0(\overline{\mathfrak X}_+,{\mathcal V})$. 

For $x\in{\mathbb F}_p^{\times}$ we compute
$T_{h_{s_d}(x)}e_i=\lambda^{[i]}(h_{s_d}(x^{-1}))e_i=x^{-k_i}e_i$. On the other hand, in Lemma
\ref{symcla} (with $k_i=r$ there) we have
$T_{{h_s}(x)}e=x^{-k_i}e$. Therefore, if we apply Lemma \ref{symcla} to
the character of ${\mathcal
  H}(\overline{\mathcal S},\overline{\mathcal U})_k$ obtained by pulling back the character $\chi_{\lambda^{[i]},{\mathcal
  J}^{[i]},b}$ of ${\mathcal H}(G,I_0)_{{\rm
    aff},k}$ along the embedding (\ref{heckeinbet}), then we obtain
$t^{k_i{{{}}}}n_{s_d}^{-1} e_{i}=k_i! e_i$ inside $M[\lambda,{\mathcal
  J},b]={\mathcal V}({\mathfrak
  e}_0)$, the latter viewed (cf. Proposition \ref{jalgheck}) as a subspace of$${\mathcal V}({\mathfrak v}_0)\cong ({\rm
  ind}^{\overline{\mathcal S}}_{\overline{\mathcal U}}{\bf 1}_k)\otimes_{{\mathcal
  H}(\overline{\mathcal S},\overline{\mathcal U})_k}M.$$A matrix computation shows $n_{s_d}s_d=\tau(-1)=\tau(-1)^{-1}\in\overline{T}$. We therefore see that $s_de_i=\lambda^{[i]}(\tau(-1))n_{s_d}^{-1}e_i$ since $\tau(-1)$ acts by $T_{\tau(-1)^{-1}}$, i.e. by $\lambda^{[i]}(\tau(-1)^{-1})$ on $k.e_i$. We get\begin{align}t^{k_i{{{}}}}s_d
  e_{i}&=\lambda^{[i]}(\tau(-1))t^{k_i{{{}}}}n_{s_d}^{-1}e_{i}\notag\\{}&=k_i!\lambda^{[i]}(\tau(-1))e_i.\label{proporz}\end{align} Formula
(\ref{hecnor}) gives $u e_{i-1}=T_{u^{-1}}e_{i-1}$. For $1\le i\le d$ we
therefore obtain\begin{gather}t^{k_i{{{}}}}\varphi e_{i-1}=t^{k_i{{{}}}}s_d u
e_{i-1}=t^{k_i{{{}}}}s_d
T_{u^{-1}}e_{i-1}=t^{k_i{{{}}}}s_de_i\label{fas1}\end{gather}while for $i=0$ we obtain\begin{gather}t^{k_0{{{}}}}\varphi
e_{d}=t^{k_0{{{}}}}s_d u e_{d}=t^{k_0{{{}}}}s_d
T_{u^{-1}}e_{d}=t^{k_0{{{}}}}s_d T_{u^{-1-d}}e_{0}=
t^{k_0{{{}}}}s_db^{-1}e_0.\label{fas2}\end{gather}Formula (\ref{proporz})
together with formula (\ref{fas1}), resp. with formula (\ref{fas2}), gives formula
(\ref{monouni1}), resp. formula (\ref{monouni2}). Combining formulae
(\ref{monouni1}) and (\ref{monouni2}) gives formula (\ref{idemmonouni}); for this
observe $\varphi t=t^p\varphi$ and $\lambda(-{\rm
  id})=\prod_{i=0}^d\lambda^{[i]}(\tau(-1))$. To see formula
(\ref{auchdisplay}) notice that the action of $\gamma(x)\in\Gamma$ is given by
the action of $\tau(x)\in {T}$, and that $\tau(x)$ acts by $T_{\tau(x)^{-1}}$, i.e. by $\lambda^{[i]}(\tau(x))$ on $k.e_i$.\hfill$\Box$\\

For $0\le j\le d+1$ let $i_j=p-1-k_{d+1-j}$ and $h_j=\sum_{i=0}^{j-1}i_{d+i+1-j}p^i$. (Attention: In general, $0=h_0$ need not be equal to $h_{d+1}$. The $h_j$
must not be confused with the $h_{s_j}$.) Define $0\le s\le p-2$ by the condition $x^{-s}=\lambda(\tau(x))$ for all $x\in{\mathbb F}_p^{\times}$. Let $\beta\in ({k}^{\rm alg})^{\times}$ be such that $\beta^{d+1}=\delta(\lambda,{\mathcal J})^{-1}b$.

\begin{satz}\label{langloc} $h=h_{d+1}/(p-1)$ is an integer, and if $k_i>0$ for some $i$ then we have an isomorphism of $k$-linear ${\rm Gal}_{{\mathbb Q}_p}$-representations$$W({\bf D}({\mathcal V}))\cong{\rm
  ind}(\omega_{d+1}^h)\otimes\omega^s\mu_{\beta}.$$  
\end{satz}

{\sc Proof:} That $h$ is an integer follows from Lemma \ref{normint}. Alternatively,
it follows from the divisibility of $\sum_{0\le i\le
 d}k_i{{{}}}$ by $p-1$, which is a consequence of$$x^{\sum_{i}k_i{{{}}}}=\prod_ix^{k_i{{{}}}}=\prod_i\lambda^{[i]}(h_{s_d}(x))=\lambda(\prod_ih_{s_i}(x))=\lambda({\rm
 id})=1$$for all $x\in{\mathbb F}_p^{\times}$. We now conclude with Proposition \ref{galstacom}, using the formulae of Proposition \ref{abslang}. (By Lemma \ref{normint}, the formula $x^{-s-h_j}=\lambda^{[j]}(\tau(x))$  for $x\in{\mathbb F}_p^{\times}$ and {\it all} $0\le
j\le d$ follows from the case $j=0$; alternatively it can be verified
by a straightforward calculation showing $x^{h_j}=(\lambda/\lambda^{[j]})(\tau(x))$.)\hfill$\Box$\\

\begin{satz}\label{h0emerton} (a) If for any $1\le j\le d$ there is some $0\le i\le
  d$ such that $k_i{{{}}}\ne k_{i+j}{{{}}}$ then $H_0(\overline{\mathfrak X}_+,{\mathcal V})$ is irreducible as a
  ${{k_{\mathcal E}^+}}[\varphi]$-module.

(b) If $(\lambda,{\mathcal
  J})$ satisfies (\ref{diff}) then $H_0(\overline{\mathfrak X}_+,{\mathcal V})$ is irreducible as a
  ${{k_{\mathcal E}^+}}[\varphi,\Gamma]$-module.

(c) If $(\lambda,{\mathcal
  J})$ satisfies (\ref{diff}) then $D({\mathcal V})$ is an irreducible $(\psi,\Gamma)$-module over ${{k_{\mathcal E}^+}}$, and ${\bf D}({\mathcal
  V})$ is an irreducible $(\varphi,\Gamma)$-module over ${{k_{\mathcal E}}}$. The integer $h=h_{d+1}/(p-1)$ is primitive and the ${\rm Gal}_{{\mathbb Q}_p}$-representation $W({\bf D}({\mathcal V}))\cong{\rm
  ind}(\omega_{d+1}^h)\otimes\omega^s\mu_{\beta}$ is irreducible. 
\end{satz}

{\sc Proof:} (In the case $d=1$ the argument was given in Theorem 5.1 of \cite{emert}.) We use the formulae in Proposition \ref{abslang}. Statement (a) follows
immediately from Proposition \ref{stnle}. 

For statement (b) we first claim that for any $1\le j\le d$ violating the hypothesis in (a), i.e. such that for all $0\le i\le d$ we have $k_i{{{}}}= k_{i+j}{{{}}}$, we have $\lambda^{[i]}\circ\tau\ne\lambda^{[i+j]}\circ\tau$ for all $0\le i\le d$.

Indeed, $k_i{{{}}}= k_{i+j}{{{}}}$ for all $0\le i\le d$ implies ${\mathcal
  J}^{[i]}={\mathcal
  J}^{[i+j]}$ for all $0\le i\le d$. Thus, by hypothesis (\ref{diff}) we have $\lambda^{[i]}\ne\lambda^{[i+j]}$ for all $0\le i\le d$. Now $k_i{{{}}}=
k_{i+j}{{{}}}$ for all $0\le i\le d$ says that the characters $\lambda^{[i]}$
and $\lambda^{[i+j]}$ of $\overline{T}$ differ at most by some character
invariant under conjugation by the full group $W$. But then we must have $\lambda^{[i]}\circ\tau\ne\lambda^{[i+j]}\circ\tau$ and the claim is proven.

Next, since $\gamma(x)\in \Gamma$ for $x\in{\mathbb F}_p^{\times}$ acts by multiplication with
$\lambda^{[i]}(\tau(x))$ on $k.e_i$ (formula (\ref{auchdisplay})) it follows that the hypotheses of Proposition \ref{stnle} are fulfilled in order to deduce the irreducibility of $H_0(\overline{\mathfrak X}_+,{\mathcal V})$ as a ${{k_{\mathcal E}^+}}[\varphi,\Gamma]$-module.

Statement (c): By Pontryagin duality theory, the natural map $H_0(\overline{\mathfrak
  X}_+,{\mathcal V})\to(H_0(\overline{\mathfrak X}_+,{\mathcal
  V})^*)^*=D({\mathcal
  V})^*$ is an isomorphism of (topological) $k$-vector spaces. It is checked that it is also an isomorphism of ${{k_{\mathcal E}^+}}[\varphi,\Gamma]$-modules. The same
check shows that the irreducibility of the ${{k_{\mathcal E}^+}}[\varphi,\Gamma]$-module $H_0(\overline{\mathfrak
  X}_+,{\mathcal V})$ implies the irreducibility of the $(\psi,\Gamma)$-module $D({\mathcal
  V})=H_0(\overline{\mathfrak
  X}_+,{\mathcal V})^*$, which in its turn implies the irreducibility of the $(\varphi,\Gamma)$-module ${\bf D}({\mathcal
  V})$ (cf. e.g. \cite{vien} Proposition 3.3.25). The irreducibility of ${\bf D}({\mathcal
  V})$ implies that of ${\rm
  ind}(\omega_{d+1}^h)\otimes\omega^s\mu_{\beta}$ (Theorem \ref{langloc}), and hence the primitivity of $h$ because only primitive $h$ give rise to irreducible ${\rm
  ind}(\omega_{d+1}^h)$.\hfill$\Box$\\

\begin{satz}\label{einerich} For any $0\le h<(p^{d+1}-1)/(p-1)$, any $0\le s\le p-2$ and any $\beta\in({k}^{\rm alg})^{\times}$ with $\beta^{d+1}\in k^{\times}$ there are $\lambda$, ${\mathcal
  J}$ and $b$ as before such that we have an isomorphism of irreducible $k$-linear ${\rm Gal}_{{\mathbb Q}_p}$-representations $W({\bf D}({\mathcal V}))\cong{\rm
  ind}(\omega_{d+1}^h)\otimes\omega^s\mu_{\beta}$ for ${\mathcal V}=\Theta_*{\mathcal V}_{M[\lambda,{\mathcal
  J},b]}$.

 If $h$ is primitive then $(\lambda,{\mathcal
  J})$ satisfies (\ref{diff}).   
\end{satz}

{\sc Proof:} Write $h(p-1)=i_0+pi_{1}+\ldots+p^di_d$ with $0\le i_j\le
p-1$, then put $k_{j}=p-1-i_{d+1-j}$ for $1\le j\le d+1$ and $k_0=k_{d+1}$. As $p-1$ divides $h(p-1)$ and hence $\sum_{j=0}^di_j$, it also divides
$\sum_{j=0}^dk_j$. As $\overline{T}$ is generated by the images of $\tau$ and
all the $h_{s_{j}}$, subject to the relation $\prod_{j=0}^dh_{s_{j}}=1$, it therefore follows that there exists a unique character
$\lambda:\overline{T}\to k^{\times}$ with
$\lambda^{[i]}(h_{s_d}(x))=\lambda(h_{s_{i-1}}(x))=x^{k_i}$ and with
$\lambda(\tau(x))=x^{-s}$ for all $x\in{\mathbb F}_p^{\times}$, all $0\le i\le
d$. Let ${\mathcal J}=\{s_i\in S\,|\,k_{i+1}=p-1\}$. Then ${\mathcal J}\subset
S_{\lambda}$. Let $b=\delta(\lambda,{\mathcal J})\beta^{d+1}\in k^{\times}$.

That $W({\bf D}({\mathcal V}))\cong{\rm
  ind}(\omega_{d+1}^h)\otimes\omega^s\mu_{\beta}$ for this triple $(\lambda,{\mathcal
  J},b)$ follows from Theorem \ref{langloc}. If $h$ is primitive then ${\rm
  ind}(\omega_{d+1}^h)\otimes\omega^s\mu_{\beta}$ is irreducible, hence ${\bf D}={\bf D}({\mathcal V})$ is irreducible. It follows that $D^{\sharp}$, the unique (by \cite{col} par. II.4 and II.5) non degenerate surjective $(\psi,\Gamma)$-module over $k_{\mathcal E}^+$ giving rise to ${\bf D}$, is irreducible. [Given a non-zero $(\psi,\Gamma)$-submodule $D_1^{\sharp}$ of $D^{\sharp}$ we
obtain a non-zero ${{k_{\mathcal E}}}$-sub vector space
$D_1^{\sharp}\otimes_{{{k_{\mathcal E}^+}}}{{k_{\mathcal E}}}$ of ${\bf D}$ stable under
$\varphi_{\bf D}$ and $\Gamma$. Thus $D_1^{\sharp}\otimes_{{{k_{\mathcal E}^+}}}{{k_{\mathcal E}}}={\bf D}$ by Proposition II.3.5 of \cite{col}, hence $D_1^{\sharp}=D_1^{\sharp}\otimes_{{{k_{\mathcal E}^+}}}{{k_{\mathcal E}}}\cap D^{\sharp}=D^{\sharp}$. See also \cite{vien} Proposition 3.3.26.] Thus $H_0(\overline{\mathfrak
  X}_+,{\mathcal V})^*$ and hence $H_0(\overline{\mathfrak
  X}_+,{\mathcal V})$ are irreducible. The irreducibility
of the
  ${{k_{\mathcal E}^+}}[\varphi,\Gamma]$-module $H_0(\overline{\mathfrak
  X}_+,{\mathcal V})$ implies the
irreducibility of the ${\mathcal H}(G,I_0)_k$-module $M[\lambda,{\mathcal J},b]$. Indeed, a proper ${\mathcal H}(G,I_0)_k$-sub module of $M[\lambda,{\mathcal J},b]$ would induce a proper sub coefficient system of ${\mathcal V}$, and Theorem \ref{finite} applied to these two coefficient systems would then induce a proper ${{k_{\mathcal E}^+}}[\varphi,\Gamma]$-sub module of $H_0(\overline{\mathfrak
  X}_+,{\mathcal V})$. Hence $(\lambda,{\mathcal
  J})$ satisfies (\ref{diff}).\hfill$\Box$\\

{\bf Remark:} Property (\ref{diff}) for the pair $(\lambda,{\mathcal
  J})$ in Theorem \ref{einerich}, and conversely, the primitivity of $h$ in Theorem
\ref{langloc} are proven here very indirectly. It looks quite cumbersome
trying to give direct, purely combinatorial proofs.\\

{\bf Remark:} The Hecke operator $T_{u^{d+1}}$ for $u^{d+1}\in\{p\cdot{\rm id},p^{-1}\cdot{\rm id}\}$ acts on
$M[\lambda,{\mathcal J},b]$ through the scalar $b$. The determinant of the
action of the geometric Frobenius on $W({\bf D}(\Theta_*{\mathcal
  V}_{M[\lambda,{\mathcal J},b]}))$ is $\delta(\lambda,{\mathcal
  J})^{-1}b$. If $d=1$ then $\delta(\lambda,{\mathcal J})=1$ in $k^{\times}$
so that this determinant is equal to $b$. If $d>1$ then
$\delta(\lambda,{\mathcal J})$ is {not} necessarily equal to $1$ in
$k^{\times}$ and is {not} even independent on $\lambda$.\\

Consider the composed functor \begin{gather}M\mapsto W({\bf
    D}(\Theta_*{\mathcal V}_M)).\label{hauptfu}\end{gather}from the category
${\rm Mod}^{\rm fin}({\mathcal
  H}(G,I_0)_{k})$
of finite dimensional ${\mathcal H}(G,I_0)_{k}$-modules into the category of
${\rm Gal}_{{\mathbb Q}_p}$-representations over $k$.

\begin{satz}\label{supersinghaupt} The functor (\ref{hauptfu}) induces a bijection between 

(a) the set of isomorphism classes of absolutely simple supersingular
${\mathcal H}(G,I_0)_k$-modules of dimension $d+1$ and 

(b) the set of isomorphism classes of smooth irreducible representations of ${\rm Gal}_{{\mathbb Q}_p}$ over $k$ of dimension $d+1$.
\end{satz}

{\sc Proof:} By \cite{berger} Corollary 2.1.5, the discussion preceeding
\cite{berger} Lemma 2.2.1 and the beginning of \cite{berger} section 3.1, any smooth irreducible representation of ${\rm Gal}_{{\mathbb Q}_p}$ over $k$
of dimension $d+1$ is of the form ${\rm
  ind}(\omega_{d+1}^h)\otimes\omega^s\mu_{\beta}$ with $1\le h\le(p^{d+1}-1)/(p-1)-1$. Therefore it is isomorphic with $W({\bf D}(\Theta_*{\mathcal V}_M))$ for an absolutely simple supersingular
${\mathcal H}(G,I_0)_k$-module $M$ of dimension $d+1$, by Theorem \ref{einerich}.

 By Proposition \ref{marfrarac} any absolutely simple supersingular ${\mathcal H}(G,I_0)_k$-module $M$ of dimension $d+1$ is of the form $M[\lambda,{\mathcal
  J},b]$ with $(\lambda,{\mathcal
  J})$ satisfying (\ref{diff}). It remains to show that the isomorphism class
of the simple ${\mathcal H}(G,I_0)_k$-module $M[\lambda,{\mathcal
  J},b]$ can be recovered from the isomorphism class of ${\bf D}({\mathcal V})$ for ${\mathcal
  V}=\Theta_*{\mathcal
  V}_{M[\lambda,{\mathcal J},b]}$.

{\it Step 1: The isomorphism class of $M[\lambda,{\mathcal
  J},b]$ can be recovered from the isomorphism class of $H_0(\overline{\mathfrak X}_+,{\mathcal V})$ as a
  ${{k_{\mathcal E}^+}}[\varphi,\Gamma]$-module.}

 Recall the $w_i$, formula (\ref{wiede}). Let $\check{w}={\rm min}\{w_i\,|\,0\le i\le d\}$. For $w\in{\mathbb Z}_{\ge0}$ let $F_w={\rm
  ker}(t)\cap {\rm
  ker}(t^{w+1}\varphi^{d+1})$. By formula
(\ref{idemmonouni}) we may recover $\check{w}$ as $\check{w}={\rm
  min}\{w\in{\mathbb Z}_{\ge0}\,|\,F_w\ne0\}$. Let
$I=\{i\,|\,w_i=\check{w}\}=\{i\,|\,F_{w_i}=F_{\check{w}}\}$. Let $\overline{T}_0=\overline{T}\cap{\rm
  SL}_{d+1}({\mathbb F}_p)$. For $i\in I$ we claim that we can recover
$\lambda^{[i]}|_{\overline{T}_0}$ and that this is the same for all $i\in I$. Indeed,
$\overline{T}_0$ is generated by the $u^{-m}h_{s_d}(x)u^m$ for $1\le m\le d+1$
and $x\in{\mathbb F}_p$, we have
$\lambda^{[i]}(u^{-m}h_{s_d}(x)u^m)=\lambda^{[i+m]}(h_{s_d}(x))=x^{k_{i+m}}$
and $k_{i+m}$ can be read off from $\check{w}$ as the coefficient of
$p^{d+1-m}$, by formula (\ref{wiede}).

Next, pick some eigenvector $e\in F_{\check{w}}$ for the action of $\Gamma$ on
$F_{\check{w}}$. By formula (\ref{auchdisplay}) there is some $i\in I$
such that the action of $\Gamma$ on $k.e$ is given by
$\lambda^{[i]}\circ\tau$. With the above this allows us to recover
$\lambda^{[i]}$ because the domain $\overline{T}$ of $\lambda^{[i]}$ is
generated by $\overline{T}_0$ and the image of $\tau$. As a result we see that we can recover $\lambda$ up to conjugation
by a power of $u$. Next, up to cyclic permutation we can recover the $k_i$ as the
digits in the expansion of $\check{w}$ in base $p$, as already remarked. Since these determine the ${\mathcal
  J}^{[i]}$ we see that together with what has been said we can recover the pairs $(\lambda^{[i]},{\mathcal
  J}^{[i]})$ up to cyclic permutation, or what is the same, we can recover the pair $(\lambda,{\mathcal
  J})$ up to conjugation by a power of $u$. Finally, knowing the $k_i$ (up to
cyclic permutation) allows us to recover $\prod_{i=0}^dk_i!$. Moreover,
knowing $\lambda$ (up to conjugation by a power of $u$) allows us to recover
$\lambda(-{\rm id})=\prod_{i=0}^d\lambda^{[i]}(\tau(-1))$; namely, this is the product of the
eigenvalues of $\tau(-1)$ acting on ${\rm
  ker}(t)$ through the action of $\Gamma$ (formula (\ref{auchdisplay})). Hence we may recover $\delta(\lambda,{\mathcal J})$. But
then we may recover $b$ from formula (\ref{idemmonouni}) which holds true with
$w_i$ replaced by $\check{w}$ and with $e_i$ replaced by any $e\in F_{\check{w}}$.

We are done because the isomorphism class of the ${\mathcal H}(G,I_0)_k$-module $M[\lambda,{\mathcal
  J},b]$ is given by $b$ and the pair $(\lambda,{\mathcal
  J})$, the latter taken up to conjugation by powers of $u$.

{\it Step 2: The isomorphism class of $H_0(\overline{\mathfrak
  X}_+,{\mathcal V})$ can be recovered from ${\bf D}({\mathcal
  V})$. }

The $(\varphi,\Gamma)$-module ${\bf D}({\mathcal
  V})$ is irreducible (cf. e.g. \cite{vien} Proposition 3.3.25). Therefore (and because ${\bf D}({\mathcal
  V})$ is of dimension $>1$) ${\bf D}({\mathcal
  V})$ contains a unique compatible $(\psi,\Gamma)$-submodule over ${{k_{\mathcal E}^+}}$ on which the $\psi$-operator is surjective (\cite{col} par. II.4 and II.5). By uniqueness, this must be $D({\mathcal
  V})$. Thus, the $(\psi,\Gamma)$-module $D({\mathcal
  V})$ can be recovered from ${\bf D}({\mathcal
  V})$. Since $H_0(\overline{\mathfrak
  X}_+,{\mathcal V})\to(H_0(\overline{\mathfrak X}_+,{\mathcal
  V})^*)^*=D({\mathcal
  V})^*$ is an isomorphism of ${{k_{\mathcal E}^+}}[\varphi,\Gamma]$-modules we are done.\hfill$\Box$\\

{\bf Remark:} A numerical version of Theorem \ref{supersinghaupt} was
proven in \cite{vigneras} Theorem 5: there it
was shown that for fixed $b$, the number of absolutely simple supersingular ${\mathcal H}(G,I_0)_k$-modules of the form $M[\lambda,{\mathcal
  J},b]$ is the same as the number of smooth irreducible representations of ${\rm Gal}_{{\mathbb Q}_p}$ over $k$ of dimension $d+1$ with a fixed determinant of the Frobenius.\\

{\bf Remark:} ${\mathcal H}(G,I_0)_k$-modules of the form $M[\lambda,{\mathcal
  J},b]$ with $(\lambda,{\mathcal
  J})$ {\it violating} (\ref{diff}) contain supersingular ${\mathcal
  H}(G,I_0)_k$-modules $M$ of dimensions $l$ which are proper divisors of
$d+1$ (and all supersingular ${\mathcal
  H}(G,I_0)_k$-modules arise in this way, by \cite{oll}). The associated ${\rm Gal}_{{\mathbb Q}_p}$-representations $W({\bf
    D}(\Theta_*{\mathcal V}_M))$ are also of dimension $l$,
  irreducible for simple $M$. One may either detect them inside $W({\bf
    D}(\Theta_*{\mathcal V}_{M[\lambda,{\mathcal
  J},b]}))$ (using Theorem \ref{langloc} and the exactness of our functor). Alternatively one may observe that $M$ is structured in the same
fashion as $M[\lambda,{\mathcal
  J},b]$ (but of 'perimeter' $l$ instead of $d+1$) (simplicity for $M$ can be
characterized by a non-periodicity property analogous to (\ref{diff})), hence one can compute $W({\bf
    D}(\Theta_*{\mathcal V}_M))$ and test its irreducibility
exactly along the lines we did with $W({\bf
    D}(\Theta_*{\mathcal V}_{M[\lambda,{\mathcal
  J},b]}))$.

\subsection{Filtrations on the Weyl group}

\label{filtrsec}

Let $\ell:W\to{\mathbb Z}_{\ge0}$ be the length function with respect
to the set of Coxeter generators $\{s_1,\ldots,s_d\}$ of $W$. Setting
$\overline{u}=s_d\cdots s_1$ we have
$s_i=\overline{u}^{d-i}s_d\overline{u}^{i-d}$ for $1\le i\le
d$. 

Let $W^{s_d}=\{w\in W\,|\,\ell(ws_d)>\ell(w)\}$. Suppose that we are given a map
$\sigma:W^{s_d}\to\{-1,0,1\}$. In the following, for $w\in W$ and $i\in\{-1,-0,1\}$ we write
$\sigma(w)=i$ as a shorthand for [$w\in W^{s_d}$ and $\sigma(w)=i$]. \\

{\bf Definition:} (a) For a subset $W'$ of $W$ we define a self mapping$$(.)_+^{W'}:W\longrightarrow W,\quad w \mapsto w_+^{W'}$$ by the formula\begin{align}w_+^{W'}&=\left\{\begin{array}{l@{\quad:\quad}l}w\overline{u}^{-1}&
      \sigma(w\overline{u}^{-1})=0\mbox{ or
      }\sigma(w\overline{u}^{-1}s_d)=0 \\{}&[\sigma(w\overline{u}^{-1})=-1\mbox{ or
  }\sigma(w\overline{u}^{-1}s_d)=1]\mbox{ and }w\overline{u}^{-1}s_d\overline{u}\notin
  W'\\{}&[\sigma(w\overline{u}^{-1})=1\mbox{
    or }\sigma(w\overline{u}^{-1}s_d)=-1]\mbox{ and }w\overline{u}^{-1}s_d\overline{u}\in W'\\w\overline{u}^{-1}s_d  &[\sigma(w\overline{u}^{-1})=1\mbox{
    or }\sigma(w\overline{u}^{-1}s_d)=-1]\mbox{ and }w\overline{u}^{-1}s_d\overline{u}\notin W'\\{}&[\sigma(w\overline{u}^{-1})=-1\mbox{
    or }\sigma(w\overline{u}^{-1}s_d)=1]\mbox{ and }w\overline{u}^{-1}s_d\overline{u}\in
  W'.\end{array}\right.\end{align}

(b) A $\sigma$-admissible filtration of $W$ is
a filtration by subsets\begin{gather}\emptyset=W_0\subset W_1\subset
W_2\subset\ldots\subset
W_{q}=W\label{wfi}\end{gather}of $W$ (for some $q\in{\mathbb N}$) such that
for each $1\le i\le q$ the map $(.)_+^{W_{i-1}}$ respects both $W-W_{i-1}$
and
$W_i-W_{i-1}$ and is bijective on the latter (i.e. restricts to
a permutation of $W_i-W_{i-1}$).\\

\begin{pro}\label{osnabrueckwart} $\sigma$-admissible filtrations exist for any $\sigma$.
\end{pro}

{\sc Proof:} Suppose we are given a sequence of subsets$$\emptyset=W_0\subset
W_1\subset W_2\subset\ldots\subset W_{i-1}\subset W$$(some $i\ge1$) such that for all $j<i$
and all $x\in W_j-W_{j-1}$ there is a $y\in W_j-W_{j-1}$ with
$x=y_+^{W_{j-1}}$.

{\it Claim: For any $w\in W-W_{i-1}$ we have $w_+^{W_{i-1}}\in W-W_{i-1}$.}

Assume that, on the contrary, $w_+^{W_{i-1}}\in W_{i-1}$,
i.e. $w_+^{W_{i-1}}\in W_{j}-W_{j-1}$ for some $j<i$. Thus
$w_+^{W_{i-1}}=v_+^{W_{j-1}}$ for some $v\in W_{j}-W_{j-1}$. As $w\notin
W_{i-1}$ we have $w\ne v$ and this forces that
either\begin{gather}w\overline{u}^{-1}=w_+^{W_{i-1}}=v_+^{W_{j-1}}=v\overline{u}^{-1}s_d\label{casea}\end{gather}or\begin{gather}w\overline{u}^{-1}s_d=w_+^{W_{i-1}}=v_+^{W_{j-1}}=v\overline{u}^{-1}.\label{caseb}\end{gather}and
moreover that neither $\sigma(w\overline{u}^{-1})=0$ nor
$\sigma(w\overline{u}^{-1}s_d)=0$.  

Suppose first that [$\sigma(w\overline{u}^{-1}s_d)=-1$ or $\sigma(w\overline{u}^{-1})=1$]
and we are in case (\ref{casea}). Then [$\sigma(v\overline{u}^{-1})=-1$ or
$\sigma(v\overline{u}^{-1}s_d)=1$] and by the definition of $v_+^{W_{j-1}}$ we
get $v\overline{u}^{-1}s_d\overline{u}\in W_{j-1}$, i.e. $w\in W_{j-1}$, contradiction.

Suppose next that [$\sigma(w\overline{u}^{-1}s_d)=1$ or $\sigma(w\overline{u}^{-1})=-1$]
and we are in case (\ref{caseb}). Then [$\sigma(v\overline{u}^{-1})=1$ or
$\sigma(v\overline{u}^{-1}s_d)=-1$] and by the definition of $v_+^{W_{j-1}}$ we
get $v\overline{u}^{-1}s_d\overline{u}\in W_{j-1}$, i.e. $w\in W_{j-1}$, contradiction.

Now suppose that [$\sigma(w\overline{u}^{-1}s_d)=1$ or $\sigma(w\overline{u}^{-1})=-1$]
and we are in case (\ref{casea}). Then $w\overline{u}^{-1}s_d\overline{u}\notin W_{i-1}$ by
the definition of $w_+^{W_{i-1}}$, i.e. $v\notin W_{i-1}$, contradiction. 

Finally, suppose that [$\sigma(w\overline{u}^{-1}s_d)=-1$ or $\sigma(w\overline{u}^{-1})=1$]
and we are in case (\ref{caseb}). Then $w\overline{u}^{-1}s_d\overline{u}\notin W_{i-1}$ by
the definition of $w_+^{W_{i-1}}$, i.e. $v\notin W_{i-1}$, contradiction. 
 
The claim is proven. It implies that the map $(.)^{W_{i-1}}_+$ respects $W-W_{i-1}$. Let $W_i-W_{i-1}$
be the maximal subset of $W-W_{i-1}$ on which $(.)^{W_{i-1}}_+$ is bijective; as $W-W_{i-1}$ is a finite set, $W_i-W_{i-1}$ is non-empty. Put
$W_i=(W_i-W_{i-1})\cup W_{i-1}$.

It is clear that the sequence $\emptyset=W_0\subset
W_1\subset W_2\subset\ldots\subset W_{i}\subset W$ again satisfies the above
condition. Proceeding inductively we thus obtain an admissible filtration of $W$.\hfill$\Box$\\

{\bf Remark:} The specific filtration constructed in the proof of
Proposition \ref{osnabrueckwart} may be regarded as the 'socle' filtration.\\

{\bf Examples:} (a) If $\sigma={\bf 0}$, i.e. if $\sigma(w)=0$ for any $w\in W^{s_d}$, then $\emptyset=W_0\subset
W_1=W$ is a $\sigma$-admissible filtration: the permutation
$(.)^{\emptyset}_+$ on $W$ is right multiplication by $\overline{u}^{-1}$.

(b) In Theorem \ref{sigma1fil} below we describe an explicit $\sigma$-admissible
filtration (in fact: the above 'socle' filtration) for $\sigma={\bf 1}$, i.e. when $\sigma(w)=1$ for any $w\in
W^{s_d}$.

(c) Assume $d=2$ so that $W=\{1,s_1,s_2,s_1s_2,s_2s_1,w_0\}$ with
$w_0=s_1s_2s_1=s_2s_1s_2$, and $W^{s_d}=W^{s_2}=\{1,s_1,s_2s_1\}$.

(c1) Assume $\sigma(s_2s_1)=\sigma(s_1)=1$ and $\sigma(1)=0$. Then
$(.)^{\emptyset}_+$ acts on $W_1=\{s_1s_2,w_0\}$ as $w_0\mapsto s_1s_2\mapsto
w_0$ while $(.)^{W_1}_+$ acts on $W-W_1$ as
$s_2s_1\mapsto 1\mapsto s_1\mapsto s_2\mapsto s_2s_1$.

(c2) Assume $\sigma(s_2s_1)=0$ and $\sigma(s_1)=\sigma(1)=1$. Then
$(.)^{\emptyset}_+$ acts on $W_1=\{s_2,s_2s_1,s_1s_2,w_0\}$ as $w_0\mapsto
s_1s_2\mapsto s_2s_1\mapsto s_2\mapsto w_0$ while $(.)^{W_1}_+$ acts on $W-W_1$ as
$s_1\mapsto 1\mapsto s_1$.\\

We assume that $d>1$. Let $\langle s_1,\ldots,s_{d-1}\rangle$ denote the subgroup of $W$ generated
by the elements $s_1,\ldots,s_{d-1}$, and similarly define ${\mathcal R}=\langle s_1,\ldots,s_{d-2}\rangle$. Put $W_0=\emptyset$ and
$W_{d+1}=W$, and for $1\le i\le d$ define the union of
cosets$$W_{i}=\bigcup_{1\le j\le i}s_j\cdots s_d\langle s_1,\ldots,s_{d-1}\rangle.$$

\begin{lem}\label{symgp}

(a) For any $w\in W_{i}\cap W^{s_d}$ we have $w\overline{u}\in W_{i-1}$.

(b) For any $w\in (W-W_{i-1})\cap W^{s_d}s_d$ we have $w\overline{u}\notin W_{i}$.

(c) For any $w\in W_{i}$ we have $|\{0\le j\le
d-1\,;\,w(\overline{u}^{-1}s_d)^j\in W^{s_d}\}|=i-1$. 

(d) The set ${\mathcal R}$ is a set of representatives in $\langle s_1,\ldots,s_{d-1}\rangle$ for the right cosets of the cyclic subgroup $(\overline{u}^{-1}s_d)^{\mathbb Z}=\{(\overline{u}^{-1}s_d)^j\,|\,j=0,\ldots, d-1\}$ of $\langle s_1,\ldots,s_{d-1}\rangle$, i.e. $\langle s_1,\ldots,s_{d-1}\rangle=\coprod_{v\in {\mathcal R}}v(\overline{u}^{-1}s_d)^{\mathbb Z}$. 
\end{lem} 

{\sc Proof:} Identify $W$ with the group of permutations of the set
$\{0,\ldots,d\}$, in such a way that $s_j=(j-1,j)$ (transposition) for $1\le
j\le d$. Then $W_{i}=\{w\in W\,|\,w(d)\le i\}$ and $W^{s_d}=\{w\in
W\,|\,w(d-1)<w(d)\}$, moreover $\overline{u}(j)=j-1$ for $1\le j\le d$ and
$\overline{u}(0)=d$. Statements (a), (b) and (c) are now easily read off. Statement (d) is clear.\hfill$\Box$\\

For $0\le j\le d-1$, for $1\le i\le d+1$ and for $v\in{\mathcal R}$ let $$w_j(i,v)=s_i\cdots
s_dv(\overline{u}^{-1}s_d)^j\in W$$and let $w_d(i,v)=w_0(i,v)$. Let $W_{i,v}=\{w_0(i,v),\ldots,w_{d-1}(i,v)\}$.

\begin{satz}\label{sigma1fil} If $\sigma={\bf 1}$ then the filtration $\emptyset=W_0\subset\ldots\subset W_{d+1}=W$ defined above is
$\sigma$-admissible. More precisely, for any $1\le i\le d+1$ we
have\begin{gather}W_{i}-W_{i-1}=\coprod_{v\in{\mathcal R}}W_{i,v}\label{klarfil}\end{gather}and
$w_j(i,v)_+^{W_{i-1}}=w_{j+1}(i,v)$ for all $0\le j\le d-1$. We have \begin{gather}|\{0\le j\le
d-1\,;\,w_{j}(i,v)\in W^{s_d}\}|=i-1.\label{carw}\end{gather}   
\end{satz}

{\sc Proof:} Lemma \ref{symgp} (d) implies the decomposition
(\ref{klarfil}). To see that $w_j(i,v)_+^{W_{i-1}}=w_{j+1}(i,v)$ assume first that
we are in the case where $w_j(i,v)\overline{u}^{-1}s_d\in W^{s_d}$. Then, as
$w_j(i,v)\overline{u}^{-1}s_d\in W_{i}$, Lemma \ref{symgp} (a) tells us
$w_j(i,v)\overline{u}^{-1}s_d\overline{u}\in W_{i-1}$, and the hypothesis $\sigma={\bf 1}$ therefore gives
$w_j(i,v)_+^{W_{i-1}}=w_{j+1}(i,v)$. If on the other hand
$w_j(i,v)\overline{u}^{-1}\in W^{s_d}$ then, as
$w_j(i,v)\overline{u}^{-1}s_d\notin W_{i-1}$, Lemma \ref{symgp} (b) tells us $w_j(i,v)\overline{u}^{-1}s_d\overline{u}\notin W_{i-1}$, and again the hypothesis $\sigma={\bf 1}$ gives
$w_j(i,v)_+^{W_{i-1}}=w_{j+1}(i,v)$. Finally, formula (\ref{carw})
follows from Lemma \ref{symgp} (c). \hfill$\Box$\\

\subsection{Reduced standard ${\mathcal H}(G,I_0)_k$-modules}
\label{unprse}

The justification for defining reduced standard modules over ${\mathcal H}(G,I_0)_k$ as we do it below is given in \cite{wty}.

$\overline{u}=s_d\cdots s_1$ is the image of $u\in
N(T)$ in $W=N(T)/T$. Our specific choices in (\ref{exam}) resp. (\ref{exam1}) allow us to regard $W$ as the subgroup of $N(T)$
generated by $s_1,\ldots,s_d$. \\

{\bf Definition:} We say that an ${\mathcal H}(G,I_0)_k$-module
$M$ is a {\it reduced standard module} (or: is of {\it of $W$-type}) if
it is of the following form $M=M(\theta,\sigma,\epsilon_{\bullet})$. First, a $k$-vector space basis of $M$
is the set of formal symbols $g_w$ for $w\in W$. The ${\mathcal
  H}(G,I_0)_k$-action on $M$ is characterized by a character
$\theta:\overline{T}\to k^{\times}$ (which we also read as a character of $T\cap I$ by inflation), a map $\sigma:W^{s_d}\to\{-1,0,1\}$ and a
set $\epsilon_{\bullet}=\{\epsilon_{w}\}_{w\in W}$ of units $\epsilon_w\in k^{\times}$. Namely, for $w\in W$ we define
$\kappa_w=\kappa_w(\theta)=\theta(wn_{s_d}s_dw^{-1})\in\{\pm 1\}$. Then the following
formulae are required for $t\in\overline{T}$ and $w\in W$: \begin{align}T_{t}(g_w)&=\theta(wt^{-1}w^{-1})g_{w},\notag\\T_{u^{-1}}(g_w)&=\epsilon_wg_{w\overline{u}^{-1}},\notag\end{align}$$T_{s_d}(g_w)=\left\{\begin{array}{l@{\quad:\quad}l}g_{ws_d}&\quad
      \mbox{  }[\sigma(ws_d)=-1\mbox{ and }{\theta}(wh_{s_d}w^{-1})\ne{\bf 1}]\mbox{ or }\sigma(w)=1       
\\-\kappa_wg_{w}&\quad \mbox{
}\sigma(ws_d)\in\{0,1\}\mbox{ and }{\theta}(wh_{s_d}w^{-1})={\bf 1}\\
g_{ws_d}-\kappa_wg_{w}&\quad \mbox{
}\sigma(ws_d)=-1\mbox{ and }{\theta}(wh_{s_d}w^{-1})={\bf
  1}\\0&\quad \mbox{
}\mbox{ all other cases }\end{array}\right.$$Here the conditions involving $\theta(wh_{s_d}w^{-1})=\theta(wh_{s_d}(.)w^{-1})$ and ${\bf 1}$
compare the homomorphism ${\mathbb F}_p^{\times}\to k^{\times}$, $x\mapsto \theta(wh_{s_d}(x)w^{-1})$
with the constant homomorphism $x\mapsto {\bf 1}(x)=1$.\\

Fix a reduced standard ${\mathcal H}(G,I_0)_k$-module $M=M(\theta,\sigma,\epsilon_{\bullet})$ and let ${\mathcal V}=\Theta_*{\mathcal V}_{M(\theta,\sigma,\epsilon_{\bullet})}$. Choose 
a $\sigma$-admissible filtration $\emptyset=W_0\subset W_1\subset\ldots\subset
W_q=W$. For $1\le i\le q$ decompose $W_i-W_{i-1}$ into its
$(.)_+^{W_{i-1}}$-orbits, i.e. write$$W_i-W_{i-1}=\coprod_{v\in{\mathcal
  R}_i}\{w_0(i,v), w_1(i,v),\ldots, w_{t_v-1}(i,v)\}$$where $t_v$ is the
cardinality of the $(.)_+^{W_{i-1}}$-orbit with index $v\in {\mathcal R}_i$, such that
$w_j(i,v)_+^{W_{i-1}}=w_{j+1}(i,v)$ for all $0\le j\le t_v-1$; here we write
$w_{t_v}(i,v)=w_{0}(i,v)$.

For $v\in{\mathcal
  R}_i$ let $N_{i,v}$ be the $k_{\mathcal
    E}^+[\varphi,\Gamma]$-submodule of $H_0(\overline{\mathfrak X}_+,{\mathcal
    V})$ generated by the elements $g_w\in
  M\subset H_0(\overline{\mathfrak
    X}_+,{\mathcal V})$ for all $w\in W_{i-1}\cup\{w_0(i,v), w_1(i,v),\ldots,
  w_{t_v-1}(i,v)\}$. Let $N_i=\sum_{v\in{\mathcal
  R}_i}N_{i,v}$, i.e. the $k_{\mathcal
    E}^+[\varphi,\Gamma]$-submodule of $H_0(\overline{\mathfrak X}_+,{\mathcal
    V})$ generated by the elements $g_w\in
  M\subset H_0(\overline{\mathfrak
    X}_+,{\mathcal V})$ for all $w\in W_{i}$.

\begin{satz}\label{newstafi} Assume that there is no
$v\in W^{s_d}$ with $\sigma(v)=-1$ and $\theta(vh_{s_d}v^{-1})={\bf
  1}$. Then$$\emptyset=N_0\subset N_1\subset\ldots\subset N_q=H_0(\overline{\mathfrak
    X}_+,{\mathcal V})$$is a filtration by $k_{\mathcal
    E}^+[\varphi,\Gamma]$-submodules, its subquotients admit direct sum decompositions$$N_i/N_{i-1}=\bigoplus_{v\in{\mathcal
  R}_i}N_{i,v}/N_{i-1}$$with $N_{i,v}/N_{i-1}$ standard cyclic of
perimeter $t_v$.

The ${\rm Gal}_{{\mathbb Q}_p}$-representation $W({\bf
  D}({\mathcal V}))$ admits a filtration $W({\bf
  D}({\mathcal V}))=V_0\supset V_1\supset\ldots\supset V_q=0$ and direct sum
decompositions $V_{i-1}/V_i=\oplus_{v\in{\mathcal
  R}_i}\overline{V}_{i-1,v}$ such that $\overline{V}_{i-1,v}$ is the ${\rm
Gal}_{{\mathbb Q}_p}$-representation corresponding to the $k_{\mathcal
    E}^+[\varphi,\Gamma]$-module $N_{i,v}/N_{i-1}$. 
\end{satz}

The proof of Theorem \ref{newstafi} is based on the following Lemma
\ref{fmlueber}.

Let $w\in W$. If $\theta(wh_{s_d}^{-1}w^{-1})={\bf
  1}$ we put $\lfloor k(w)\rfloor=0$ and $\lceil k(w)\rceil=p-1$ (but do
not define $k(w)$). If $\theta(wh_{s_d}^{-1}w^{-1})\ne{\bf
  1}$ we define $1\le k(w)\le p-2$ by
$\theta(wh_{s_d}^{-1}(x)w^{-1})=x^{-k(w)}$ for all $x\in{\mathbb F}_p^{\times}$, and we put $\lfloor k(w)\rfloor=\lceil k(w)\rceil=k(w)$.

\begin{lem}\label{fmlueber} Let $v\in W^{s_d}$. In $H_0(\overline{\mathfrak
    X}_+,{\mathcal V})$ we have the following identities.\\If $\sigma(v)=1$ then\begin{align}t^{p-1}n_{s_d}^{-1}g_v&=\kappa_{v}g_{vs_d}\label{ly1}&{}\\t^{\lceil k(vs_d)\rceil}n_{s_d}^{-1}g_{vs_d}&=\lceil k(vs_d)\rceil!g_{vs_d}\label{ly2}&{}\\n_{s_d}^{-1}\kappa_{v}g_{vs_d}&-g_v&\in\sum_{i\ge0}k.t^i{s_d}g_v\label{ly3}\\t^{p-1-\lceil k(vs_d)\rceil}n_{s_d}^{-1}g_v&-(p-1-\lceil k(vs_d)\rceil)!g_v&\in\sum_{i\ge0}k.t^is_dg_{vs_d}.\label{ly4}\end{align} If $\sigma(v)=0$ then\begin{gather}t^{\lfloor k(v)\rfloor}n_{s_d}^{-1}g_v=\lfloor
    k(v)\rfloor!g_{v}\quad\quad\quad\mbox{ and }\quad\quad\quad t^{\lceil
      k(vs_d)\rceil}n_{s_d}^{-1}g_{vs_d}=\lceil
    k(vs_d)\rceil!g_{vs_d}.\label{ly5}\end{gather}If $\sigma(v)=-1$ and $\theta(vh_{s_d}v^{-1})\ne{\bf 1}$ then\begin{align}t^{p-1}n_{s_d}^{-1}g_{vs_d}&=\kappa_{vs_d}g_{v}\label{ly6}&{}\\t^{\lceil
      k(v)\rceil}n_{s_d}^{-1}g_{v}&=\lceil
    k(v)\rceil!g_{v}&{}\label{ly7}\\n_{s_d}^{-1}\kappa_{vs_d}g_{v}
    &-g_{vs_d}&\in\sum_{i\ge0}k.t^i{s_d}g_{vs_d}\label{ly8}\\t^{p-1-\lceil
      k(v)\rceil}n_{s_d}^{-1}g_{vs_d}&-(p-1-\lceil
    k(v)\rceil)!g_{vs_d}&\in\sum_{i\ge0}k.t^is_dg_{v}\label{ly9}.\end{align}
\end{lem}

{\sc Proof:} We view $M$ as an ${\mathcal H}(\overline{\mathcal
  S},\overline{\mathcal U})_k$-module by means of the embedding
(\ref{heckeinbet}). As in the proof of Proposition \ref{abslang} we use the corresponding embedding\begin{gather}({\rm
  ind}^{\overline{\mathcal S}}_{\overline{\mathcal U}}{\bf 1}_k)\otimes_{{\mathcal
  H}(\overline{\mathcal S},\overline{\mathcal U})_k}M\cong{\mathcal V}({\mathfrak v}_0)\hookrightarrow H_0(\overline{\mathfrak X}_+,{\mathcal V}).\label{locanal}\end{gather}Regarded as an ${\mathcal H}(\overline{\mathcal
  S},\overline{\mathcal U})_k$-module, $M$ is the
direct sum, indexed by all $v\in
  W^{s_d}$, of the two dimensional ${\mathcal H}(\overline{\mathcal
  S},\overline{\mathcal U})_k$-modules $\langle g_v,g_{vs_d}\rangle$. If $\sigma(v)=0$ this summand splits up further as$$\chi_{\lfloor k(v)\rfloor}\bigoplus\chi_{\lceil
  k(vs_d)\rceil}\cong \langle g_v,g_{vs_d}\rangle$$with the ${\mathcal H}(\overline{\mathcal
  S},\overline{\mathcal U})_k$-modules $\chi_{\lfloor k(v)\rfloor}$, $\chi_{\lceil
  k(vs_d)\rceil}$ considered in Lemma \ref{symcla}. Namely, send $e\in\chi_{\lfloor k(v)\rfloor}$ to $g_v$ and send $e\in\chi_{\lceil
  k(vs_d)\rceil}$ to $g_{vs_d}$. Thus, if $\sigma(v)=0$ we conclude with Lemma
\ref{symcla}. Next, if $\sigma(v)=1$ we have an isomorphism of ${\mathcal H}(\overline{\mathcal
  S},\overline{\mathcal U})_k$-modules$$M_{\lceil
  k(vs_d)\rceil}\cong \langle g_v,g_{vs_d}\rangle$$sending $e\in M_{\lceil
  k(vs_d)\rceil}$ to $g_v$ and sending $f\in M_{\lceil
  k(vs_d)\rceil}$ to
$T_{n_{s_d}}g_{v}=T_{{s_d}}T_{n_{s_d}s_d}g_{v}=\kappa_{v}g_{vs_d}$. Thus,
in this case we conclude with Lemma \ref{sl2fpind}. Finally, if $\sigma(v)=-1$ we have an isomorphism of ${\mathcal H}(\overline{\mathcal
  S},\overline{\mathcal U})_k$-modules$$M_{\lceil
  k(v)\rceil}\cong \langle g_v,g_{vs_d}\rangle$$sending $e\in M_{\lceil
  k(v)\rceil}$ to $g_{vs_d}$ and sending $f\in M_{\lceil
  k(v)\rceil}$ to $T_{n_{s_d}}g_{vs_d}=T_{{s_d}}T_{n_{s_d}s_d}g_{vs_d}$. The
latter equals $\kappa_{vs_d}g_v$ as we assume $\theta(vh_{s_d}v^{-1})\ne{\bf 1}$
(it equals
$\kappa_{vs_d}g_v-g_{vs_d}$ if $\theta(vh_{s_d}v^{-1})={\bf 1}$). Thus, also in
this case we conclude with Lemma \ref{sl2fpind}.\hfill$\Box$\\ 

{\sc Proof of Theorem \ref{newstafi}:} As in Proposition \ref{abslang} we see that ${\rm ker}(t|_{H_0(\overline{\mathfrak X}_+,{\mathcal V})})$ generates $H_0(\overline{\mathfrak X}_+,{\mathcal V})$ as a $k_{\mathcal E}^+[\varphi,\Gamma]$-module and that $H_0(\overline{\mathfrak X}_+,{\mathcal V})$ is a torsion $k_{{\mathcal E}}^+$-module.

{\it Step 1:} We claim that for $1\le i\le q$, for $v\in{\mathcal
  R}_i$ and for $0\le j\le t_v-1$ there are $n=n(i,v,j)\in{\mathbb Z}_{\ge0}$ and
$\varrho=\varrho(i,v,j)\in k^{\times}$ with
\begin{gather}t^{n}\varphi g_{w_{j-1}(i,v)}=\varrho g_{w_{j}(i,v)}\mbox{ modulo }N_{i-1}.\label{modnin}\end{gather}To prove this we begin by computing$$\varphi
g_{w_{j-1}(i,v)}=s_dug_{w_{j-1}(i,v)}=s_dT_{u^{-1}}(g_{w_{j-1}(i,v)})=\epsilon_{w_{j}(i,v)}s_dg_{{w_{j-1}(i,v)}\overline{u}^{-1}}.$$Now
 $s_dg_{{w_{j-1}(i,v)}\overline{u}^{-1}}$ is related to $n_{s_d}^{-1}g_{{w_{j-1}(i,v)}\overline{u}^{-1}}$ by a factor $\pm 1$, and we
 have ${w_{j}(i,v)}={{w_{j-1}(i,v)}_+^{W_{i-1}}}$. Thus
we need to show that for $w=w_{j-1}(i,v)$ there are $n\in{\mathbb Z}_{\ge0}$ and $\varrho'\in
k^{\times}$ with$$t^{n}n_{s_d}^{-1}g_{w\overline{u}^{-1}}=\varrho' g_{w_+^{W_{i-1}}}\quad\mbox{
  modulo }N_{i-1}.$$To do this, we turn back to the definition of $w_+^{W_{i-1}}$: 

In the cases where $\sigma(w\overline{u}^{-1})=0$ or
$\sigma(w\overline{u}^{-1}s_d)=0$ this follows from formula (\ref{ly5}). 

In the cases where [$\sigma(w\overline{u}^{-1})=-1$ or
$\sigma(w\overline{u}^{-1}s_d)=1$] and $w\overline{u}^{-1}s_d\overline{u}\notin W_{i-1}$ this follows from formula (\ref{ly2})
resp. (\ref{ly7}). 

In the cases where [$\sigma(w\overline{u}^{-1})=1$ or
$\sigma(w\overline{u}^{-1}s_d)=-1$] and $w\overline{u}^{-1}s_d\overline{u}\in W_{i-1}$ this follows from formula (\ref{ly4})
resp. (\ref{ly9}) since $\sum_{i\ge0}k.t^is_d
g_{w\overline{u}^{-1}s_d}=\sum_{i\ge0}k.t^i\varphi
g_{w\overline{u}^{-1}s_d\overline{u}}$ is contained in $N_{i-1}$.

In the cases where [$\sigma(w\overline{u}^{-1})=1$ or
$\sigma(w\overline{u}^{-1}s_d)=-1$] and $w\overline{u}^{-1}s_d\overline{u}\notin W_{i-1}$ this follows from formula (\ref{ly1})
resp. (\ref{ly6}). 

In the cases where [$\sigma(w\overline{u}^{-1})=-1$ or
$\sigma(w\overline{u}^{-1}s_d)=1$] and $w\overline{u}^{-1}s_d\overline{u}\in W_{i-1}$ this
follows from formula (\ref{ly3})
resp. (\ref{ly8}) since $\sum_{i\ge0}k.t^is_d
g_{w\overline{u}^{-1}s_d}=\sum_{i\ge0}k.t^i\varphi
g_{w\overline{u}^{-1}s_d\overline{u}}$ is contained in $N_{i-1}$. 

{\it Step 2:} We claim that for any $i\ge1$ the
  subquotient $N_i/N_{i-1}$ is a direct sum, indexed by $v\in{\mathcal R}_i$, of standard cyclic $k_{\mathcal
    E}^+[\varphi,\Gamma]$-modules $N_{i,v}/N_{i-1}$ of perimeter $t_v$, and that $N_i$ is a direct
  summand of $H_0(\overline{\mathfrak X}_+,{\mathcal
    V})$ as a $k_{\mathcal
    E}^+$-module. We proceed by induction on $i$.

By induction hypothesis, $N_{i-1}$ is a direct summand of $H_0(\overline{\mathfrak X}_+,{\mathcal
    V})$ as a $k_{\mathcal
    E}^+$-module. Therefore, and since
  $M={\rm ker}(t|_{H_0(\overline{\mathfrak X}_+,{\mathcal
    V})})$ by Theorem \ref{finite}, the classes of the $g_{w_j(i,v)}$ (for
$v\in{\mathcal R}_i$ and $0\le j\le t_v-1$) form a $k$-basis of ${\rm
  ker}(t|_{N_{i}/N_{i-1}})$. Together with step 1 we deduce our claim
concerning $N_i/N_{i-1}$ and the $N_{i,v}/N_{i-1}$. Next, if $N_{i,v}/N_{i-1}$ is $t$-divisible then
  with $N_{i-1}$ also $N_{i,v}$ is a direct summand of $H_0(\overline{\mathfrak X}_+,{\mathcal
    V})$ as a $k_{\mathcal E}^+$-module. If however $N_{i,v}/N_{i-1}$ is not
  $t$-divisible then it must be entirely contained in ${\rm
    ker}(t|_{H_0(\overline{\mathfrak X}_+,{\mathcal
    V})/N_{i-1}})$, cf. Proposition \ref{stnle}. This means
  that the numbers $n(i,v,j)$ appearing in step 1 are zero for all $0\le j\le t_v-1$. In view of the formulae in Lemma \ref{fmlueber}, the case by case distinction in step 1 in the passage from $w_j(i,v)$ to $w_{j+1}(i,v)=w_j(i,v)_+^{W_{i-1}}$ then reveals that for any $0\le j\le t_v-1$ we have $w_{j}(i,v)\overline{u}^{-1}s_d\overline{u}\in W_{i-1}$ or $\sigma(w_j(i,v)\overline{u}^{-1})=0$ or $\sigma(w_j(i,v)\overline{u}^{-1}s_d)=0$. It follows that for all $w\in
  W-W_i$ and all $l\ge i$ we have
  $w_+^{W_l}=w_+^{W_{l}-\{w_0(i,v),\ldots,w_{t_v-1}(i,v)\}}$. Therefore, and as
  $N_{i,v}/N_{i-1}\subset {\rm
    ker}(t|_{H_0(\overline{\mathfrak X}_+,{\mathcal
    V})/N_{i-1}})$, step 1 shows that the $k_{\mathcal
    E}^+[\varphi]$-submodule of $H_0(\overline{\mathfrak X}_+,{\mathcal
    V})/N_{i-1}$ generated by the $g_w$ for $w\in W-\{w_0(i,v),\ldots,w_{t_v-1}(i,v)\}$ is a $k_{\mathcal
    E}^+[\varphi]$-module complement of $N_{i,v}/N_{i-1}$. In particular,
  $N_{i,v}/N_{i-1}$ is a direct summand of $H_0(\overline{\mathfrak X}_+,{\mathcal
    V})/N_{i-1}$ as a $k_{\mathcal
    E}^+$-module, hence $N_{i,v}$ is a direct
  summand of $H_0(\overline{\mathfrak X}_+,{\mathcal
    V})$.\hfill$\Box$\\ 

{\bf Remark:} We leave to the reader the necessary modifications of the proof of Theorem
\ref{newstafi} in the general case, i.e. where we allow $v\in
W^{s_d}$ with $\sigma(v)=-1$ and $\theta(vh_{s_d}v^{-1})={\bf 1}$. Following
Lemma \ref{sl2fpind} one first must refine the formulae (\ref{ly6}), (\ref{ly7}), (\ref{ly8}), (\ref{ly9}): If $\sigma(v)=-1$ let us put
  ${\bf F}_v=1$ if $\theta(vh_{s_d}v^{-1})={\bf 1}$ and ${\bf F}_v=0$ if
  $\theta(vh_{s_d}v^{-1})\ne{\bf 1}$. Then formulae (\ref{ly6}), (\ref{ly7}),
  (\ref{ly8}), (\ref{ly9}) remain true for general $v\in
W^{s_d}$ with $\sigma(v)=-1$ if each occurence of $\kappa_{vs_d}g_{v}$
  is replaced by $\kappa_{vs_d}g_{v}-{\bf F}_vg_{vs_d}$. In addition, by formulae
  (\ref{homcae}), (\ref{homcaf}), if $\sigma(v)=-1$ and $\theta(vh_{s_d}v^{-1})={\bf 1}$
  we have the formulae $n_{s_d}^{-1}g_v=g_v$ and $t^{p-1}n_{s_d}^{-1}g_{vs_d}+g_{vs_d}\in\sum_{i\ge0}k.t^i{s_d}g_v$.\\

{\bf Remark:} Apparently, the above proof of Theorem \ref{newstafi} allows us to
precisely compute the parameters $n(i,v,j)\in{\mathbb Z}_{\ge0}$ and
$\varrho(i,v,j)\in k^{\times}$, hence the standard cyclic $k_{\mathcal
    E}^+[\varphi,\Gamma]$-modules $N_{i,v}/N_{i-1}$, hence their associated ${\rm
Gal}_{{\mathbb Q}_p}$-representations. In the following Theorem \ref{descmodprincser} we do this in the modular principal
  series case $\sigma={\bf 1}$.\\ 

Let ${\mathcal R}\subset W$ and $w_j(i,v)\in W$ be as in Theorem \ref{sigma1fil}.

\begin{satz}\label{descmodprincser} Suppose that $\sigma={\bf 1}$. Then $H_0(\overline{\mathfrak X}_+,{\mathcal V})$ admits a $k_{\mathcal E}^+[\varphi,\Gamma]$-module filtration $0=N_0\subset N_1\subset\ldots\subset N_{d+1}=H_0(\overline{\mathfrak X}_+,{\mathcal V})$ together with direct sum decompositions\begin{gather}{N_i}/{N_{i-1}}=\bigoplus_{v\in {\mathcal R}}{N_{i,v}}/{N_{i-1}}\label{dirmodpro}\end{gather}such that $N_{i,v}/N_{i-1}$ for $1\le i\le d+1$ and $v\in {\mathcal R}$ is a standard cyclic $k_{\mathcal E}^+[\varphi,\Gamma]$-module of perimeter $d$. More precisely we have: For fixed $i$ and $v\in {\mathcal R}$ put $w_j=w_j(i,v)$ for $0\le j\le d$; then there is a $k$-basis $e_0,\ldots,e_{d-1}$ of ${\rm ker}(t|_{N_{i,v}/N_{i-1}})$ such that, setting $e_d=e_0$, for $0\le j\le d-1$ we have\begin{align}t^{p-1}\varphi e_j&=\epsilon_{w_j}e_{j+1}&\mbox{ if }w_{j+1}\notin W^{s_d},\label{phi1modpri}\\\varphi e_j&=\epsilon_{w_j}e_{j+1}&\mbox{ if }w_{j+1}\in W^{s_d},\label{phi2modpri}\\\gamma(x)e_j&=\theta(w_j\tau(x)w_j^{-1})e_j&\mbox{ for all }x\in{\mathbb F}_p^{\times}.\label{gammodpri}\end{align}The distribution between the cases (\ref{phi1modpri}) and (\ref{phi2modpri}) is given by formula (\ref{carw}). In particular, if $1\le i\le d$ then the ${\rm
Gal}_{{\mathbb Q}_p}$-representation associated with $N_{i,v}/N_{i-1}$ has dimension $d$, whereas if $i=d+1$ it has dimension $0$.  
\end{satz}

{\sc Proof:} This follows from Theorem \ref{sigma1fil} and Theorem \ref{newstafi}. To see formula (\ref{gammodpri}) recall that
$\gamma(x)$ acts as $\tau(x)$, and that $\tau(x)$ acts as the Hecke operator
$T_{\tau(x)^{-1}}$. For the formulae (\ref{phi1modpri}) and (\ref{phi2modpri}) we must inspect the above proof of formula (\ref{modnin}); namely, it is enough to prove that for all $w\in W$ we
have \begin{align}t^{p-1}\varphi
  g_w&=\epsilon_{w}g_{w\overline{u}^{-1}s_d}&\mbox{ if
  }w\overline{u}^{-1}\in
  W^{s_d},\label{slphi1modpri}\\\varphi
  g_w-\epsilon_{w}g_{w\overline{u}^{-1}s_d}&\in \sum_{n\ge0}k.t^n\varphi
  g_{w\overline{u}^{-1}s_d\overline{u}}&\mbox{ if }w\overline{u}^{-1}\notin
  W^{s_d}\label{slphi2modpri}\end{align}in $H_0(\overline{\mathfrak X}_+,{\mathcal V})$. As before we begin with$$\varphi
g_w=s_dug_w=s_dT_{u^{-1}}(g_w)=\epsilon_ws_dg_{w\overline{u}^{-1}}.$$Suppose first that $w\overline{u}^{-1}\in W^{s_d}$. Define $1\le r\le p-1$ by
the
condition $$\theta(w\overline{u}^{-1}h_{s_d}^{-1}(x)\overline{u}w^{-1})=x^r$$
for all $x\in{\mathbb F}_p^{\times}$. Then, since
$\sigma(w\overline{u}^{-1})=1$, we find that, as an ${\mathcal H}(\overline{\mathcal
  S},\overline{\mathcal U})_k$-module, $\langle
g_{w\overline{u}^{-1}s_d}, g_{w\overline{u}^{-1}}\rangle$ is isomorphic with
$M_r$ as considered in Lemma \ref{sl2fpind} --- take $e=g_{w\overline{u}^{-1}}$ there. We compute \begin{align}t^{p-1}s_d g_{w\overline{u}^{-1}}&=t^{p-1}s_dn_{s_d}n_{s_d}^{-1}g_{w\overline{u}^{-1}}\notag\\{}&\stackrel{(i)}{=}s_dn_{s_d}t^{p-1}n_{s_d}^{-1}g_{w\overline{u}^{-1}}\notag\\{}&\stackrel{(ii)}{=}s_dn_{s_d}T_{n_{s_d}}g_{w\overline{u}^{-1}}\notag\\{}&=s_dn_{s_d}T_{s_dn_{s_d}}T_{{s_d}}g_{w\overline{u}^{-1}}\notag\\{}&=g_{w\overline{u}^{-1}s_d}\notag\end{align}where $(i)$ uses $t^{p-1}s_dn_{s_d}=s_dn_{s_d}t^{p-1}$ while $(ii)$ uses formula (\ref{homcaa}).

Now suppose that $w\overline{u}^{-1}\notin W^{s_d}$. Define $1\le r\le p-1$ by
$$\theta(w\overline{u}^{-1}s_dh_{s_d}^{-1}(x)s_d\overline{u}w^{-1})=x^r$$ for
all $x\in{\mathbb F}_p^{\times}$. Then, since
$\sigma(w\overline{u}^{-1}s_d)=1$, we find that $\langle
g_{w\overline{u}^{-1}s_d}, g_{w\overline{u}^{-1}}\rangle$ is isomorphic with
$M_r$ as considered in Lemma \ref{sl2fpind} --- this time take
$e=g_{w\overline{u}^{-1}s_d}$ there. We
compute\begin{align}s_dg_{w\overline{u}}&=s_dT_{s_d}(g_{w\overline{u}^{-1}s_d})\notag\\{}&=s_dT_{n_{s_d}^{-1}{s_d}}T_{n_{s_d}}(g_{w\overline{u}^{-1}s_d})\notag\\{}&=n_{s_d}T_{n_{s_d}}(g_{w\overline{u}^{-1}s_d})\notag\\{}&\equiv
  g_{w\overline{u}^{-1}s_d}\notag\end{align}where
the last congruence uses formula (\ref{homcac}) and is to be understood as an identity modulo $$\sum_{n\ge0}k.t^nn_{s_d}g_{w\overline{u}^{-1}s_d}=\sum_{n\ge0}k.t^ns_d g_{w\overline{u}^{-1}s_d}=\sum_{n\ge0}k.t^n\varphi  g_{w\overline{u}^{-1}s_d\overline{u}}.$$\hfill$\Box$\\

 {\bf Further examples:} (a) The easiest case is the generic case $\sigma={\bf 0}$, i.e. $\sigma(w)=0$ for all $w\in
  W^{s_d}$. ('Most' reduced standard ${\mathcal H}(G,I_0)_k$-modules arising by reduction from a locally unitary principal series representation are of this sort, cf. \cite{wty}.) To describe it, fix a set of representatives ${\mathcal Q}$ in $W$ for the right cosets of the cyclic group
$\overline{u}^{\mathbb Z}=\{\overline{u}^j\,|\,j=0,\ldots,
d\}$, i.e. such that $W=\coprod_{v\in
  {\mathcal Q}}v\overline{u}^{\mathbb Z}$. For
example, one may take ${\mathcal Q}=\langle
s_1,\ldots,s_{d-1}\rangle$. Then $H_0(\overline{\mathfrak X}_+,{\mathcal V})$ admits a direct sum
  decomposition $$H_0(\overline{\mathfrak X}_+,{\mathcal V})=\oplus_{v\in
      {\mathcal Q}}N_{1,v}$$ such that $N_{1,v}$ for each $v\in
    {\mathcal Q}$ is standard cyclic of perimeter $d+1$. More precisely, it can be
    described as follows. The set with $d+1$ elements $\{g_{v}, g_{v\overline{u}},\ldots,
 g_{v\overline{u}^{d}}\}$ is a $k$-basis of ${\rm ker}(t|_{N_{1,v}})$. It
 generates $N_{1,v}$ as a $k_{\mathcal E}^+[\varphi,\Gamma]$-module. In $H_0(\overline{\mathfrak
    X}_+,{\mathcal V})$ we have the identities (the proof proceeds as in Theorem \ref{descmodprincser}, but is easier)\begin{align}t^{\lceil k(v\overline{u}^{j-1})\rceil}\varphi g_{v\overline{u}^j}&=\epsilon_{v\overline{u}^{j}}\kappa_{v\overline{u}^{j-1}}\lceil
  k(v\overline{u}^{j-1})\rceil ! g_{v\overline{u}^{j-1}}&\mbox{ if }v\overline{u}^{j-1}\notin W^{s_d},\notag\\t^{\lfloor k(v\overline{u}^{j-1})\rfloor}\varphi g_{v\overline{u}^j}&=\epsilon_{v\overline{u}^{j}}\kappa_{v\overline{u}^{j-1}}\lfloor
  k(v\overline{u}^{j-1})\rfloor ! g_{v\overline{u}^{j-1}}&\mbox{ if
  }v\overline{u}^{j-1}\in
  W^{s_d},\notag\\\gamma(a)g_{v\overline{u}^j}&=\theta(v\overline{u}^j\tau(a)\overline{u}^{-j}v^{-1})g_{v\overline{u}^j}&\mbox{ for all }a\in{\mathbb Z}_p.\notag\end{align}

(b) Assume $d=2$. Theorem \ref{newstafi} and our examples in section \ref{filtrsec} show that, for suitable choices of $\sigma$, the $k_{\mathcal E}^+[\varphi,\Gamma]$-module $H_0(\overline{\mathfrak X}_+,{\mathcal V})$ admits a two-step filtration with standard cyclic subquotients of perimeters $2$ and $4$ (resp. $4$ and $2$). At least for generic $\theta$, the corresponding ${\rm
Gal}_{{\mathbb Q}_p}$-representation admits a two-step filtration with irreducible subquotients of dimensions $2$ and $4$ (resp. $4$ and $2$).\\

\begin{flushleft} \textsc{Humboldt-Universit\"at zu Berlin\\Institut f\"ur Mathematik\\Rudower Chaussee 25\\12489 Berlin, Germany}\\ \textit{E-mail address}:
gkloenne@math.hu-berlin.de \end{flushleft} 
\begin{thebibliography}{abcdefgh} 

\bibitem{baho}{\it P. N. Balister and S. Howson}, Note on Nakayama's Lemma for compact $\Lambda$-modules, Asian J. Math. {\bf 1}, no. 2, 224--229 (1997)

\bibitem{berger}{\it L. Berger}, On some modular representations of the Borel subgroup of ${\rm GL}_2({\mathbb Q}_p)$.
    Compositio Math. {\bf 146} (2010), no. 1, 58--80.

\bibitem{calu}{\it R. W. Carter, G. Lusztig}, Modular representations of
  finite groups of Lie type, Proc. London Math. Soc. (3) {\bf 32} (1976),
  347--348

\bibitem{col}{\it P. Colmez}, $(\varphi,\Gamma)$-modules et repr\'{é}sentations du mirabolique de ${\rm GL}_2({\mathbb Q}_p)$, Ast\'{e}risque {\bf 330} (2010), 61-153.

\bibitem{colhaupt}{\it P. Colmez}, Repr\'{e}sentations de ${\rm GL}_2({\mathbb Q}_p)$ et $(\varphi,\Gamma)$-modules, Ast\'{e}risque {\bf 330} (2010), 281-509

\bibitem{emert}{\it M. Emerton}, On a class of coherent rings, with
  applications to the smooth representation theory of ${\rm GL}_2(Q_p)$ in
  characteristic $p$, preprint 2008

\bibitem{fon}{\it J. M. Fontaine}, Repr\'{e}sentations $p$-adiques des corps locaux. I. The Grothendieck Festschrift, Vol. II,
 249--309, Progr. Math. {\bf 87}, Birkh\"auser Boston, Boston, MA, 1990. 

\bibitem{jalg}{\it E. Grosse-Kl\"onne}, $p$-torsion coefficient systems for
  ${\rm SL}_2({\mathbb Q}_p)$ and ${\rm GL}_2({\mathbb Q}_p)$, Journal of
  Algebra {\bf 342}, 97--110 (2011)

\bibitem{wty}{\it E. Grosse-Kl\"onne}, Locally unitary principal series representations of ${\rm GL}_{d+1}(F)$, M\"unster Journal of Mathematics, special issue dedicated to Peter Schneider, to appear 

\bibitem{ihg2}{\it E. Grosse-Kl\"onne}, From pro-$p$ Iwahori-Hecke modules to
    $(\varphi,\Gamma)$-modules II, preprint

\bibitem{oll}{\it R. Ollivier}, Parabolic induction and Hecke
modules in characteristic $p$ for $p$-adic ${\rm GL}_n$, Algebra and Number
Theory {\bf 4}, no.6 (2010), 701–-742.

\bibitem{os}{\it R. Ollivier, P. Schneider}, Pro-$p$ Iwahori-Hecke algebras are Gorenstein, preprint

\bibitem{ollsec}{\it R. Ollivier, V. S\'{e}cherre}, Modules universels de
  ${\rm GL}_3$ sur un corps $p$-adique en caracteristique $p$, preprint

\bibitem{vien}{\it M. Vienney}, Construction de $(\varphi,\Gamma)$-modules en caract\'{e}ristique $p$, Lyon 2012.

\bibitem{vigneras}{\it M.F. Vign\'{e}ras}, Pro-$p$-Iwahori Hecke ring and supersingular $\overline{\mathbb F}_p$-representations,  Math. Ann.  {\bf 331}  (2005),  no. 3, 523--556 and Math. Ann. {\bf  333}  (2005),  no. 3, 699--701. 

\end{thebibliography}
\end{document}